\documentclass[11pt]{article} 

\usepackage{jmlr2e}
\usepackage{soul}
\usepackage{color}
\usepackage[utf8]{inputenc} 
\usepackage[T1]{fontenc}    
\usepackage{mathtools}
\usepackage{algpseudocode,algorithm2e}
\usepackage{enumitem}
\usepackage{url}
\usepackage{natbib}
\usepackage[dvipsnames]{xcolor}
\usepackage{xspace}

\input{macromine01.tex}

\usepackage{floatrow}
\newfloatcommand{capbtabbox}{table}[][\FBwidth]

\usepackage{soul}

\usepackage{booktabs}
\usepackage{multirow}
\usepackage{times}
\usepackage{marginnote}
\usepackage{comment}
\usepackage[linecolor=blue!60!,backgroundcolor=blue!10!,textwidth=5.5cm]{todonotes}

\usepackage{hyperref}
\usepackage[nameinlink,noabbrev,capitalize]{cleveref}
\hypersetup{colorlinks,
            linkcolor=blue,
            citecolor=blue,
            urlcolor=magenta,
            linktocpage,
            plainpages=false}

\title{Online non-parametric regression with kernels}

\author{
    Oleksandr Zadorozhnyi and Pierre Gaillard \\
    \bfseries{and Sebastien Gerchinovitz and Alessandro Rudi}}
\date{\today}

\newcommand{\algo}{KAAR\xspace}

\ShortHeadings{Online nonparametric regression with kernels}{}
\firstpageno{1}

\begin{document}
	
	\title{Online nonparametric regression with Sobolev kernels}
	
	\author{\name Oleksandr Zadorozhnyi \email oleksandr.zadorozhnyi@uni-potsdam.de\\
		\addr Insitute of Mathematics\\
		University of Potsdam\\
		14471 Potsdam, Germany
		\AND
		\name Pierre Gaillard \email pierre.gaillard@inria.fr \\
		\addr  Centre de Recherche INRIA de Paris \\
		2 rue Simone Iff \\
		75012 Paris France 
	\AND 
		\name S\'ebastien Gerchinovitz \email sebastien.gerchinovitz@irt-saintexupery.com \\
		\addr  IRT Saint Exup\'ery \& Institut de Math\'ematiques de Toulouse\\
		3 rue Tarfaya \\
		31400 Toulouse
	\AND 
		\name Alessandro Rudi \email alessandro.rudi@inria.fr \\
		\addr Centre de Recherche INRIA de Paris \\
		2 rue Simone Iff  \\
	  75012 Paris France 
	}
	\editor{ Unknown Editor}

\maketitle

\begin{abstract}
In this work we investigate the variation of the online kernelized ridge regression algorithm in the setting of $d-$dimensional
	adversarial nonparametric regression. We derive the regret upper
	bounds on the classes of Sobolev spaces $W_{p}^{\beta}(\cX)$, $p\geq
	2, \beta>\frac{d}{p}$. The upper bounds are supported by the minimax regret
	analysis, which reveals that in the cases $\beta> \frac{d}{2}$ or
	$p=\infty$ these rates are (essentially) optimal. Finally, we compare the performance of the kernelized ridge regression forecaster to the known non-parametric forecasters in terms of the regret rates and their computational complexity as well as to the excess risk rates in the setting of statistical (i.i.d.) nonparametric regression.
\end{abstract}


\section{Introduction}
\label{sec:intro}
We consider the online least-squares regression framework \citep{Bianchi:06} as a game between the environment and the learner where the task is to sequentially predict the environment's output $y_{t}$ given the current input  $x_{t}$ and the observed history $\{ (x_{i},y_{i})\}_{i =1}^{t-1}$. Specifically, let $\cX \subset \mbr^d$ be an input space, $\cY  \subset \mbr$ a label space, and $\hat{\cY}\subset \mbr$ a target space. Before the game starts, the environment secretly produces a sequence of input-output pairs  $(x_1,y_1), (x_2,y_2),\dots$ in  $\cX \times \cY$ over some (possibly infinite) time horizon.

At each round $t\geq 1$, the environment first reveals an input $x_{t} \in \cX$; the learner predicts $\smash{\hat{y}_{t} \in \hat{\cY}}$ based on past information $(x_1,y_1),\dots,(x_{t-1},y_{t-1}) \in \cX \times \cY$ and on the current input $x_{t}$, which is considered the estimate of the true label $y_t \in \cY$. The true label $y_{t}$ is then revealed, the learner suffers the squared loss $({y_{t}}-\hat{y}_{t})^{2}$ and round $t+1$ starts. The problem is to design an algorithm which minimizes  the learner's cumulative regret
\begin{align}
\label{eq:regret_f}
R_{n}\paren{\cF} \bydef \sup_{f \in \cF}R_{n}\paren{f}\,, \qquad \text{where} \quad  R_n(f) \bydef   \sum_{t=1}^{n}\paren{y_{t} - \hat{y}_{t}}^{2} -  \sum_{t=1}^{n}\paren{y_{t} -f\paren{x_{t}} }^{2} \,,
\end{align}
over $n \geq 1$ rounds with respect to the best fixed prediction rule from some reference functional class $\cF \subset \mbr^{\cX}$. 

Unlike in the standard statistical learning framework where the data stream is assumed to be generated from some underlying stochastic process, usually with an independent noise component in the setting of adversarial online learning no stochastic assumption on the nature of the datasample $\{x_{s},y_{s}\}_{s=1}^{T}$ is posed. The problem of online learning with arbitrary (adversarial) data goes back to the work of \cite{Foster:91}. A lot of theoretical research has been done since then for parametric models \citep[see example][]{Azoury:01,Bianchi:99,Vovk:97}. 
The amount of data and the complexity of current machine learning problems have led the community to explore the more general problem of online-learning with methods based on nonparametric decision rules and with the reference classes being bounded functional sets of continuous functions (see ex. \cite{Vovk06Metric}, \cite{Rakhlin:14}). 
Much effort has been devoted to the regret analysis with respect to functional classes that include Sobolev spaces \citep{Rakhlin:14,Rakhlin:14a,Vovk06Metric,Vovk:07}. 
Surprisingly, only a few explicit algorithms have been designed to address the regression problem \citep{Vovk06Metric,Vovk:06,Vovk:07,Gaillard:15}. While having optimal (or close to optimal) regret rates, they have the disadvantage of either being computationally intractable or of providing suboptimal regret upper bounds (see Table for computational complexities of some known algorithms \ref{tab:main_tabl}). 
For more details on previous work, we refer the reader to Section~\ref{sect:comparison}.

In this work we consider the framework of online adversarial regression over the benchmark classes $\cF$ being the bounded balls of continuous representatives in  Sobolev spaces (\citep[see e.g.,][]{Adams:03}) $W_{p}^{\beta}\paren{\cX}$ , $p \geq 2$ and $\beta \geq \frac{d}{p}$. 

The problem is of interest since, to date, the computationally efficient algorithm (\citep[see e.g.,][]{jezequel2019efficient})  which achieves the optimal regret rates is provided in the case when the underlying Sobolev ball is included in a Sobolev RKHS. The latter corresponds to the case when $\beta > \frac{d}{2}$ and $p=2$.  

\paragraph{Overview of the main results and outline of the paper} The aim of this paper is to provide a deeper analysis of the regret achieved a version of online kernel ridge regression algorithm, Kernel Aggregating Algorithm Regression, (\algo, see \cite{Gammerman2004}). In particular the key contribution is the analysis of the robustness of the \algo which returns an element in RKHS while competing against a function from Sobolev class which does not belong to a RKHS. We notice that (on the contrary to many known nonparametric schemes, see for example \cite{Rakhlin:14}, \cite{Vovk:07}) this algorithm is computationally tractable. Comparison of the performance of \algo to the known procedures (both in regret rates and computational efficiency) is summarized in Table~\ref{tab:main_tabl}. Furthermore, we also prove lower bounds for minimax regret (which is defined as the infimum over all admissible strategies of a supremum of all data-sequences), which assures that \algo reaches optimal or close to optimal (up to an arbitrary small polynomial factor in the number of rounds) regret rates on bounded balls of Sobolev spaces $W_p^\beta(\cX)$ with $p \geq 2$ and $\beta > \frac{d}{2}$ or when $p=\infty$.

\begin{table}[t!]
\begin{center}
\renewcommand{\arraystretch}{1.3}
\resizebox{\textwidth}{!}{%
\begin{tabular}{@{}lllcllclllcll@{}}
\toprule
{} & \multicolumn{2}{c}{ {{ \textbf{\algo} \eqref{eq:kern_ridge}} }} & \phantom{ab}& \multicolumn{2}{c}{ \cite{Rakhlin:14}} &
\phantom{ab} & \multicolumn{3}{c}{ \cite{Gaillard:15}}&
\phantom{ab} & \multicolumn{2}{c}{ EWA by \cite{Vovk06Metric}} \\
\cmidrule{2-3} \cmidrule{5-6} \cmidrule{8-10}  \cmidrule{12-13}
& Regret\footnotemark& Cost  & & Regret & Cost  & & Regret & Cost & Cost ($d=1, p=\infty$)\footnotemark& & Regret & Cost \\ 
\midrule
$\beta > \frac{d}{2}$ & \color{blue}{\textbf{$n^{1-\frac{2\beta}{2\beta+d} + \epsilon}$}} & {\textbf{$n^3+ dn^2$}} & & $n^{1-\frac{2\beta}{2\beta+d}}$  & Non constructive   & & $n^{1-\frac{2\beta}{2\beta+d}}$  &       $\exp\paren{n}$ & \text{poly(n)} & & $n^{1 - \frac{\beta}{\beta +d}}$  &       $\exp\paren{n} +nd$          \\
$\frac{d}{p}< \beta \leq \frac{d}{2}$ &  \color{blue}{\textbf{$n^{1-\frac{\beta}{d} \frac{p - d/\beta}{p-2} + \epsilon}$}}  & {\textbf{$n^{3}+dn^2$}} && $n^{1-\frac{\beta}{d}}$    & Non constructive   && $n^{1-\frac{\beta}{d}}$               &   $\exp\paren{n}$ & $n^{ \lceil \beta\rceil \big(\frac{5\beta+2}{2\beta +1}\big)}$ & & $n^{1- \frac{\beta}{\beta +d}}$  &       $\exp\paren{n}+nd$            \\ 
$p=\infty, \beta \leq \frac{d}{2} $ &  \color{blue}{$n^{1-\frac{\beta}{d} + \epsilon}$}                    & {$n^{3} + dn^2$}             &&  $n^{1-\frac{\beta}{d}}$                        &   Non constructive && $n^{1-\frac{\beta}{d}}$               & $\exp\paren{n}$ & $n^{ \lceil \beta\rceil \big(\frac{5\beta+2}{2\beta +1}\big)}$&  & $n^{1-\frac{\beta}{\beta +d}}$  &       $\exp\paren{n} +nd$\\
\bottomrule
\end{tabular}
}
\end{center}

\caption{Regret rates and time complexity of \algo~\eqref{eq:kern_ridge} (new upper bounds from this paper are highlighted in blue) and the existing algorithms for online nonparametric regression. }
\label{tab:main_tabl}
\end{table}

More precisely, the result is threefold. On the one hand, our analysis recovers the classical result  for Sobolev spaces, i.e. when $\beta > d/2$ and $p \geq 2$.
In particular,  we show in Theorem~\ref{thm:reg_tar} that on the classes of continuous functions which belong to Sobolev RKHS of smoothness $\beta$, \algo (with properly chosen regularization parameter) achieves the optimal regret upper bound\footnotemark
\[
        R_{n}(\cF) \ \lesssim   \  n^{1- \frac{2 \beta}{2\beta+ d}} ~ \log n.
\]

\addtocounter{footnote}{-3}
\stepcounter{footnote}\footnotetext{In terms of its upper bound.}
\stepcounter{footnote}\footnotetext{\cite{Gaillard:15} only provide an efficient version of their algorithm for Sobolev spaces with $p=\infty$,  $d=1$ and $\beta \geq 1/2$. Their efficient algorithm can however be extended for any $\beta \in (0,1/2)$ with a polynomial time complexity.}
\stepcounter{footnote} \footnotetext{The notation $\lesssim$ denotes an approximate inequality which includes multiplicative constants which depend on $\cF$ and $\cX$.}
On the other hand, we consider the more challenging scenario when ${1}/{2} \geq \frac{\beta}{d} > 1/p$ which corresponds to the benchmark functional classes that cannot be embedded into a RKHS and that have smaller smoothness. We will refer to this case as the \textit{hard-learning scenario}. In Theorem~\ref{thm:higher_beta_p} we prove that in such a scenario with $\smash{\cF = B_{W_{p}^{\beta}\paren{\cX}}\paren{0,R}}$ the regret of \algo is upper-bounded by
\[
    R_{n}\paren{\cF} \ \lesssim\  n^{1- \frac{\beta}{d}\frac{p-\frac{d}{\beta}}{p-2}} ~\log n.
\]
In particular, when $p = \infty$, the regret upper bound is of order $O(n^{1-\frac{\beta}{d} +\epsilon}~ \log n)$.
The latter bound is proven to be essentially optimal (up to a constant $\epsilon$ that can be made arbitrary small) by the corresponding lower bound for minimax regret in Section~\ref{sec:lower_bound} for the lower bounds. 
Optimal regret upper bounds on the classes of bounded H\"older balls were previously derived with polynomial-time algorithms for $d~=~1$ \cite{Gaillard:15}. The case $d\geq 1$ and $\beta =1$ was also analyzed for Lipschitz and semi-Lipschitz losses in \cite{cesa2017algorithmic}. Notice that throughout the paper we do not consider the case of Sobolev spaces with $\frac{\beta}{d} \leq 1/p$. In the latter case the existence of continuous representatives for equivalence classes in $W_{p}^{\beta}\paren{\cX}$ is not guaranteed.
\begin{figure}[t!]
\begin{center}
  \includegraphics[width=.5\textwidth]{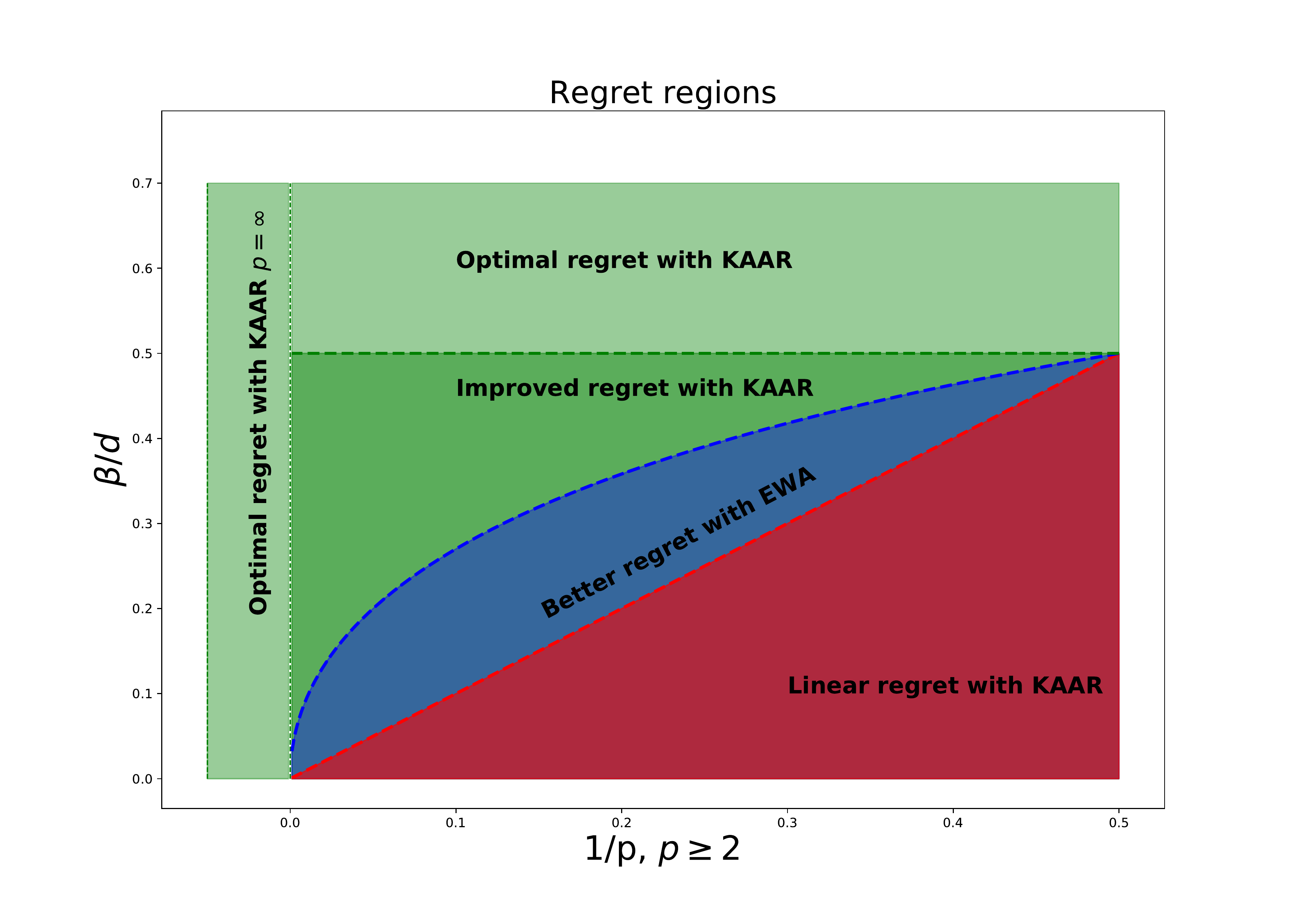}
  \hfill
  \includegraphics[width=.5\textwidth]{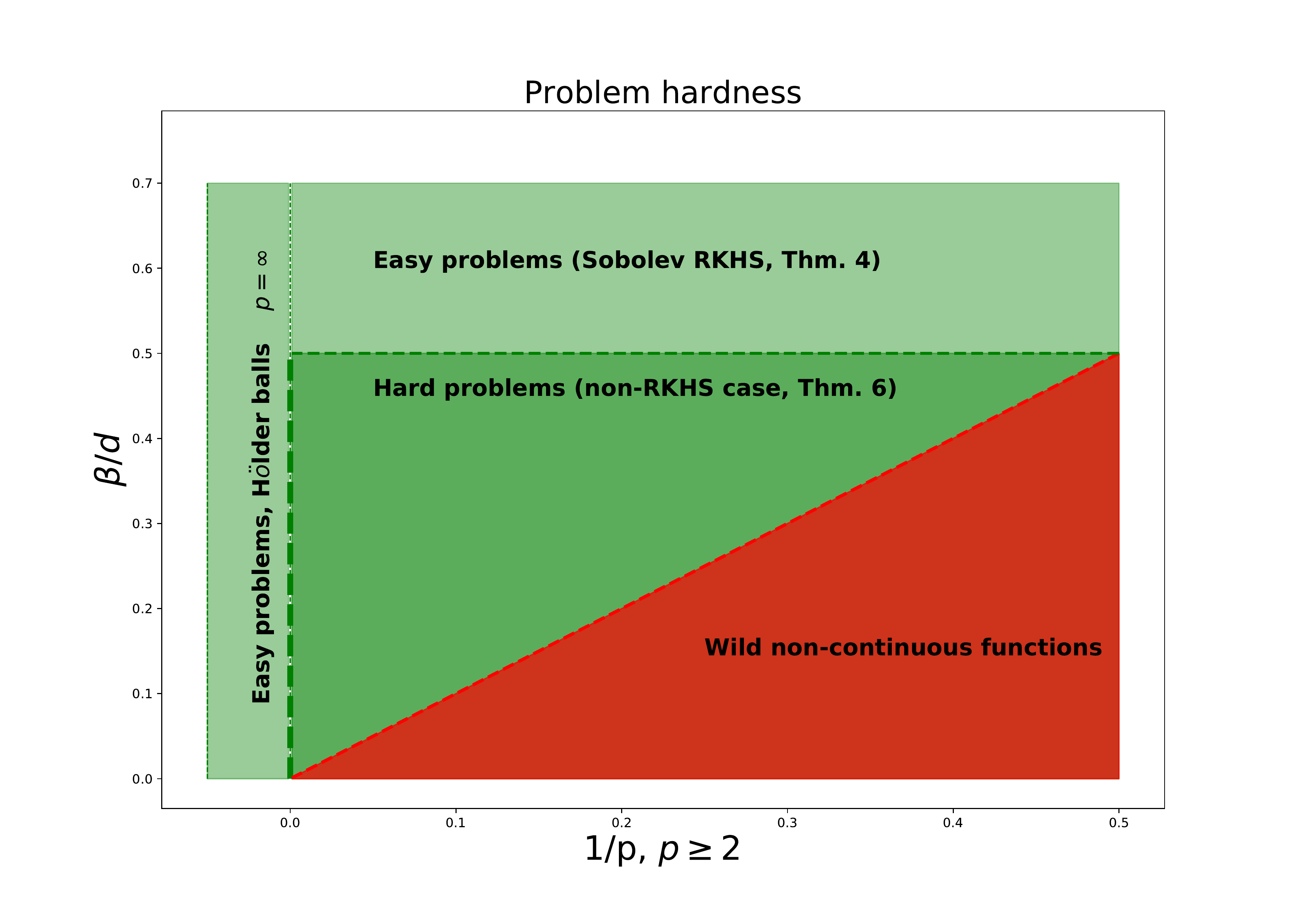}
\end{center}
\caption{(Left) Different regions in the $(\frac{1}{p}, \frac{\beta}{d})$-plane for which our new regret bound for \algo: [light green] is optimal (i.e., $\beta > d/2$ or $p=\infty$); [dark green] improves the bound of EWA by~\cite{Vovk:06}; [blue] is worse than the bound of EWA; [red] is linear in $n$ (i.e, $\beta \leq d/p$). (Right) Hardness of the problem in the $(\frac{1}{p}, \frac{\beta}{d})$ plane} 
\label{fig:opt-regions}
\end{figure}

In Figure~\ref{fig:opt-regions}, we plot the regions of the $(1/p, \frac{\beta}{d})$-plane corresponding to the different regret cases  where we obtain either the optimal rate or a suboptimal rate, but improved with respect to classical discretization algorithms in the nonparametric framework (see \cite{Vovk06Metric}). Note that the smaller $\beta/d$ and $p$ are, the harder the problem is. Additional graphs comparing the regret of \algo with the regret of EWA (see \cite{Vovk:97}) are available in Appendix~\ref{app:regret_rates_comp}. 

To complete the analysis of the \algo in the setting of online nonparametric regression over Sobolev spaces, 
we use the general results of~\cite{Rakhlin:14}, derive sharp bound on the fat-shattering dimension and establish corresponding lower bounds for the minimax regret.
More precisely, we prove that any admissible algorithm suffers at least regret of order $\smash{n^{1-2\beta/(2\beta + d)}}$ in the smooth case $\beta > d/2$, and $\smash{n^{1-\beta /d}}$ when $\beta \leq d/2$. In particular, this implies that \algo achieves optimal regret rates when $\beta > d/2$ or $p=\infty$. The regret analysis of \algo on the classes of bounded classes of continuous functions in Sobolev spaces $W_{p}^{\beta}\paren{\cX}$, $p \geq 2$ as well as lower bounds for the minimax regret for the classes of bounded balls in Sobolev spaces $W_{p}^{\beta}\paren{\cX}$ are summarized in Table~\ref{tab:bounds}. 

\begin{table}[!ht]
\begin{center}
\resizebox{.65\textwidth}{!}{
\begin{tabular}[t]{lcc}
\toprule
&Upper bound of \algo & Lower bound for minimax regret\\
\midrule
$\beta>\frac{d}{2}$& $n^{1 - \frac{2\beta}{2\beta + d} + \epsilon}\log(n)$ & $n^{1 - \frac{2\beta}{2\beta +d}}$\\
$\frac{d}{p} < \beta \leq \frac{d}{2}$& $n^{1-\frac{\beta}{d} \frac{p - d/\beta}{p-2}+\epsilon}\log\paren{n}$ &$n^{1-\frac{\beta}{d}}$\\
$p=\infty$, $\beta\leq \frac{d}{2}$ &$n^{1-\frac{\beta}{d} + \epsilon}\log\paren{n}$ &$n^{1-\frac{\beta}{d}}$\\
\bottomrule 
\end{tabular} }
\end{center}
  \caption{Regret upper bounds of \algo and the corresponding lower bound on the classes of bounded subsets of $W_{p}^{\beta}\paren{\cX}$, $\beta \in \mbr$, $p\geq 2$. Here $\epsilon>0$ is an arbitrary small number.}%
  \label{tab:bounds}
\end{table}

The outline of the rest of the paper is as follows. In Section~\ref{sec:not_backgrnd}, we fix the notation and recall the definition of Sobolev spaces, reproducing kernel Hilbert spaces (RKHS) and their data-sample based effective dimension. Furthermore, we describe \algo therein. In Section~\ref{sec:main_res}, we provide our regret upper bounds for \algo and in Section~\ref{sec:lower_bound} we present the corresponding lower bounds. 
Finally, in Section~\ref{sect:comparison}, we make more detailed comparisons with existing work both in the adversarial online regression setting and in the standard statistical framework with i.i.d. observations. We discuss the optimality of the rates, compare it to the excess risk analysis in the statistical case and comment on the aspect of computational complexity by showing that \algo is superior to the known nonparametric schemes in terms of runtime and storage complexities. All the proofs as well as technical details on Sobolev spaces and kernels are given in the Appendices. 

\section{Notation and background}
\label{sec:not_backgrnd}

\subsection{Kernels and effective dimension}
 \label{sec:kern_lrn}
We recall below some notations on reproducing kernel Hilbert Spaces (RKHS). An in-depth survey on this topic can be found in \cite{Smola:02} and \cite{Steinwart:08}.
A  Hilbert space of functions $\cH := \{f: \cX \mapsto \mbr \}$ equipped with an inner product $\inner{\cdot,\cdot}_{\cH}$ is called RKHS if for every $x \in \cX$ the evaluation functional $\delta_{x}\paren{f} := f\paren{x}$ is continuous in $f$.
Furthermore, we say that function $k\paren{\cdot,\cdot}: \cX \times \cX \mapsto \mbr$ over domain $\cX$ is a real-valued kernel if every kernel matrix $K_{n} \bydef \paren{k(x_{i},x_{j})}_{i,j=1}^{n}$ is positive semi-definite. It is known that value of kernel can be represented as an inner product in some Hilbert space $H$, namely $k\paren[1]{x,x^{'}} = \inner[1]{\phi\paren{x}, \phi\paren{x'}}_{H}$, where we call $H$ a feature space and $\phi\paren{\cdot}$ a feature map of kernel $k$.  Lastly we say that RKHS $\cH$ is generated by kernel $k\paren{\cdot,\cdot}$ if for every $x$ it holds: $k_{x} \bydef k\paren{x,\cdot} \in \cH$ and $f\paren{x}= \inner[1]{f,k_x}_{\cH}$ for every $f \in \cH$ and $x \in \cX$ (i.e. the so-called reproducing property holds). In this case we write $\cH_{k}$ to denote the RKHS generated by kernel $k\paren{\cdot,\cdot}$ and say that kernel $k\paren{\cdot,\cdot}$ is a reproducing kernel of $\cH_{k}$. 
Denote also $\lambda_{j}(K_{n})$ to be the $j-$th largest eigenvalue of the matrix $K_{n}$. We give below the definition of the effective dimension, which measures the complexity of the underlying RKHS based on a given data sample. It plays a key role in our regret analysis of \algo.

\begin{definition}[Effective dimension]
	\label{def:eff_dim}
	Let $k: \cX \times \cX \to \mbr$ be a kernel function, $\cD_{n}=\{ \bx_{i}\}_{i=1}^{n} \in \cX^n$ be a sequence of inputs and $\tau>0$.
	The effective dimension associated with the sample 
	$\cD_{n}$, the kernel $k$ and a scale parameter $\tau$ is defined as
	\begin{align}
	\label{eq:effect_dimension}
	d_{eff}^{n}\paren{\tau} := \tr \paren[1]{{\paren{K_{n}+\tau I}^{-1}K_{n}}} = \sum_{j=1}^{n} \frac{\lambda_{j}\paren{K_n}}{\lambda_{j}\paren{K_n} + \tau},
	\end{align}
	where $I : \mbr^n \mapsto \mbr^n$ is the identity matrix and $K_{n}$ is the kernel matrix associated to the kernel $k$.
\end{definition}

In statistical learning, it has been shown (\cite{Zhang:05},  \cite{Rudi:15}, and \cite{Blanchard:17}) that the effective dimension characterizes the generalization error of kernel-based algorithms. 
This is a decreasing function of the scale parameter $\tau$ and $\smash{d_{eff}^{n}(\tau) \to 0}$ when $\tau \to \infty$. 
On the other side, as $\tau \to 0$, it converges to the rank of $K_{n}$, which can be interpreted as the "physical" dimension of the points $(k_{x_i})_{1\leq i\leq n}$.
Additional definitions and related notations on kernels (which are used in the proofs) are given in Appendix~\ref{app:kernels}. 

\subsection{Sobolev Spaces} Let $\beta \in \mbn_{*}$, $2\leq p < \infty$ and $\cX \bydef [-1,1]^{d}$, where we use standard notation for $\mbn \bydef \{0,1,2,\ldots,\}$,$\mbn_{*} \bydef \{1,2\ldots,\}$. We denote by $L_{p}\paren{\cX}$ the space of equivalence classes of $p$-integrable functions with respect to the Lebesgue measure $\lambda$ on the Borel $\sigma-$algebra $\cB(\cX)$ and by $[f]_{\lambda}$ the $\lambda-$equivalence class to some function $f: \cX \mapsto \mbr$. We denote $|\gamma|_{1} \bydef \sum_{i=1}^{n}\abs{\gamma_{i}}$ for $\gamma \in \mbn^{d}$ and we write  $D^{\gamma}f$ for the multidimensional weak derivative (see section 5.2.1, page 242 in \cite{Evans:98}) of the function $f: \cX \mapsto \mbr$ of order $\gamma \in \mbn_*^{d}$.  We denote by $C^{m}\paren{\cX}$ the space of all $m$-times differentiable functions $f$ with multidimensional derivative $D^{\gamma}f$ (  $\abs{\gamma}_{1} \leq m$, $\gamma \in \mbn^{d}$) that are continuous on $\cX$ 
and let $C(\cX)$ denotes the standard space of continuous functions equipped with the norm $\|f\|_{C(\cX)} = \max_{x \in \cX}|f(x)|$ (we write it simply  $\|f\|$  when no confusion can arise). For the normed space $(\cG,\|\cdot\|)$ we use $B_{\cG}(x,R)$ and $\smash{\overline{B}_{\cG}(x,R)}$ to denote respectively the open and the closed ball of radius $R$ centered at the point $x$.

We recall that the Sobolev space (see Chapter 3 in \cite{Adams:03})  $\smash{W_{p}^{\beta}(\cX)}$ is the space of all equivalence  classes of functions  $ [f]_{\lambda} \in L_{p}\paren{\cX}$ such that 
$$
\norm{f}_{W^{\beta}_{p}\paren{\cX}} \bydef 
    \left\{ 
    \begin{array}{ll}
        \paren[1]{ \sum_{\abs{\gamma}_{1} \leq \beta}\norm{D^{\gamma}f}_{L_{p}\paren{\cX}}^{p}}^{\frac{1}{p}} & \text{if}\quad  p < \infty \\
         \sup_{\abs{\gamma}_{1}\leq \beta} \norm{D^{\gamma}f}_{L_{\infty}\paren{\cX}} &  \text{if}\quad  p = \infty  
    \end{array} 
    \right. 
$$
is finite. 

The notion of Sobolev spaces is then extended to the case of any real $\beta>0$ (see Appendix~\ref{app:sob_spaces} for the details) by means of the Gagliardo semi-norms. In the case $p=2$ it can be shown to be equivalent to the known approach of fractional Sobolev spaces defined via Fourier transform.

\paragraph{Sobolev Reproducing Kernel Hilbert Spaces.} 
\label{sec:sob_spaces} 
We recall here known results on embedding characteristics of fractional Sobolev spaces, which are essential in our analysis. Let  $s \in \mbr_{+}$ and consider the Sobolev space ${W_{2}^{s}\paren{\cX}}$ with $\cX \subset \mbr^d$. It is a separable Hilbert space (see Chapter 7 in \cite{Schaback:07}) with the inner product ${\langle{f,g}\rangle=\sum_{\norm{\gamma}_1 \leq s} \langle{D^{\gamma}f,D^{\gamma}g}\rangle_{L_2\paren{\cX}}}$. By Sobolev Embedding Theorem (see Theorem 7.34 in \cite{Adams:03} for the case $s \in \mbr_{+}$, $s> {d}/{2}$) we have that $\smash{W_{2}^{s}\paren{\cX} \hookrightarrow C\paren{\cX}}$. The latter embedding is to be understood in the sense that there exists $C_1>0$, such that each equivalence class has a unique element $f \in C\paren{\cX}$ such that ${\norm{f}_{C\paren{\cX}} \leq C_{1}\norm{f}_{W_{2}^{s}\paren{\cX}}}$. We refer to the set of continuous representatives of all equivalence classes in $W_{2}^{s}\paren{\cX}$ as to \textit{Sobolev RKHS} and denote it as $W^{s}\paren{\cX}$. It can be shown (see paragraph 7.5 and Theorem 7.13 in \cite{Schaback:07} ) that $W^{s}\paren{\cX}$ is indeed a RKHS. Furthermore  (see part (c)Theorem~7.34 ) when $p \geq 2$, $W_{p}^{s}\paren{\cX}$ is embedded into the space of continuous functions $C\paren{\cX}$ if $s > d/p$ while if $s< \frac{d}{2}$ is not (and not embeddable into) a RKHS.


Furthermore, (see Chapter 7 in \cite{Schaback:07}), Sobolev RKHS $W^{s}\paren{\cX}$ is generated by the translation invariant kernel, which is a restriction to $\cX$ of the kernel $k_{s}$ of $\smash{W^{s}\paren{\mbr^d}}$ (see also Corollary 10.48 on page 170 in \cite{Wendlandt:05}). It is a continuous, bounded and measurable kernel (see general Lemma 4.28 and 4.25 in \cite{Steinwart:08} ) which is defined for all $x,x^{'} \in \cX$ by 
\begin{equation}
    \label{eq:Sobolev_kernel}
    k\paren{x,x'} \bydef \frac{2^{1-s}}{\Gamma \paren{s}} \norm{x-x'}_{2}^{s-\frac{d}{2}}K_{\frac{d}{2}-s}\paren{\norm{x-x'}_{2}} \,,
\end{equation}  
where $K_{d/2-s}\paren{\cdot}$ is a modified Bessel function of the second kind (see Chapter 5.1 in \cite{Wendlandt:05} for more details on Bessel function). 
Alternatively, the kernel function $k\paren{\cdot}$ of Sobolev RKHS $W^{s}\paren{\cX}$ can be described by
its Fourier transform, which equals $\cF(k)\paren{\omega} = (1+\|\omega\|_2^{2})^{-s}$. 
We refer the reader to the Chapters 10-11 in \cite{Wendlandt:05} as well as to \cite{Novak:17} for more details on the kernel functions of Sobolev RKHS. 

\subsection{\algo} 
In this work, we analyse the regret achieved by \algo (\cite{Gammerman2004}), 
over the (Sobolev) RKHS $W^{s}\paren{\cX}$. The regret is measured with respect to the benchmark classes of bounded Sobolev balls $W_{p}^{\beta}\paren{\cX}$ which may have different regularity, i.e. we consider the case when $\beta \neq s$. 

\begin{algorithm}[!th]
\SetAlgoLined
\textbf{Parameters:} $d \geq 1$, $s > d/2$, and  $\tau >0$ \\
\textbf{Initialization:} define $k(\cdot,\cdot)$ as in~\eqref{eq:Sobolev_kernel}\;
 \While{$t \geq 1$}{
    observe $x_t \in \cX$; \\
    $\tilde y_t \bydef (y_1,\dots,y_{t-1},0)^\top$;\\
    $\tilde k(x_t) \bydef \big(k(x_1,x_t),\dots,k(x_{t-1},x_t),k(x_t,x_t)\big)$;\\
    $K_t \bydef \big(k(x_i,x_j)\big)_{1\leq i,j \leq t}$;\\
    forecast $\hat y_t \bydef \tilde{y}^{\top}_{t}\paren{K_{t}+\tau I_t}^{-1} \tilde k\paren{x_{t}}$; \\
    observe $y_t$;
 }
 \caption{\algo \citep{Gammerman2004} on Sobolev RKHS}
 \label{alg:kaar}
\end{algorithm}

\algo (see Algorithm~\ref{alg:kaar}) was first introduced in the case of adversarial sequential linear regression by \cite{Vovk:01} and \cite{Azoury:01}; further it was analyzed in \cite{Bianchi:06}, \cite{Rakhlin:14}, \cite{Gaillard:19} and applied to concrete forecasting problems including electricity~\citep{Devaine:13}, air quality~\citep{Mallet:09} and exchange rate~\citep{Amat:18} forecasting. It was extended to the case of general reproducing Hilbert spaces in \cite{Gammerman2004}, while \cite{jezequel2019efficient} provide a variation of the algorithm with the same regret and reduced computational complexity. 
In the case of Sobolev spaces, KAAR (Alg.~\ref{alg:kaar}) reads as follows. 
Let $\tau >0$, $s>{d}/{2}$ at round $t \geq 1$;
\algo predicts $\smash{\hat y_t \bydef \hat f_{\tau,t}(x_t)}$, where 

\begin{equation}
\label{eq:kern_ridge}
\hat{f}_{\tau,t} \bydef \argmin_{f \in W^s(\cX) } \bigg\{ {\sum_{j=1}^{t-1}\paren[1]{y_{j}-f\paren[0]{x_{j}}}^{2}} + \tau\norm{f}_{W_2^s(\cX)}^{2}  + f^{2}\paren{x_{t}}\bigg\}.
\end{equation}

The prediction $\hat y_t = \hat f_{\tau,t}(x_t)$ can be computed in the closed form by Algorithm~\ref{alg:kaar} in $\cO(n^3 + n^2d)$ operations (see Section~\ref{sec:complexity} for details on the computational complexity). This improves computational complexity over other known nonparametric online regression algorithms, which achieve optimal regret with respect to Sobolev spaces in dimension $d$. 

\begin{remark}
	We remark that the right-hand-side of~\eqref{eq:kern_ridge} depends on the input $x_t$, so while $\hat{f}^{\tau}_{t} \in W^{s}\paren{\cX}$, the prediction function  $\smash{x_{t} \mapsto \hat f_{\tau,t}\paren{x_{t}}}$ is a measurable function which in general not necessarily belongs to the space $W_2^s(\cX)$, thus the prediction map not necessarily belongs to the benchmark class against which the algorithm is competing with.
	This corresponds to the so-called case of improper learning (see more details in \citep{Rakhlin:15,Hazan:18}). Furthermore, a sequential version of kernel ridge regression  was considered by~\cite{Zhdanov:10}. It removes the term $f^{2}(x_t)$ in the r.h.s. of~\eqref{eq:kern_ridge} and clips the prediction, by forecasting  $ \hat{y}^{M}_{t} \bydef \smash{\hat y_t \bydef \min(\max(-M,\tilde{f}_{\tau,t}(x_t)),M)}$, where $\tilde{f}_{t}$ is the solution to the Problem \ref{eq:kern_ridge} without $f^{2}\paren{x_{t}}$ term. In the case of nonlinear estimator \algo, for the clipped version of the \algo forecaster $\hat{y}^{M}_{t}$, since for every $y_{t} \in [-M,M]$ we have $\paren{y_{t} - \hat{y}^{M}_{t}}^{2} \leq \paren{y_{t}-\hat{y}_{t}}^{2}$ so the upper bound regret analysis for \algo can directly be applied to its clipped version. 
\end{remark}

We emphasize that throughout the paper $\beta$ and $p$ refer to the parameters of the benchmark Sobolev space and $s >{d}/{2}$ refers to the smoothness parameter of RKHS $W^{s}\paren{\cX}$ used in \algo.

\section{Main results: Upper bound on the regret of \algo on the classes of Sobolev balls.}
\label{sec:main_res}

In this section, we present regret upper bounds of \algo on the reference classes of bounded balls in $W_{p}^{\beta}\paren{\cX}$, $\beta>\frac{d}{p}$. 
{ By Sobolev embedding Theorem (see \cite{Adams:03} Theorem 7.34 or ,say, Equation 10 on page 60 in \cite{Edmunds:96}), condition $\beta>\frac{d}{p}$  implies that every equivalence class in  $W_{p}^{\beta}(\cX)$ has a continuous representative. In our analysis under $R_{n}\paren{B_{W_{p}^{\beta}\paren{\cX}}\paren{0,R}}$ we always understand regret with respect to the correspondent ball of continuous representatives bounded in the norm of the space ${W_{p}^{\beta}\paren{\cX}}$} \citep{Adams:03}. 
We consider the framework of online adversarial regression with the label space $\smash{\cY \bydef [-M,M]}$, target space $\hat{\cY} \subset \mbr$, the input space being $\smash{\cX = [-1,1]^d}$ and the reference class $\cF \bydef B_{W_{p}^{\beta}\paren{\cX}}\paren{0,R}$ being an open ball in Sobolev space $W_{p}^{\beta}\paren{\cX}$ of radius $R > 0$ with $d \in \mbn_{*}, \beta \in \mbr_{+}$ and $p \geq 2$, where we use standard notation for $\mbn_{*} = \{1,2\ldots,\}$. We remark that the assumption on the input space is given for simplicity and can be weakened to any bounded domain in $\mbr^d$ with Lipschitz boundary ( see Chapter 4 in \cite{Adams:03} on more details on Lipschitz boundaries).

\subsection{Key preliminary result and the upper bound on the effective dimension.}\label{sec:eff_dim}


We start by recalling a general upper bound on the regret of \algo on the bounded balls of the general separable  RKHS in terms of the effective dimension. It is a direct extension of the upper bound of \algo in \cite{Vovk:01,Azoury:01} from finite dimensional linear regression to kernel regression and can be retrieved from Theorem~2 in \cite{Gammerman2004} (see also  Proposition~1 and 2 in \cite{jezequel2019efficient} for the next statement) for the case of Sobolev RKHS, as the underlying kernel function is continuous. The regret of \algo on any $f \in W_2^s(\cX)$ is upper-bounded as
\begin{equation}
\label{eq:reg_bound}
R_{n}\paren{f} \leq \tau \norm{f}^{2}_{W_2^s(\cX)} +M^{2}\paren[3]{ 1 +\log\paren[2]{1 + \frac{n\kappa^2}{\tau}}}d^{n}_{eff}\paren{\tau},
\end{equation}
where $\kappa >0$ is such that $\sup_{x} k(x,x) \leq \kappa^2$ and $d^{n}_{eff}\paren{\tau}$ is the effective dimension 
as given in Definition~\ref{def:eff_dim}. The regret bound~\eqref{eq:reg_bound} will be used as a starting point to prove different upper bounds in the next subsection. 

To apply the equation above, in the next theorem we provide a bound on the effective dimension for the Sobolev space $W^s\paren{\cX}$.

\begin{theorem}[Upper bound for the effective dimension of Sobolev RKHS]
	\label{prop:eff_dim_sob}
	Let $\epsilon \in (0,1/4)$, $d\geq 1$, $n > \log_2\paren{d}$ and $s > {d}/{2}$. Consider the Sobolev RKHS $W^{s}\paren{\cX}$ with $\smash{\cX \bydef [-1,1]^{d}}$. For any sequence of inputs $x_1,\dots,x_n \in \cX$, the effective dimension $d_{eff}\paren{\tau}$ is upper-bounded as 
	\[
	    d_{eff}^{n}\paren{\tau} \leq  C  \paren{ \Big(\frac{n}{\tau}\Big)^{\frac{d}{2s} + \epsilon_{1}} +1 } \,,
	\]
	where $\epsilon_1 = d\epsilon/s^2$, and $C$\footnote{Throughout the paper, we refer to constants $C, C_{1}, etc$ which may depend on the properties of the domain $\cX$, the functional class $\cF$ or other quantities (such as $\epsilon$) but are always independent of $n$. We refer also to $\epsilon,\epsilon_1$ as to some infinitesimal numbers (possibly zeros). Their exact values are omitted and may differ from a statement to another, but we will specify this dependency in case this will be necessary for analysis.} is a constant which depends on $d,s,R,K,M, \cX,\epsilon$, but is independent of $n$. Furthermore, if $s \in \mbn_*$, then $\epsilon = 0$. 
\end{theorem}

The proof of this statement is presented in \cref{sec:d_eff_bounds}. It is based on some known properties of low rank projections in Sobolev spaces which are recalled in \cref{app:sob_spaces}.

\subsection{Regret upper bound for the Sobolev RKHS ($\frac{\beta}{d} > \frac{1}{2}$).}
\label{sec:bound_sob_rkhs}

 Notice that when $p \geq 2$ and $\beta \geq d/2$ we have $W_{p}^{\beta}\paren{\cX} \subseteq W_{2}^{\beta}\paren{\cX}$ and (by Sobolev embedding Theorem) $W_{p}^{\beta}\paren{\cX} \mapsto C\paren{\cX}$. The space of continuous representatives of $W_{p}^{\beta}\paren{\cX}$ is a closed subspace thus it is a RKHS (since it is a subspace of $W^{\beta}\paren{\cX}$ which is RKHS).  Using \algo with $s=\beta$ and putting the upper bound for the effective dimension of $W^{\beta}\paren{\cX}$ into the regret upper bound~\eqref{eq:reg_bound} 
with the proper choice of the parameter $\tau := \tau_{n}$, we obtain the following result. 

\begin{theorem}
	\label{thm:reg_tar} 
    Let $\cX \bydef [-1,1]^d$, $\beta \in (\frac{d}{2},+\infty)$, $p\geq 2$,  $M > 0$, and $n > \log_{2}\paren{d}$, $n \in \mbn$. 
    Then for any datasample $\{\bx_{t},y_{t}\}_{t = 1}^{n} \in \paren{\cX \times \cY}^{n}$, any $\epsilon>0$ regret of the \algo  with
    $$s = \beta, \quad \smash{\tau_{n} \bydef n^{\frac{d}{2\beta+d}}},$$
    on the benchmark class $\smash{\cF \bydef B_{W^{\beta}_{p}(\cX)}\paren{0,R}}$ satisfies the following upper bound
	\begin{equation*}
	    \frac{R_{n}\paren{\cF}}{n} \leq C n^{-\frac{2\beta}{2\beta +d} + \epsilon}\log\paren{n}\,,
	\end{equation*}
	where constant $C$ may depend on $d,s,R,K,M, \cX$, and $\epsilon$, but not on $n$.  
\end{theorem}
Proof of Theorem~\ref{thm:reg_tar} is given in \cref{app:proof_higher_beta}
\begin{remark}
	In the lower-bound section we prove that the upper bound of Theorem~\ref{thm:reg_tar} matches the minimax optimal for $\beta> d/2$ on the class of bounded Sobolev balls ( modulo a constant $\epsilon$ in the exponent that can be made arbitrarily small and logarithmic term in the number of observations). This rate was achieved by \cite{Rakhlin:14}  by a non-constructive procedure. An explicit forecaster has been proposed in \cite{Gaillard:15};  it can be calculated efficiently when $p = \infty$ and $d = 1$ and in general has exponential time and storage complexity. We believe that Theorem~\ref{thm:reg_tar} is the first (essentially) optimal regret upper bound for the classes of bounded balls in Sobolev spaces $W_p^\beta(\cX)$ with $d\geq 1$, $\beta > \frac{d}{2}$, and $p \geq 2$ that is achieved by a computationally efficient procedure. 	
\end{remark}

\subsection{Regret upper bound over Sobolev spaces when $ \frac{d}{p}< \beta \leq \frac{d}{2}$, $p\geq 2$.}
\label{sec:bound_sob_dif}
In this part we consider \algo over the benchmark classes of bounded balls $\smash{W_p^\beta(\cX)}$ when $ \frac{d}{p} < \beta \leq \frac{d}{2}$, $p >2$ and refer to this case as to "hard learning" scenario. When $\frac{\beta}{d} \leq \frac{1}{2}$the Sobolev space $\smash{W_p^\beta(\cX)}$ is not included into any Sobolev reproducing kernel Hilbert space so using \algo in this case we need to control the error due to using the element $\hat{f}^{\tau}_{t} \in W^{s}\paren{\cX}$ when competing against any function from  $W_{p}^{\beta}\paren{\cX}$.
 In this case the regret analysis can be decomposed into two parts: approximation of any function $f \in W_{p}^{\beta}\paren{\cX}$ by some element $f_{\epsilon} \in W^{s}\paren{\cX}$ and regret of \algo with respect to bounded balls in $W^{s}\paren{\cX}$.
 Intuitively, the smaller the approximation error between $f$ and  $f_{\epsilon}$, the larger the norm of the approximation function $f_\epsilon$ should be, which implies the larger regret upper bound of \algo with respect to $f_\epsilon$ (see bound \eqref{eq:reg_bound}). Therefore, in this case one has to control a trade-off between the approximation error of $f \in W_p(\cX)$ by means of some $f_\epsilon \in W^s(\cX)$ and the regret suffered with respect to $f_\epsilon$. We have the following result.

\begin{theorem}
\label{thm:higher_beta_p}
    Let $\cX = [-1,1]^d$, $p > 2$, $\beta \in \mbr_{+}$, $d/p < \beta  \leq d/2$,  $M >0$, $\epsilon >0$, $n\geq \log_{2}d$ and $\{(\bx_t,y_{t})\}_{t=1}^{n} \in (\cX \times [-M,M])^n$ be arbitrary sequence of observations. Then by choosing $s = \frac{d}{2}+ \epsilon$ and
$$\tau_{n} = n^{1- \frac{d\paren{1-p^{-1}} -\beta^{'}}{d\paren{1-2p^{-1}}} }$$
 where $\beta^{'}= \beta-\epsilon$ is sufficiently close to $\beta$ decision rule~\ref{alg:kaar} of \algo
satisfies the following regret upper bound
\[
    R_{n}\paren{\cF} \leq Cn^{1- \frac{\beta}{d}\frac{p-\frac{d}{\beta}}{p-2} + \epsilon \theta}\log\paren{n},
\]
where ${\cF = B_{W_{p}^{\beta}\paren{\cX} }\paren{0,R}}$, and $R>0$. Constant $C$ depends on $d,s,R,\beta,M, \cX$, and $\epsilon$, but not on $n$ and constant $\theta = \frac{p}{\paren{p-2}d}$.
\end{theorem}
The proof of the Theorem is given in \cref{app:proof_higher_beta}. The Theorem and its implications are discussed in Section~\ref{sect:comparison}. Here we want just to provide two remarks that help to interpret the result.
\begin{remark}
    In the proof we provide the regret upper bound for any choice $s >\frac{d}{2}$, however the rate for $ \frac{d}{p} \leq \beta \leq \frac{d}{2}$ is minimized by the choice $s>\frac{d}{2}$ as small as possible. Therefore in this situation we choose $s \bydef d/2 + \epsilon$ with an arbitrary small $\epsilon>0$. Furthermore the result of Theorem~\ref{thm:higher_beta_p} is of asymptotic nature as it includes condition $n \geq n_{0}$, where $n_{0}$ depends exponentially on $d$. To the best of our knowledge this dependence is unavoidable when using techniques which we use in this work.
\end{remark}

\begin{remark}
   Notice that in an interesting particular case of Theorem~\ref{thm:higher_beta_p} when $p = \infty$ and $\beta \in \mbr_{+}$ the space $W^\beta_\infty(\cX)$ corresponds to functions with derivatives up to order $\lfloor \beta \rfloor$ bounded in supremum norm and $\lfloor \beta \rfloor$-th derivatives are Hölder continuous of order $\alpha  \in (0,1)$ \citep{Adams:03}. Then the regret of Theorem~\ref{thm:higher_beta_p} leads to a regret upper bound of order $O(n^{1-\frac{\beta}{d} + \epsilon} \log n)$. This upper bound is optimal on the class $W_{\infty}^{\beta}\paren{\cX}$, up to a negligible factor $\epsilon$ that can be made arbitrary small (see Section~\ref{sec:lower_bound}).
\end{remark}


\section{Lower bounds}
\label{sec:lower_bound}

In this section, we present lower-bounds on the regret of any algorithm on the bounded closed balls in Sobolev spaces $W_{p}^{\beta}\paren{\cX}$ with $\beta > {d}/{p}$, $p \geq 2$. 
We define the minimax regret for the problem of online nonparametric regression on the functional class $\cF$ as 
\begin{align}
\label{def:minimax_regret_al}
\smash{\tilde R_n(\cF) \bydef \inf_{\cA} \sup_{(\bx_s,y_s)_{s \leq n} \in  \paren{\cX \times \cY}^{n}}  R_n(\cF)},
\end{align}
where $\cA = \paren{\cA_{s}}_{s\geq 1}$ is any \textit{admissible} forecasting rule, i.e. such that at time $t \in \mbn$ outputs a prediction $\hat y_t \in \hat{\cY}$ based on past predictions $\paren{\hat{y}_{s}}_{s \leq t-1}$ and data-sample $ \paren{\{x_{s},y_{s}\}_{s \leq t-1} \cup x_{t}}$. More formally, we assume $\cA$ is such that the map $\cA_{t} : \paren{\hat{\cY}^{t-1} \times \paren{\cX \times \cY}^{t-1} \times \cX } \mapsto \hat{\cY}$ is measurable and call such algorithm admissible. The most important element of this technical assumption is that the forecaster cannot use the future outcomes for making current decisions. Notice that in this setting we consider the oblivious adversary meaning that all outputs $\paren{x_{t},y_{t}}_{t\geq 1}$ are fixed in advance. With this notation we have the following result.
\begin{theorem}
	\label{thm:lower_bounds_2} 
	Let $M >0$, $p\geq 1$, $\beta > d/p$ and $\cX \bydef [-1,1]^{d}$ as before. 
	Consider the problem of online adversarial nonparameteric regression with $ y_{t} \in [-M,M]$ over the benchmark class $\smash{\cF \bydef B_{W_{p}^{\beta}\paren{\cX}}\paren{0,M}}$. Then minimax regret from Equation~\eqref{def:minimax_regret_al} is lower-bounded as
	\begin{align*}
			\tilde{R}_{n}\paren{\cF} \geq \left\{ 
				\begin{array}{ll}
					C_{1}n^{1-\frac{\beta}{d}} & \text{if} \quad \beta \leq \frac{d}{2} \\
					C_{2}n^{1- \frac{2\beta}{2\beta +d}} & \text{if} \quad \beta > \frac{d}{2}
				\end{array}\right. \,,
	\end{align*}
	where $C_{1}$ and $C_2$ are constants which depend on $M,\cX,d,\beta$, and $p$, but are independent of $n$. 
\end{theorem}
The proof is based on the general minimax lower bounds of~\cite{Rakhlin:14} and is given in \cref{app:proof_lower_bound}.
\begin{table}[H]
	\begin{center}
		\renewcommand{\arraystretch}{1.3}
		\resizebox{\textwidth}{!}{
			\begin{tabular}{@{}lllclcl}
				\toprule
				{} & \multicolumn{2}{c}{ Statistical i.i.d. regression} & \phantom{ab}& \multicolumn{2}{c}{ Adversarial online nonparametric regression} \\
				\cmidrule{2-3} \cmidrule{5-6} 
				& Best known excess risk upper bound  & Lower bound  & & Best known upper bound for $R_{n}\paren{\cF}/n$ & Lower bound  \\ 
				\midrule
				$\beta > \frac{d}{2}$ & \centering{\textbf{$n^{-\frac{2\beta}{2\beta+d} }$}} & \centering{\textbf{$n^{-\frac{2\beta}{2\beta+d} }$}} & & \centering { \color{blue}{\textbf{$n^{-\frac{2\beta}{2\beta+d} }$}} } & \color{blue}{\textbf{$n^{-\frac{2\beta}{2\beta+d} }$}}  \\
				$\frac{d}{p}< \beta \leq \frac{d}{2}$ &  {$n^{-\frac{2\beta}{2\beta+d}}$}  & {\textbf{$n^{-\frac{1}{2}}$}} && \textbf{$n^{-\frac{\beta}{d} }$}  & {$n^{-\frac{\beta}{d}}$} \\ 
				$p=\infty, \beta \leq \frac{d}{2} $ &  {$n^{-\frac{2\beta}{2\beta + d}}$}                    & {$n^{-\frac{2\beta}{2\beta+d}}$}&&  \color{blue}{$n^{-\frac{\beta}{d}}$}                        &   \color{blue}{$n^{-\frac{\beta}{d}}$} \\
				\bottomrule
			\end{tabular}
		}
	\end{center}
	
	\caption{Best known regret and excess risk upper and lower bounds on the classes of Sobolev balls. Results achieved by \algo are highlighted with blue color.}
	\label{tab:main_tabl2}
\end{table}

\begin{remark}
	In Table~\eqref{tab:main_tabl2} we compare the best known lower and upper regret bounds on the classes of (continuous representatives) of Sobolev balls in the settings of adversarial online regression to the correspondent bounds for the excess risk in the statistical i.i.d. scenario.  Interestingly, on the classes of Sobolev balls in spaces $W_{p}^{\beta}\paren{\cX}$, $\beta \geq \frac{d}{2}$ and H\"older balls $W_{\infty}^{\beta}\paren{\cX}$ rates for the (normalized) regret and for the excess risk are optimal and archived by the regularized empirical risk minimization procedure (for example by regularized least squares estimators in the statistical learning scenario, see \cite{Fischer:17}) and \algo as shown in this work.
\end{remark}

\section{Discussion}
\label{sect:comparison}
In this part we compare regret rates of \algo with the existing algorithms in the adversarial online nonparametric regression in terms of regret bounds and computational complexity. Furthermore we compare the regret bounds to the excess risk bounds for the known algorithmic schemes in the statistical least-squares regression scenario. We point out on interesting consequences for the gap in the rate which arises due to adversarial data.

\subsection{General comparison to the setting of statistical nonparametric regression}

 To unify settings we always consider the normalized regret of class $\cF$, $n^{-1}{R_{n}\paren{\cF}}$. In the setting of statistical learning we assume a sample $\cD_{n} = \paren{(x_i,y_i)}_{i = 1}^{n}$ is generated independently from the distribution $\nu_{x,y}$ of a pair of random variables $X,Y$ over a probability space $\paren{\cX \times \cY, \cB\paren{\cX\times \cY},\nu}$ and let $f_{D_{n}} : \cX \mapsto \cY$ be some data-dependent estimator of the output $y$. Denote $\nu\paren{\cdot|x}$ to be a regular conditional probability distribution of $Y$ given $\{X=x\}$ and $\mu$ to be the $X-$marginal of $\nu$. In statistical nonparametric regression scenario the performance measure of data-dependent estimator $f_{\cD_{n}}$ is controlled through the excess risk  $\cE\paren{f_{\cD}} - \inf_{f \in \cF}\cE\paren{f}$, where $\cE\paren{f} = \ee{}{\paren{Y-f\paren{X}}^{2}}$. If $\cF$ is dense in $L_{2}\paren{\cX,\cB\paren{\cX},\mu}$ the latter is equivalent to $\norm[1]{f_{\nu} - f_{\cD}}_{L_{2}\paren{\cX,\cB\paren{\cX},\mu}}$, where $f_{\nu}(x)$ is for $\mu-$almost all $x$ is a version of conditional expectation of $y$ with respect to $v\paren{\cdot|\cdot}$. 
For comparison we consider data-dependent estimators with values in a Sobolev RKHS $\cH \bydef W^{s}\paren{\cX}$ and assume that $f_{\nu} \in \cF \subset W_{p}^{\beta}\paren{\cX,\mu}$ but in general that $f_{\nu} \notin \cH$. We denote $W_{p}^{\beta}\paren{\cX,\mu}$ for the Sobolev space with respect to the measure $\mu$. To avoid technical difficulties with threatening weak-derivatives with respect to arbitrary Borel measure, we assume $\mu$ to have upper and lower bounded Radon-Nikodym derivative with respect to Lebesgue measure over $\cX$. In this case $W_{p}^{\beta}\paren{\cX,\mu}$ is equivalent to the standard Sobolev space $W_{p}^{\beta}\paren{\cX}$. If other is not stated, the focus will be on the excess risk upper bounds in high probability, namely under $\norm{f_{\cD_{n}} - f_{\nu}}_{L_{2}\paren{\cX}} \leq C(\delta)\psi(n)$ we understand inequality which holds with probability at least $1 - \exp(-\delta)$ for some $\delta>0$, where $C\paren{\cdot} : \mbr_{+} \mapsto \mbr_{+} $ and $\psi\paren{\cdot} : \mbr_{+} \mapsto \mbr_{+}$ are some functions. We consider $\cX = [-1,1]^{d}$ as before, however all the subsequent results in the statistical regression scenario can be reformulated for any bounded subset of $\mbr^{d}$ with Lipschitz boundary.
 
 We start with the case, in which $f_{\nu} \in \cF \subset \cH$, $\cH$ is the Sobolev RKHS. Theorem~1 in \cite{Caponentto:06} implies (by taking $b=\frac{2\beta}{d}$ and $c=1$ therein) that for $\cH= W_{}^{\beta}\paren{\cX}$, $\beta > \frac{d}{2}$, $f_{\nu} \in B_{W_{2}^{\beta}\paren{\cX,\mu}}\paren{0,1}$ and $f_{\cD_{n}} \in \cH$ being a regularized least-squares estimator, we obtain that it holds $\norm{f_{\cD_{n}} - f_{\nu}}_{L_{2}\paren{X,\nu}} \leq Cn^{-\frac{2\beta}{2\beta +d}}$ which is the optimal rate in the setting of nonparametric regression (see \cite{Tsybakov:09} and \cite{Gyorfi:02} for matching lower bounds). Under the same conditions, optimal excess risk rates on $B_{W_{2}^{\beta}\paren{\cX,\mu}}\paren{0,1}$ can be deduced from Corollary~6 in \cite{Lin:18} using the decision rule based on the spectral kernel algorithms or stochastic gradient descent. 
 It follows that the regret rates of \algo on classes $W_{2}^{\beta}\paren{\cX}$ essentially match (disregarding arbitrary small polynomial factor) the optimal known rates for the excess risk in the i.i.d. scenario on classes $W_{2}^{\beta}\paren{\cX,\mu}$.
 
 The setting in which the underlying RKHS is a subspace of reference class of regular functions is studied in several works. In the particular case of $H := H_{\gamma}\paren{\cX}$ being a Gaussian RKHS over $\cX$, 
 $f_{\nu}\paren{\cdot} \in W_{2}^{\beta}\paren{\cX} \cap L_{\infty}\paren{\cX}$, $\beta \in \mbr_{+}$, Corollary~2 in \cite{Eberts:11} implies that the Gaussian kernel ridge regression estimator with the proper choice of {both} regularization parameter $\lambda$ {and} band-width $\gamma$ achieves essentially optimal rates for excess risk of order $n^{-\frac{2\beta}{2\beta+d} + \epsilon }$ when $\beta > \frac{d}{2}$, $\epsilon>0$. This rate hold when $\beta \leq \frac{d}{2}$ under additional condition $\cY = [-M,M]$ (which implies $\nu-$a.s. boundedness of $f_{\nu}$ that is not ensured unless $\beta > \frac{d}{2}$) however in this case it remains to be unknown whether the given rate is optimal on the given classes. 

In the case $f_{\nu} \in W_{\infty}^{\beta}\paren{\cX}$ (i.e. it has partial derivatives of order up to $ \lfloor \beta \rfloor$, and every partial derivative of order $\lfloor \beta \rfloor-$ is $\beta - \lfloor \beta \rfloor$ H\"older continuous) when $\beta \leq \frac{d}{2}$, excess risk upper bound of order $n^{-\frac{2\beta}{2\beta +d} + \epsilon}$ in the statistical i.i.d. scenario is essentially optimal (see Chapter 3.2, Theorem 3.2 in \cite{Gyorfi:02} for a lower minimax rate of convergence which implies the lower bound). This rate is better in comparison to the normalized regret rate of \algo ( $n^{-\frac{\beta}{d} + \epsilon}$) which in it's turn is 
essentially optimal in the adversarial Sobolev Regression setting. 

The latter two facts uncover an interesting consequence, namely that the gap between the optimal rates for regret in the setting of adversarial regression and the excess risk in the statistical setting on classes of bounded balls in $W_{\infty}^{\beta}\paren{\cX}$ is due purely to the adversarial nature of the data.

If $f_{\nu} \in W_{2}^{\beta}\paren{\cX,\mu}$ Corollary~6 in \cite{Steinwart:09} and their discussion afterwards implies that excess risk upper bounds of kernel ridge-regression least squares estimator based on the Sobolev kernel of finite smoothness
$s\geq \beta> \frac{d}{2}$ are of order $n^{-\frac{2\beta}{2\beta +d}}$ and thus optimal. Notice that in this case we do not need to know the smoothness parameter $\beta$ but only the (possibly crude) upper bound $s$. Similarly, Theorem $1$ and Example $2$ in \cite{Vivien:18} imply the excess risk rates (in expectation) for the stochastic gradient descent estimator with values in Sobolev space $W_{}^{s}\paren{\cX}$ over the class $W_{2}^{\beta}\paren{\cX,\mu}$,$ \frac{d}{2} < \beta < s $. They are optimal under an additional assumption $s-\beta \geq \frac{d}{2}$. Corollary 4.4 in \cite{Lin:20} implies risk upper bound for a general spectral kernel decision rule of order $n^{-\frac{2\zeta}{ (2\zeta + \gamma) \vee 1}}$  where parameter $\zeta$ is the power of the so-called source condition (see \cite{Engl:00} also see \cite{Blanchard:07} for the statistical perspective)
and $\gamma $ is the decay rate of effective dimension. Adapting this to the case of Sobolev regression over the space $W_{2}^{\beta}\paren{\cX}$ with decision rules valued in $W^{s}\paren{\cX}$ we get $\zeta = \frac{\beta}{2s}$ and $\gamma = \frac{d}{2s}$, $\beta \leq \frac{d}{2}$, $s>\frac{d}{2}$) . If $s>d$ we have the excess risk upper bound of order $n^{-\frac{\beta}{s}}$ which is worse than the rate $n^{-\frac{2\beta}{2\beta + d}}$. If $\frac{d}{2} < s \leq d$ the excess risk upper rate is $n^{-\frac{2\beta}{2\beta+d}}$ when $s-\frac{d}{2} < \beta \leq \frac{d}{2}$ and $n^{-\frac{\beta}{s}}$ when $0 < \beta <s- \frac{d}{2}$. In the latter case, on the classes of H\"older balls, the rate is better then the lower bound on the minimax regret, being worth then $n^{-\frac{2\beta}{2\beta + d}}$ achieved, as stated above, by, for example, regularized least squares estimator with Gaussian kernels. In the worse case scenario ($\beta < \frac{d}{2}$, $\beta +\frac{d}{2} <s$) one also observes the gap between upper rates for the excess risk in the statistical learning scenario achieved by general spectral regularization methods ($n^{-\beta/s}$) and the lower bounds for the minimax regret  ($n^{-\frac{\beta}{d}}$) in the online regression setting. 

 A broader analysis of the quantity $f_{\cD_{n}} - f_{\nu}$ 
 in the norms of the interpolation Hilbert spaces (which in its particular case uncovers the excess risk rates) 
 which ranges between $\cH$ and $L_{2}\paren{\cX}$ for the kernel ridge-regression estimator is given in \cite{Fischer:17}. Corollary 4.1 therein and inclusion between Sobolev spaces allow to deduce the excess risk upper bounds of order $n^{-\frac{2\beta}{2\beta+d} + \epsilon }$ for $f_{\nu} \in W_{p}^{\beta}\paren{\cX}$, $\beta >0$, $p \geq 2$. If $\frac{\beta}{d} \in (\frac{1}{p}, \frac{1}{2}]$ and $p \geq 2$ then the aforementioned excess risk rates are better then the regret upper bounds obtained by \algo on the same classes. To the best of our knowledge, the best known lower bounds in probability on the excess risk on the classes of balls in the Sobolev spaces are of order $n^{-\frac{1}{2}}$ (see  Corollary 4.2 in \cite{Fischer:17} with $t =0$, $f_{\nu} \in W_{p}^{\beta}\paren{\cX} \subset W_{2}^{\beta}\paren{\cX}$ and notice that  $f_{\nu}$ is bounded on $\cX$ by Sobolev embedding and Bolzano-Weierstrass theorem).  
 
\subsection{ Comparison in the setting of adversarial nonparametric regression.}

\paragraph{Previous works on online nonparametric regression and optimal rates.}

The setting of online nonparametric regression is definitely not new. The standard idea is to use an $\epsilon$-net of the bounded functional space and exploit the exponential weighted average (EWA) forecaster for a finite class of experts which will be the element of the $\epsilon-$net (see Chapter 1 in the monograph \cite{Bianchi:06} for the finite EWA and \cite{Vovk06Metric} for its application in the nonparametric case). This leads to the regret rate of order $n^{\frac{d}{\beta+d}}$, which were known to be suboptimal in the case $\beta> \frac{d}{2}$. \textit{Minimax} regret analysis in terms of (sequential) entropy growth rates of the underlying functional classes was provided by \cite{Rakhlin:14}. In particular, the optimal rates of order $n^{\frac{d}{2\beta +d}}$ (up to a logarithmic terms) when the reference class is Sobolev RKHS ( $\beta> \frac{d}{2}$ ) and of order $n^{1-\frac{\beta}{d}}$ on the classes of H\'older balls (which correspond to classes $B_{W_{\infty}^{\beta}\paren{\cX}}(0,R)$) can be achieved by using the generic forecaster with Rademacher complexity as a relaxation (for more details see Example 2, Theorems~2,3 and Section~6 in \cite{Rakhlin:14}).  
Although the relaxation procedure ensures minimax optimality, it is not constructive in general. An explicit forecaster, which designs an algorithm based on a multi-scale exponential weighted average algorithm (called Chaining EWA), has been provided in \cite{Gaillard:15} . The latter achieves an optimal rate when competing against functional classes of uniformly bounded functions which have certain (sharp) growth condition on the sequential entropy (see \cite{Rakhlin:14}). This condition implies optimal rates, for example on classes where sequential entropy is of order of metric entropy (which holds for example if for Sobolev classes $W^{\beta}\paren{\cX}$). 

Chaining EWA has been shown to be computationally efficient on the class of H\"older balls ($p=\infty$) with $d=1$. In general, the Chaining EWA forecaster is computationally prohibitive (as it has exponential time complexity in the number of rounds). \cite{Vovk:06} analyzes the regret when competing against a general reproducing kernel Hilbert space defined on an arbitrary set $\cX \subset \mbr$ 
and proves in this case the existence of an algorithm (which is based on the so-called idea of defensive forecasting and requires the knowledge of the feature kernel map) with the regret of order $\cO(\sqrt{n})$ over
unit balls within the general reproducing Hilbert space. \cite{Vovk:07} extends the analysis to the more general framework of Banach spaces, which is described through the decay rate of the so-called modulus of convexity of the underlying Banach space (originally introduced by \cite{Clarkson:36}).  As a particular example it includes Sobolev spaces where the parameter $p$ of the modulus of convexity being the parameter of the $p-$ from the definition of $W_{\beta}^{p}\paren{\cX}$. 

We notice that all aforementioned approaches have the disadvantage of either having suboptimal regret bounds or having (in general) the prohibitive computational complexity.

\paragraph{Comparison with Exponential Weighted Average (EWA) forecaster.}
The idea of using of the EWA forecaster in the nonparametric setting over bounded benchmark functional class $\cW$ is to consider the $\epsilon-$ net $\cW_{\epsilon}$ of the smallest cardinality: 
\begin{align*}
    \cW_{\epsilon} \subset \cW, \cW_{\epsilon} = \min_{K} \{ f_{1},f_{2},\ldots,f_{K} : \forall f \in \cW \exists i \in \{1,\ldots,K\}, \text{ s.t.} \norm{f-f_{i}}_{\infty} \leq \epsilon \}
\end{align*}
and to use the (finite) EWA forecaster (see \cite{Bianchi:06}) on the set $\cW_{\epsilon}$. It was introduced in \cite{Vovk06Metric} (see also discussions in \cite{Rakhlin:14} and \cite{Gaillard:15}) and leads to the composed  regret upper bound of order $n\epsilon + \log\paren{\cN_{\infty}\paren{\epsilon,\cF}}$, where the last term is the metric entropy of class $\cF$ on scale $\epsilon$. It is known (see \cite{Edmunds:96}) that for the benchmark class of Sobolev spaces $W^{\beta}_{p}\paren{\cX}$ (with $p \geq 2$ and $\beta > {d}/{p}$), metric entropy is of order $\epsilon^{-\frac{d}{\beta}}$. 
Balancing the terms by a proper choice of $\epsilon$, it results in an upper bound of order $n^{{d}/(\beta+d)}$ (see also Corollary~8 in \cite{Vovk:06}). 
As is illustrated in
Figure~\ref{fig:opt-regions} in the $(\frac{\beta}{d}, p^{-1})$ plane, regret upper-bounds of \algo are smaller than that of EWA as soon as $\frac{\beta}{d}$ is large enough. More precisely, EWA outperforms KAAR when $\frac{\beta}{d} \in [\frac{1}{p}, \frac{\sqrt{1+4p}-1}{2p}]$. 
The latter is not surprising since \algo, which 
outputs prediction rules in Sobolev RKHS (i.e. functions of sufficiently high regularity), performs worse on the when competing against functions of small regularity. EWA does not have this drawback, as it acts through the space discretization.   

In the case $p\geq 2$ and $\frac{\beta}{d} \leq \frac{1}{p}$ it is generally not true that there exists a continuous representative for each equivalence class in $W_{p}^{\beta}\paren{\cX}$. In the case of additional continuity assumption (i.e. considering bounded balls in $W_{p}^{\beta}\paren{\cX} \cap \cC\paren{\cX}$ as a benchmark class instead) the best (known) upper bound for minimax regret (and thus for regret itself) is of order $n^{1-\frac{1}{p}}$ (see Example~2 in \cite{Rakhlin:14}). It is achieved by a non-constructive algorithm based on the notion of relaxation of sequential Rademacher complexity. Notice that EWA can be also applied over classes  $B_{W_{p}^{\beta}\paren{\cX} \cap \cC\paren{\cX}}\paren{0,1}$, $\frac{\beta}{d} \leq \frac{1}{p}$; here it provides the same rate $n^{\frac{d}{\beta+d}}$ which is worth than $n^{1-\frac{1}{p}}$.

\paragraph{Comparison with defensive forecaster by \cite{Vovk:07}.} \cite{Vovk:07} describes the algorithms that are based on the defensive forecasting schemes in general Banach spaces. The benchmark classes are irregular but continuous functions, particularly including Sobolev spaces. By transferring the results given in Equations (6) and (11) in \cite{Vovk:07} to the setting of this work, defensive forecaster BBK29 (see pages 19--20 in \cite{Vovk:07}) achieves for a unit ball $\cF= B_{W_{p}^{\beta}\paren{\cX}}\paren{0,1}$ the following regret bound
\[
    R_n(\cF) \leq \left\{ 
        \begin{array}{lll}
C n^{1-\frac{\beta}{d}+\epsilon} & \text{if}\quad   p = \infty \\
C n^{1 -\frac{1}{p}} & \text{if} \quad  2 \leq p < \infty & \text{and} \quad  \frac{d}{p} \leq \beta. 
\end{array}
\right. \,.
\]
Therefore, in the first case, which corresponds to H\"older balls in $W_{\infty}^{\beta}\paren{\cX}$ and $0<\beta\leq 1$, we recover the same rate as Theorem~\ref{thm:higher_beta_p} but for the range $\beta > 0$. The rate is optimal, as stated in  Theorem~\ref{thm:lower_bounds_2}. In the second case ($p\geq 2$ and $\frac{d}{p}<\beta<1$), the upper bounds provided by Theorem~\ref{thm:higher_beta_p} (if $\beta > \frac{d}{p}$) or Theorem~\ref{thm:reg_tar}  (if $d=1$ and $\beta >1/2$) are always better then the correspondent bounds of~\cite{Vovk:07}. 

\subsection{Computational complexity}
\label{sec:complexity}
Here we consider an optimal computational scheme for \algo and compare its costs with those of the known nonparametric algorithms (in terms of both runtime and storage complexity). 

Recall that \algo for any $x_{t} \in \cX$, $\paren{x_{s},y_{s}}_{s\leq t-1} \in \paren{\cX \times \cY }^{t-1}$ computes
\[
\hat{y}_{t} = \hat{f}_{\tau,t}\paren{x_t} = \inner[1]{\hat{f}_{\tau,t},k_{x_t}}_{\cH_{k}} = \sum_{s=1}^{t}k\paren{x_t,x_s}c_s \,,
\]
where $c \in \mbr^{t}, c= \paren{K_{t}+\tau \mbi}^{-1}\tilde{y}_{t} $, $\tilde{y}_{t}^{\top} = \paren{Y^{\top}_{t-1},0}$ and $K_t=\paren{k\paren{x_i,x_j}}_{i,j \leq t}$ is the kernel matrix at step $t$. A naive way to compute the value of \algo at the input $x_t$ 
is by computing the inverse of matrix $K_{t}+\tau \mbi_{t}$. This requires $\cO\paren{t^3}$ iterations in round  $t$ and implies $\cO\paren{n^4}$ cumulative time complexity over $n$ rounds. The letter can be improved by using the Cholesky decomposition and the rank-one update of the kernel matrix. Namely, we use the approach as in Algorithm~1 (see \cite{Rudi:15}) for general RKHS. More precisely, at time $t$ we compute the Cholesky decomposition $R_{t-1}R_{t-1}^{\top} = K_{t}+\tau \mbi$; next, we denote the following quantities
\begin{align*}
b_{t}& \bydef \paren{k\paren{x_{t},x_1},\ldots,k\paren{x_{t},x_{t-1}}} \qquad
\alpha_t \bydef K_{t-1}^{\top}b_{t} + \tau b_t\\
\gamma_t &\bydef a_{t}^{\top}a_t + \tau k\paren{x_t,x_t} \qquad \qquad \quad \quad g_t \bydef \sqrt{1+\gamma_t},
\end{align*}
and $u_{t}=(\frac{\alpha_{t}}{1+g_t},g_t)$,$v_{t}=(\frac{\alpha_{t}}{1+g_t},-1)$. 
Using this, we compute an update of $R_t$:
\[
R_{t} \bydef \begin{pmatrix}
R_{t-1} & 0  \\
0 & 0 
\end{pmatrix},
\\
R_{t}  \bydef \text{\textsc{cholupdate($R_{t}$,$u_{t}$,'+')}} , 
\\
R_{t}  \bydef \text{\textsc{cholupdate($R_{t}$,$v_{t}$,'-')}}
\]
and calculate the solution's coefficients $c_{t} = R_{t}^{-1}\paren{R_{t}^{\top}}^{-1}K_{t}\tilde{y}_{t}$. Notice that the procedure \textsc{cholup $(R,a,"+")$} returns the upper triangular Cholesky factor of $R+a^{\top}a$, whereas \textsc{cholup $(R,a," -")$} returns the upper triangule update of $R-a^{\top}a$. 
At round $t$ ($t \leq n$) its computational cost is at most $ \cO\paren{t^2}$. Taking into the account that at the end we compute kernel matrix $K_{n} = (k(x_i,x_j))_{i,j \leq n}$ for a $d-$dimensional input $x_t$, which adds $dn^2$ to the total computational complexity we obtain, that the total computational costs is of the order of $\smash{\cO(n^3 +n^2 d)}$ operations. The latter complexity can be further improved when $\smash{\beta > d/(2\sqrt{2}-2)}$ (which implies $\beta > d/2$) to ${\cO(n^{1+ \frac{2d/\beta}{(1-\paren{d/(2\beta)}^{2})}})}$ by using Nyström projection \citep{jezequel2019efficient} while retaining the optimal regret. In particular, it converges to linear runtime complexity when $\beta \to \infty$. \cite{jezequel2019efficient} also provides additional improvements to the complexity if features $x_t$ are revealed to the learner beforehand. 

As was mentioned before, most existing work in online nonparametric regression on Sobolev spaces ( in particular \citep{Rakhlin:14,Vovk06Metric,Vovk:06,Vovk:07}) does not provide efficient  (i.e., polynomial in time) algorithms. Work by \cite{Rakhlin:14} provides an optimal minimax analysis; however, they do not develop constructive procedures. 
More precisely, they require knowledge of the (tight) upper bounds for the so-called \textit{relaxations}. To obtain the latter ones, in general, one must compute the offset Rademacher complexity, which is numerically infeasible. The approach of using EWA in nonparametric setting (\cite{Vovk06Metric}) has non-optimal rates and suffers from prohibitive computational complexity because it must update the weights of the experts in the $\epsilon-$net. For Sobolev balls its size is of order $\cO\paren{\exp\paren{n}}$ (given that the number of experts scales as $\paren{\cN\paren{\cF}}$ with $\log{\cN\paren{\cF}}$ being the metric entropy of the class $\cF$, which is polynomial in the number of rounds) so that the total time complexity will be $\cO\paren{\exp{n} + nd}$ (where $nd$ comes from the aggregation of observations $x_t \in \cX \subset \mbr^d$ over $n$ rounds). The defensive forecasting approaches by \citep{Vovk:06,Vovk:07} require the knowledge of the so-called Banach feature map, which is typically inaccessible in the computational design of the algorithm.

To the best of our knowledge, the only algorithm that addresses the problem of computational cost in online nonparametric regression is the Chaining EWA forecaster ( \cite{Gaillard:15}). On class $W_{\infty}^{\beta}\paren{\cX}$ with $\beta = r+\alpha$, $\alpha \in (0,1]$, $r \in \mbn_*$, the Chaining EWA forecaster can be efficiently implemented through piecewise polynomial approximation ---see Lemma~12 and Appendix~C in \cite{Gaillard:15}. Its time and storage total complexities are of order:
\[
\text{Storage:}  \quad \cO\paren[1]{n^{2r+4 + \frac{\beta(r-1) +1}{2\beta +1}}\log\paren{n}}\,, \qquad \text{Time:} \quad \cO\paren[1]{n^{(r+1)(2+\frac{\beta}{2\beta +1})} \log\paren{n}} \,.
\]
Notice that storage complexity of \algo is $\cO\paren{n^2}$ and it is uniformly better for any $\beta=r+\alpha>0$ than of Chaining EWA. Furthermore, its time complexity is better for all $\beta \geq 1$ (and worth for $0< \beta <1$) than that of the efficient implementation of the Chaining EWA. As was mentioned in \cite{Gaillard:15}, in most of the cases the direct implementation of the Chaining EWA forecaster requires $\exp\paren{d poly\paren{n}}$ time (due to the exponentially many updates of the expert's coefficients). 

\paragraph{Acknowledgements}

Oleksandr Zadorozhnyi would like to acknowledge the full support of the Deutsche Forschungsgemeinschaft (DFG) SFB 1294 and the mobility support due to the UFA-DFH through the French-German Doktorandenkolleg CDFA 01-18.

The authors acknowledge the Franco-German University (UFA) for its support through the bi-national \emph{Coll\`ege Doctoral Franco-Allemand} CDFA 01-18.

\bibliography{advkern01}

\begin{thebibliography}{52}
\providecommand{\natexlab}[1]{#1}
\providecommand{\url}[1]{\texttt{#1}}
\expandafter\ifx\csname urlstyle\endcsname\relax
  \providecommand{\doi}[1]{doi: #1}\else
  \providecommand{\doi}{doi: \begingroup \urlstyle{rm}\Url}\fi

\bibitem[Adams and Fournier(2003)]{Adams:03}
H.~Adams and J.~Fournier.
\newblock \emph{Sobolev spaces}.
\newblock Academic Press, 2003.

\bibitem[Amat et~al.(2018)Amat, Michalski, and Stoltz]{Amat:18}
C.~Amat, T.~Michalski, and G.~Stoltz.
\newblock Fundamentals and exchange rate forecastability with simple machine
  learning methods.
\newblock \emph{Journal of International Money and Finance}, 88:\penalty0
  1--24, 2018.

\bibitem[Azoury and Warmuth(2001)]{Azoury:01}
K.~Azoury and M.~Warmuth.
\newblock Relative loss bounds for on-line density estimation with the
  exponential family of distributions.
\newblock \emph{Machine learning}, 43:\penalty0 211--246, 2001.

\bibitem[Blanchard and Muecke(2017)]{Blanchard:17}
G.~Blanchard and N.~Muecke.
\newblock Optimal rates of regularization of statistical inverse learning
  problems.
\newblock \emph{Foundations of Computational Mathematics}, 18:\penalty0
  971--1013, August 2017.

\bibitem[Blanchard et~al.(2007)Blanchard, Bousquet, and Zwald]{Blanchard:07}
G.~Blanchard, O.~Bousquet, and L.~Zwald.
\newblock Statistical properties of kernel principal component analysis.hal
  hal-00373789.
\newblock \emph{Machine Learning}, 3:\penalty0 259--294, 2007.

\bibitem[Brezis and Mironescu(2018)]{Brezis:17}
H.~Brezis and P.~Mironescu.
\newblock Gagliardo-nierenberg inequalities and non-inequalities.
\newblock \emph{Annales de l'Institut de Henri Poincare}, 1:\penalty0
  1355--1376, 2018.

\bibitem[Caponetto and E.De.Vito(2006)]{Caponentto:06}
A.~Caponetto and E.De.Vito.
\newblock Optimal rates for the regularized least-squares algorithm.
\newblock \emph{Foundations of Computational Mathematics}, pages 331--368,
  2006.

\bibitem[Cesa-Bianchi(1999)]{Bianchi:99}
N.~Cesa-Bianchi.
\newblock Analysis of two gradient-based algorithms for online regression.
\newblock \emph{Journal Computational System Sci.}, pages 392--411, 1999.

\bibitem[Cesa-Bianchi and Lugosi(2006)]{Bianchi:06}
N.~Cesa-Bianchi and G.~Lugosi.
\newblock \emph{Prediction, Learning and Games}.
\newblock Cambridge University Press, 2006.

\bibitem[Cesa-Bianchi et~al.(2017)Cesa-Bianchi, Gentile, and
  Gerchinovitz]{cesa2017algorithmic}
P.~Cesa-Bianchi, N.and~Gaillard, C.~Gentile, and S.~Gerchinovitz.
\newblock Algorithmic chaining and the role of partial feedback in online
  nonparametric learning.
\newblock \emph{arXiv preprint arXiv:1702.08211}, 2017.

\bibitem[Clarkson(1936)]{Clarkson:36}
J.A. Clarkson.
\newblock Uniformly convex spaces.
\newblock \emph{Transactions of the American Mathematical Society},
  40:\penalty0 396--414, 1936.

\bibitem[Devaine et~al.(2013)Devaine, Gaillard, Goude, and Stoltz]{Devaine:13}
M.~Devaine, P.~Gaillard, Y.~Goude, and G.~Stoltz.
\newblock Forecasting electricity consumption by aggregating specialized
  experts - a review of the sequential aggregation of specialized experts, with
  an application to slovakian and french country-wide one-day-ahead
  (half-)hourly predictions.
\newblock \emph{Machine Learning}, 90\penalty0 (2):\penalty0 231--260, 2013.

\bibitem[Di~Nezza et~al.(2012)Di~Nezza, Palatucci, and Valdinoci]{Nezza:12}
E.~Di~Nezza, G.~Palatucci, and E.~Valdinoci.
\newblock Hitchhiker's guide to the fractjional sobolev spaces.
\newblock \emph{Bulletin des Sciences Mathematique}, 136:\penalty0 521--573,
  2012.

\bibitem[Eberts and Steinwart(2011)]{Eberts:11}
M.~Eberts and I.~Steinwart.
\newblock Optimal learning rates for least squares svm using gaussian kernels.
\newblock In \emph{Advances in Neural Information Processing Systems 24}, pages
  1539--1547. Curran Associates, Inc., 2011.

\bibitem[Edmunds and Triebel(1996)]{Edmunds:96}
D.~Edmunds and H.~Triebel.
\newblock \emph{Function Spaces, Entropy Numbers,Differential Operators}.
\newblock Cambridge University Press, 1996.

\bibitem[Engl et~al.(2000)Engl, Hanke, and Neubauer]{Engl:00}
H.W. Engl, M.~Hanke, and A.~Neubauer.
\newblock \emph{Regularization of inverse problems}.
\newblock Springer Netherlands, 2000.
\newblock ISBN 978-0-7923-4157-4.

\bibitem[Evans(1998)]{Evans:98}
L.C. Evans.
\newblock \emph{Partial Differential Equations}.
\newblock American Mathematical Society, 1998.

\bibitem[Fischer and Steinwart(2017)]{Fischer:17}
S.~Fischer and I.~Steinwart.
\newblock Sobolev norm learning rates for regularized least-squares algorithms.
\newblock \emph{Arxiv}, pages 1--26, 2017.
\newblock URL \url{https://arxiv.org/pdf/1702.07254.pdf}.

\bibitem[Foster(1991)]{Foster:91}
D.~Foster.
\newblock Prediction in the worst case.
\newblock \emph{Annals of Statistics}, 19:\penalty0 1084--1090, 1991.

\bibitem[Gaillard and Gerchinovitz(2015)]{Gaillard:15}
P.~Gaillard and S.~Gerchinovitz.
\newblock A chaining algorithm for online nonparametric regression.
\newblock In \emph{Proceedings of The 28th Conference on Learning Theory},
  volume~40, pages 764--796, 2015.

\bibitem[Gaillard et~al.(2019)Gaillard, Gerchinovitz, Huard, and
  Stoltz]{Gaillard:19}
P.~Gaillard, S.~Gerchinovitz, M.~Huard, and G.~Stoltz.
\newblock Uniform regret bounds over $\mathbb{R}^d$ for the sequential linear
  regression problem with the square loss.
\newblock In Aur\'elien Garivier and Satyen Kale, editors, \emph{Proceedings of
  the 30th International Conference on Algorithmic Learning Theory}, volume~98
  of \emph{Proceedings of Machine Learning Research}, pages 404--432, Chicago,
  Illinois, 22--24 Mar 2019. PMLR.
\newblock URL \url{http://proceedings.mlr.press/v98/gaillard19a.html}.

\bibitem[Gammerman et~al.(2004)Gammerman, Kalnishkan, and Vovk]{Gammerman2004}
A.~Gammerman, Y.~Kalnishkan, and V.~Vovk.
\newblock On-line prediction with kernels and the complexity approximation
  principle.
\newblock In \emph{Proceedings of the 20th conference on Uncertainty in
  artificial intelligence}, pages 170--176, 2004.

\bibitem[Gy\"orfi(2002)]{Gyorfi:02}
L.~Gy\"orfi.
\newblock \emph{A Distribution-Free theory of nonparametric regression}.
\newblock Springer, 2002.

\bibitem[Hazan et~al.(2018)Hazan, Hu, Li, and Li]{Hazan:18}
E.~Hazan, W.~Hu, Y.~Li, and Z.~Li.
\newblock Online improper learning with an approximation oracle.
\newblock In S.~Bengio, H.~Wallach, H.~Larochelle, K.~Grauman, N.~Cesa-Bianchi,
  and R.~Garnett, editors, \emph{Advances in Neural Information Processing
  Systems}, volume~31, pages 5652--5660. Curran Associates, Inc., 2018.
\newblock URL
  \url{https://proceedings.neurips.cc/paper/2018/file/ad47a008a2f806aa6eb1b53852cd8b37-Paper.pdf}.

\bibitem[J{\'e}z{\'e}quel et~al.(2019)J{\'e}z{\'e}quel, Gaillard, and
  Rudi]{jezequel2019efficient}
R.~J{\'e}z{\'e}quel, P.~Gaillard, and A.~Rudi.
\newblock Efficient online learning with kernels for adversarial large scale
  problems.
\newblock In \emph{Advances in Neural Information Processing Systems}, pages
  9427--9436, 2019.

\bibitem[Lin and Cevher(2018)]{Lin:18}
J.~Lin and V.~Cevher.
\newblock Optimal convergence for distributed learning with stochastic gradient
  methods and spectral regularization algorithms.
\newblock \emph{Arxiv}, pages 1--53, 2018.
\newblock URL \url{https://arxiv.org/pdf/1801.07226.pdf}.

\bibitem[Lin et~al.(2020)Lin, Rudi, Rosasco, and V.]{Lin:20}
J.~Lin, A.~Rudi, L.~Rosasco, and Cevher V.
\newblock Optimal rates for spectral algorithms with least-squares regression
  over hilbert spaces.
\newblock \emph{Applied and Computational Harmonic Analysis}, pages 868--890,
  2020.
\newblock URL
  \url{https://www.sciencedirect.com/science/article/abs/pii/S1063520318300174}.

\bibitem[Loring(2011)]{WTu:11}
W.~Tu Loring.
\newblock \emph{An Introduction to Manifolds}.
\newblock Springer, 2011.

\bibitem[Mallet et~al.(2009)Mallet, Stoltz, and Mauricette]{Mallet:09}
V.~Mallet, G.~Stoltz, and B.~Mauricette.
\newblock Ozone ensemble forecast with machine learning algorithms.
\newblock \emph{Journal of Geophysical Research: Atmospheres}, 114\penalty0
  (D5), 2009.

\bibitem[Narcowich and Ward(2004)]{Narcowich:04b}
F.~Narcowich and J.~Ward.
\newblock Scattered-data interpolation on $\mbr^{n}$: error estimates for
  radial basis and band-limited functions.
\newblock \emph{SIAM J. MATH. ANAL}, 36:\penalty0 284--300, 2004.

\bibitem[Narcowich et~al.(2004)Narcowich, Ward, and Wendland]{Narcowich:04}
F.~Narcowich, J.~Ward, and H.~Wendland.
\newblock Sobolev bounds on functions with scattered zeros, with applications
  to radial basis function surface fitting.
\newblock \emph{Mathematics of Computation}, 74:\penalty0 743--763, 2004.

\bibitem[Novak et~al.(2017)Novak, Ulrich, Wozniakowski, and Zhung]{Novak:17}
E.~Novak, M.~Ulrich, H.~Wozniakowski, and S.~Zhung.
\newblock Reproducing kernels of sobolev spaces on $\mbr^d$ and applications to
  embedding constants and tractability.
\newblock \emph{Arxiv}, 2017.
\newblock URL \url{https://arxiv.org/pdf/1709.02568.pdf}.

\bibitem[Pagliana et~al.(2020)Pagliana, Rudi, De~Vito, and
  Rosasco]{Pagliana:20}
N.~Pagliana, A.~Rudi, E.~De~Vito, and L.~Rosasco.
\newblock Interpolation and learning with scale-dependent kernels.
\newblock \emph{Arxiv}, 2020.
\newblock URL \url{https://arxiv.org/pdf/2006.09984.pdf}.

\bibitem[Pillaud-Vivien et~al.(2018)Pillaud-Vivien, Rudi, and Bach]{Vivien:18}
L.~Pillaud-Vivien, A.~Rudi, and F.~Bach.
\newblock Statistical optimality of stochastic gradient descent on hard
  learning problems through multiple passes.
\newblock In \emph{NIPS}, 2018.

\bibitem[Rakhlin and Sridharan(2014)]{Rakhlin:14}
A.~Rakhlin and K.~Sridharan.
\newblock Online nonparametric regression.
\newblock \emph{Journal of Machine Learning Research}, pages 1--27, 2014.

\bibitem[Rakhlin et~al.(2014)Rakhlin, Sridharan, and A.Tewari]{Rakhlin:14a}
A.~Rakhlin, K.~Sridharan, and A.Tewari.
\newblock Sequential complexities and uniform martingale laws of large numbers.
\newblock \emph{Probability Theory and related random fields}, 161:\penalty0
  111--153, 2014.

\bibitem[Rakhlin et~al.(2015)Rakhlin, Sridharan, and Tewari]{Rakhlin:15}
A.~Rakhlin, K.~Sridharan, and A.~Tewari.
\newblock Online learning via sequential complexities.
\newblock \emph{Journal of Machine Learning Research}, pages 155--186, 2015.

\bibitem[Rudi et~al.(2015)Rudi, Camoriano, and Rosasco]{Rudi:15}
A.~Rudi, R.~Camoriano, and L.~Rosasco.
\newblock Less is more: Nystr{\"o}m computational regularization.
\newblock In \emph{Advances in Neural Information Processing Systems}, pages
  1657--1665, 2015.

\bibitem[Schaback(2007)]{Schaback:07}
J.~Schaback.
\newblock \emph{Kernel-based meshless methods}.
\newblock Lecture notes, 2007.

\bibitem[Smola and Sch\"olkopf(2002)]{Smola:02}
A.~Smola and B.~Sch\"olkopf.
\newblock \emph{Learning with Kernels: Support Vector Machines, Regularization,
  Optimization and Beyond}.
\newblock MIT Press, Cambridge, MA, 2002.

\bibitem[Stein(1970)]{Stein:70}
E.M. Stein.
\newblock \emph{Singular integrals and differentiability properties of
  functions}.
\newblock Princeton University Press, 1970.

\bibitem[Steinwart and Cristmann(2008)]{Steinwart:08}
I.~Steinwart and A.~Cristmann.
\newblock \emph{Support Vector Machines}.
\newblock Springer, 2008.

\bibitem[Steinwart et~al.(2019)Steinwart, Hush, and Scovel]{Steinwart:09}
I.~Steinwart, D.~Hush, and C.~Scovel.
\newblock Optimal rates for least-squares regression.
\newblock In S.~Dasgupta and A.~Klivans, editors, \emph{Proceedings of the 22nd
  Annual Conference on Learning Theory}, pages 79--93, 2019.

\bibitem[Tsybakov(2009)]{Tsybakov:09}
A.~Tsybakov.
\newblock \emph{Introduction to nonparametric estimation}.
\newblock Springer, 2009.

\bibitem[Vovk(1998)]{Vovk:97}
V.~Vovk.
\newblock Competitive online linear regression.
\newblock \emph{Proceedings of the 1997 conference on advances in neural
  information processing systems, 10}, pages 364--370, 1998.

\bibitem[Vovk(2001)]{Vovk:01}
V.~Vovk.
\newblock Competitive online statistics.
\newblock \emph{International statistical review}, 69:\penalty0 213--248, 2001.

\bibitem[Vovk(2006{\natexlab{a}})]{Vovk06Metric}
V.~Vovk.
\newblock Metric entropy in competitive online prediction.
\newblock \emph{Arxiv}, 2006{\natexlab{a}}.

\bibitem[Vovk(2006{\natexlab{b}})]{Vovk:06}
V.~Vovk.
\newblock On-line regression competitive with reproducing kernel hilbert
  spaces.
\newblock In \emph{International Conference of Theory and Application of Models
  of Computation}, volume~69, pages 452--463, 2006{\natexlab{b}}.

\bibitem[Vovk(2007)]{Vovk:07}
V.~Vovk.
\newblock Competing with wild prediction rules.
\newblock \emph{Machine Learning}, 69:\penalty0 193--212, 2007.

\bibitem[Wendlandt(2005)]{Wendlandt:05}
H.~Wendlandt.
\newblock \emph{Scattered Data Approximation}.
\newblock Cambridge University Press, 2005.

\bibitem[Zhang(2005)]{Zhang:05}
T.~Zhang.
\newblock Learning bounds for kernel regression using effective data
  dimensionality.
\newblock \emph{Neural Computation 17(9)}, pages 2077--2098, 2005.

\bibitem[Zhdanov and Kalnishkan(2010)]{Zhdanov:10}
F.~Zhdanov and Y.~Kalnishkan.
\newblock An identity for kernel ridge regression.
\newblock In \emph{Algorithmic Learning Theory}, pages 405--419. Springer,
  2010.

\end{thebibliography}

\newpage
\begin{center}
    \Large
    Appendices
    \vspace*{1cm}
\end{center}
\appendix 

\section{Notation on kernels and linear operators over reproducing kernel Hilbert spaces}
\label{app:kernels}
We complete Section~\ref{sec:kern_lrn} by providing addition notations on kernels that are used in the proofs. 
We consider kernel methods that choose forecaster $\smash{\hat f_t}$ in a reproducing kernel Hilbert space $\cH_{k}$ which is associated with a reproducing kernel $k: \cX \times \cX \mapsto \mbr$. The prediction rule $f_{t}$ at round then forecasts $\hat{f}_{t}\paren{x_{t}} = \inner{f_{t}, k_{x_{t}}}_{\cH_{k}}$. We use the following notations, which are common in the setting of kernel learning.
\paragraph{Integral and covariance operators}  Let $\paren{\cX,\cB\paren{\cX}}$ be a measurable space and $\mu$ be some measure on a Borel $\sigma-$algebra $\cB\paren{\cX}$. We define $S: \cH_{k} \mapsto L_2\paren{\cX,\mu}$ to be the restriction operator of a function $f \in \cH_k$ to its equivalence class in $L_2\paren{\cX,\mu}$. We drop the dependence of $S$ on the measure $\mu$ to simplify the notation. The correspondent adjoint $S^{\star}: L_{2}\paren{\cX,\mu} \mapsto \cH_{k}$ is then well-defined and has the form $S^{\star}f = \int_{x \in \cX} f\paren{x}k_{x}d\mu\paren{x}$ for any $f \in L_{2}\paren{\cX,\mu}$. We define the (kernel) integral operator $L = SS^{\star} : L_2\paren{\cX,\mu} \mapsto L_2\paren{\cX,\mu}$ such that for any $f \in L_2\paren{\cX,\mu}$ we have for  $\mu$ almost all $x \in \cX$
\begin{align}
\label{eq:kern_int}
    L\paren{f}\paren{x} = \int_{z \in \cX} k\paren{x,z}f\paren{z}d\mu\paren{z}.
\end{align}
The (kernel) covariance operator $T = S^{\star}S: \cH_k \mapsto \cH_k$ is defined as
\begin{align}
    T = \int_{x \in \cX }k_{x} \otimes k_xd\mu\paren{x}.
\end{align}
It is known (see e.g. Theorem 2.2 and Theorem 2.3 in \cite{Blanchard:07}) that the operators $T$ and $L$ are both positive, self-adjoint and trace-class operators. Moreover, they have the same non-zero spectrum.

\paragraph{Evaluation and empirical covariance operators} Analogous to the population case, based on the data sequence $(x_{s},y_{s})_{1\leq s\leq t}$ for each $t\in \{1,\ldots,n\} $, we define the evaluation operator $S_{t}f : \cH_{k} \mapsto \mathbb{R}^{t}$, such that for any $j \in \{1,\ldots, t\}$  
\[
\paren{S_{t}f}_{j} = \inner{f,k_{x_{j}}} = f\paren{x_{j}} \,.
\] 
Let $S_{t}^{\star} : \mbr^{t} \mapsto \cH_{k}$ be the corresponding adjoint. Then, for any $y \in \mbr^{t}$  
\[
S_{t}^{\star} y = \sum_{i=1}^{t}y_{i}k_{x_{i}} \,.
\]
Note that the kernel matrix $K_{t} \bydef \paren{k(x_i,x_j)}_{1\leq i,j \leq t}$ satisfies $K_t=S_{t}S_{t}^{\star}$. We also define $T_{t}: \cH_{k} \mapsto \cH_{k}$ the empirical covariance operator for $t\geq 1$ as
\[
T_{t} \bydef S_{t}^{\star}S_{t} = \sum_{i=1}^{t} k_{x_{i}} \otimes k_{x_{i}} \,.
\]
For any $f \in \cH_{k}$, $T_{t}f = \sum_{i=1}^{t} k_{x_{i}}\inner{k_{x_{i}},f} = \sum_{i=1}^{t}f\paren{x_i}k_{x_{i}}$. For a given $\tau>0$, we define the regularized covariance operator $A_{t}= T_{t}+\tau \mbi$, where $\mbi: \cH_{k} \mapsto \cH_{k}$ is the identity operator. Finally, we call $\lambda_{j}\paren{A}$ the $j-$th largest eigenvalue of the operator $A$ (i.e. $\lambda_{1}\paren{A} \geq \lambda_{2}\paren{A} \ldots \geq \lambda_{n}\paren{A} \geq \ldots$). It is worth pointing out that both $T_{t}$ and $K_{t}$ are positive semi-definite for all $t \in \{1,2,\ldots,n\}$. Since the kernel $k\paren{\cdot,\cdot}$ is bounded, $T_{t}$ is a trace class operator. In other words, $T_t$ is a compact operator for which a trace may be defined; i.e., in some orthonormal basis 
$\smash{\paren{\phi_{k}}_{k \in \mbn_*}}$, the trace 
$\smash{\norm{T_t}_{1} \bydef \tr{\abs{A}} \bydef \sum_{k} \big\langle (T_t^{\star}T_t)^{1/2} \phi_{k},\phi_{k} \big\rangle = \sum_{k} \sqrt{\lambda_{k}\paren{T_t^{\star}T_t}}}$
is finite. With a slight abuse of notation, we write $S_{n}f = \paren{f\paren{x_i}}_{i=1}^{n}$ for any function $f: \cX \mapsto \mbr$ and datasample $\{x_{s}\}_{s \geq 1}$.

\section{Preliminary results on Sobolev spaces}
\label{app:sob_spaces}

 In this part, we recall known results on Sobolev spaces that will be useful for our analysis. We refer the curious reader to \cite{Adams:03} for an extensive survey on Sobolev spaces and to \cite{Nezza:12} for the specific case of non-integer exponents.  
 
 \subsection{Definition and notation}
 \label{app:def_sobolev}
 Let here $\cX \subseteq \mbr^d$, $p \in [1,\infty)$ and denote  $L_{p}\paren{\cX}$ for the equivalence class of $p-$integrable functions with respect to the Lebesque measure $\lambda$ on $\cX$. We recall the definition of Sobolev spaces $W_p^r(\cX)$ when $r \geq 0$ is an integer.
 
 \paragraph{Definition of Sobolev spaces with integer $r \in \mbn_*$.}  We recall (see Section~\ref{sec:sob_spaces}) that the Sobolev spaces $W_{p}^{r}\paren{\cX}$ and $W_{\infty}^{r}\paren{\cX}$ are the vector spaces of equivalence classes of functions defined as:
 
\begin{align*}
    W_{p}^{r}\paren{\cX} \bydef \bigg \{f:\cX \to \mbr  \quad \text{s.t.} \quad  \norm{f}_{W^{r}_{p}\paren{\cX}} & \bydef \paren[1]{ \sum_{\abs{\gamma}_{1} \leq r}\norm{D^{\gamma}f}_{L_{p}\paren{\cX}}^{p}}^{\frac{1}{p}} < \infty \bigg \},
\end{align*}
and
\begin{align*}
    W_{\infty}^{r}\paren{\cX} \bydef \bigg \{f:\cX \to \mbr \quad \text{s.t.} \quad   \norm{f}_{W_{\infty}^{r}\paren{\cX}} & \bydef \sup_{\abs{\gamma}_{1}\leq r} \norm{D^{\gamma}f}_{L_{\infty}\paren{\cX}} < \infty \bigg \} \,.
\end{align*}
We also define the Sobolev semi-norm $\smash{\abs{f}_{W_{p}^{j}\paren{\cX}} \bydef \sum_{\gamma: \abs{\gamma}=j}\norm{D^{\gamma}f}_{L_{p}\paren{\cX}}}$. 

\paragraph{Definition of Sobolev spaces with non-integer smoothness exponent $\beta$.} Let $\beta \in \mbr_{+}$; for our proposes we write $\beta = r + \sigma$ with $r \in \mbn_0$ and $\sigma \in (0,1)$. 
Let $u : \cX \mapsto \mbr$ be some fixed measurable function. We define the map $\varphi_{u} : \cX \times \cX \mapsto \mbr\cup\{ \infty\} $ such that for $1 \leq p < \infty$ and all $(x,y) \in \cX \times \cX$:
\[
\varphi_{u}\paren{x,y} = \frac{\abs{u\paren{x} - u\paren{y} }}{\norm{x-y}_{2}^{\frac{d}{p}+\sigma }} \,,
\]
and denote 
\begin{align*}
    \tilde{W}_{p}^{\sigma}\paren{\cX} \bydef \{ u \in L_{p}\paren{\cX} : \norm{\varphi_{u}}_{L_{p}\paren{\cX \times \cX}} < \infty\}.
\end{align*}
The space $\tilde{W}_{p}^{\sigma}\paren{\cX}$ equipped with the norm $\norm{u}_{\tilde{W}_{p}^{\sigma}\paren{\cX}} \bydef \big( \norm{u}_{ L_{p}\paren{\cX} } + \norm{ \varphi_{u} }_{L_{p}\paren{ \cX \times \cX }} \big)^{\frac{1}{p}}$ can be shown to be a Banach space. 
With this notation, Sobolev space $W^{\beta}_{p}\paren{\cX}$, $\beta = r + \sigma$ can be defined as
\begin{equation}
\label{eq:gen_sob}
    W^{\beta}_{p}\paren{\cX} \bydef \big\{ u \in W_{p}^{r}\paren{\cX}: D^{\gamma}u \in \tilde{W}_{p}^{\sigma}\paren{\cX} \text{ for any $\gamma \in \mbn^{d}$ such that $\abs{\gamma}_1 = r$ }\big\} \,.
\end{equation}
Equipped with the norm
\begin{align}
    \label{eq:normreal}
    \norm{u}_{W_{p}^{\beta}\paren{\cX}} \bydef \bigg( \norm{u}^{p}_{W_{p}^{r}\paren{\cX}} + \sum_{\gamma: \abs{\gamma} = r} \norm{D^{\gamma}u}^{p}_{\tilde{W}_{p}^{\sigma}\paren{\cX}} \bigg)^{\frac{1}{p}} \,,
\end{align}
it becomes Banach space. In the case $\beta=m \in \mbn_*$, it matches the definition of the Sobolev space $W_{p}^{m}\paren{\cX}$ (up to a re-scaling of the norm). If $m=0$ (i.e. $r = \sigma \in [0,1)$), we find that $W_{p}^{m}\paren{\cX} = L_{p}\paren{\cX}$ so that the norm in $W_{p}^{\sigma}\paren{\cX}$ is given by
\begin{align}
    \label{eq:normless1}
    \norm{u}_{W_{p}^{r}\paren{\cX}}= \norm{u}_{\tilde{W}_{p}^{\sigma}\paren{\cX}} \bydef \big( \norm{u}_{L_{p}\paren{\cX}} + \norm{\varphi_{u}}_{L_{p}\paren{\cX \times \cX}} \big)^{\frac{1}{p}} \,.
\end{align}
In accordance with above definition of the class $W_{p}^{\beta}\paren{\cX}$, for any $\beta = r +\sigma$, $\sigma \in [0,1)$ we set 
\begin{align}
    \tilde{W}_{\infty}^{\sigma}\paren{\cX} \bydef \{ u \in L_{\infty}\paren{\cX}: \sup_{x,y \in \cX, x\neq y}\frac{\abs{u\paren{x} - u\paren{y}}}{\norm{x-y}_{2}^{\sigma}} \leq \infty \}. 
\end{align}
Now, for $\beta =r+\sigma \in \mbr$ the Sobolev space $W_{\infty}^{\beta}\paren{\cX}$ can be defined as a functional space
\begin{equation}
\label{eq:gen_sob_infty}
    W^{\beta}_{\infty}\paren{\cX} \bydef \big\{ u \in W_{\infty}^{m}\paren{\cX}: D^{\gamma}u \in \tilde{W}_{\infty}^{\sigma}\paren{\cX} \text{ for any $\gamma$ such that $\abs{\gamma}_1 = r$ }\big\} \,.
\end{equation}
equipped with a norm
\begin{align}
    \norm{u}_{W_{\infty}^{\beta}\paren{\cX}} \bydef \max\{ \norm{u}_{W_{\infty}^{r}\paren{\cX}}, \max_{\gamma: \abs{\gamma}_{1}=r}\norm{D^{r}u}_{\tilde{W}_{\infty}^{\sigma}\paren{\cX}} \}
\end{align}

\subsection{Approximation properties of the Sobolev spaces. }

We recall that $W^s(\cX)$ is a Sobolev RKHS, a space of continuous representatives from equivalence classes of functions from the Sobolev space $W_2^s(\cX)$ provided $s>\frac{d}{2}$. The goal of this section is to control the regret with respect to a ball in an arbitrary Sobolev space $\smash{W_p^\beta}\paren{\cX}$ with $p \geq 2$ and $\beta \neq s$. To do so, we need to control the approximation error of $ f \in \smash{W_p^\beta}\paren{\cX}$ by the elements from some subset $\cG \subset W_2^s\paren{\cX}$ uniformly over $f \in \smash{W_p^\beta}\paren{\cX} $. This can be achieved by considering the subset of the band limited functions (see ex. \cite{Narcowich:04}), which is in $W_{2}^{s}\paren{\cX}$ for any $s>0$. Namely, for $\sigma \in \mbr_{+}\setminus \{0\}$ we define $B_{\sigma}$ to be 
\begin{align}
    B_{\sigma} \bydef \{ f \in L_{2}\paren[1]{\mbr^{d}}\cap C_{\infty}\paren[1]{\mbr^{d}} : supp\paren{\cF(f)} \subset B\paren{0,\sigma} \},
\end{align}
where we denote $\cF\paren{f}$ for the Fourier transform of $f$ and recall that $B\paren{0,\sigma}$ is an open ball in $\mbr^{d}$ with radius $\sigma$. 

The next result is the consequence of Proposition 3.7 in \cite{Narcowich:04b}
(see also the proof of Lemma 3.7 in \cite{Narcowich:04}).
To be able to apply the aforementioned Proposition we need to extend functions $f : \cX \mapsto \mbr$, $f \in W_{}^{s}\paren{\cX}$ to functions $\tilde{f} :\mbr^{d} \mapsto \mbr$  such that $\tilde{f} \in W_{}^{s}\paren{\mbr^d}$. By Stein's Extension Theorem (see \cite{Stein:70}, page. 181) because $\cX$ is a bounded Lipschitz domain there exists a linear operator $\mathfrak{C} : W_{}^{s}\paren{\cX} \mapsto W_{}^{s}\paren{\mbr^d}$ 
which is continuous ( i.e. since it is linear we have $\norm{\mathfrak{C} f}_{W_{2}^{s}\paren{\mbr^d}} \leq \tilde{C} \norm{f}_{W_{2}^{s}\paren{\cX}}$). For this operator $\cC$, every $f \in W_{}^{s}\paren{\cX}$ and $g_{\sigma} \in B_{\sigma}$ by definition of the norm in $W^{s}_{2}\paren{\cX}$ we have $\norm{f - g_{\sigma}}_{W_{2}^{s}\paren{\cX}} \leq \norm{\mathfrak{C}f - g_{\sigma}}_{W_{2}^{s}\paren{\mbr^d}}$. Applying Lemma 3.7 in \cite{Narcowich:04} to $\mathfrak{C}f\in W_{}^{s}\paren{\mbr^d}$, and using the argument as in the proof of Theorem 3.8 in \cite{Narcowich:04} for $g_{\sigma}$ given by Lemma 3.7, we have 

$$ \norm{f-g_{\sigma}}_{W_{2}^{r}\paren{\mbr^d} } \leq c\sigma^{r-s}\norm{g_{\sigma}}_{W_{2}^{s}\paren{\mbr^d} } $$
and 
$$\norm{g_{\sigma}}_{W_{2}^{s}\paren{\cX}} \leq \norm{g_{\sigma}}_{W_{2}^{s}\paren{\mbr^d}} \leq c_2\norm{\mathfrak{C} f}_{W_{2}^{s}\paren{\mbr^d}} \leq c_3\norm{ f}_{W_{2}^{s}\paren{\cX}}.$$ 
Thus we obtain the following statement. 

\begin{proposition}
	
	\label{prop:sig_approx_sobolev}
	Let $s \geq r \geq 0$. For every $f \in W^{s}_{2}\paren{\cX}$, $\sigma >0$ there exists a function $g_{\sigma} \in B_{\sigma}$ and constants $C_0$ and $C_1$  
	which are independent of $\sigma$
	such that 
	\[
	\norm{f-g_{\sigma}}_{W_{2}^{r}\paren{\cX}} \leq C_{0} \sigma^{r-s} \norm{f}_{W_{2}^{s}\paren{\cX}} 
	\quad 
	\text{and}
	\quad 
	\norm{g_{\sigma}}_{W_{2}^{r}\paren{\cX}} \leq C_{1} \sigma^{r-s}\norm{f}_{W_{2}^{s}\paren{\cX}} \,.
	\]
\end{proposition}

We now state an upper-bound of $\smash{\norm{f}_{W_{p}^{r}\paren{\cX}}}$ when $f$ belongs to the  intermediate Sobolev spaces $\smash{W_{p_1}^{s_1}(\cX)}$ and $\smash{W_{p_2}^{s_2}(\cX)}$ for some $p_1,p_2,s_1,s_2$. This result is a Gagliardo-Nirenberg--type inequality and follows from the result originally stated in Theorem 1 in \cite{Brezis:17}.

\begin{proposition}[Theorem 1, \cite{Brezis:17}]
	\label{thm:gagl_nierenberg_ineq}
	Let $\cX \subseteq \mbr^{d}$ be a Lipschitz bounded domain. Let $0 \leq r,s_{1},s_{2} < \infty$ and $1\leq p_{1},p_{2},p \leq \infty $ be real numbers such that there exists $\theta \in (0,1)$ with
	\[
	r = \theta s_{1} + \paren{1-\theta}s_{2} \quad \text{ and } \quad \frac{1}{p} = \frac{\theta}{p_1} + \frac{1-\theta}{p_2} \,.
	\]
	Let
	$
	A \bydef \big\{ \paren{s_{1},s_{2},p_{1},p_{2}} \quad \text{ s.t. } \quad  s_2 \in \mbn_{*},\  p_2 =1,\  s_{2}-s_{1} \leq 1 - \frac{1}{p_1} \big\}.
	$
	If $\paren{s_{1},s_2,p_1,p_2} \notin A$, then there exists a constant $C >0$ which depends on $s_1, s_2, p_1, p_2, \theta$ and $\cX$ such that 
	\[
	\norm{f}_{W_{p}^{r}\paren{\cX}} \leq  C \norm{f}^{\theta}_{W_{p_{1}}^{s_{1}}\paren{\cX}} \norm{f}^{1-\theta}_{W_{p_{2}}^{s_{2}}\paren{\cX}}\,,
	\]
	for all $f \in W_{p_1}^{s_{1}}\paren{\cX} \cap W_{p_2}^{s_{2}}\paren{\cX}$.
\end{proposition}

In the next corollary we state two particular cases of Proposition~\ref{thm:gagl_nierenberg_ineq} that will prove useful. 
\begin{corollary}
	For the domain $\cX = [-1,1]^{d}$ and any $\epsilon >0$, all $p \geq 2$ and $\beta > d/p $ there exists a constant $C >0$ depending on $p$, $d$, $\epsilon$ and $\beta$ such that
	\begin{equation}
	\label{eq:spec_case}
	\norm{g}_{W_{p}^{\frac{d}{p} + \epsilon}\paren{\cX}} \leq C \norm{g}^{\frac{d}{\beta p}+\frac{\epsilon}{\beta}}_{W_{p}^{\beta}\paren{\cX}} \norm{g}^{1-\frac{d}{\beta p}-\frac{\epsilon}{\beta}}_{L_{p}\paren{\cX}}\,,
	\end{equation}
	for all function $g \in \smash{W_{p}^{\beta}\paren{\cX}}$. Furthermore, for all $\beta >0$,  $p\geq 2$ and $\epsilon >0$, there exists a constant $C > 0$ depending on $\beta$, $p$, $d$, and $\epsilon$ such that
	\begin{equation}
	\label{eq:spec_case_p}
	\norm{g}_{W_{2}^{\frac{d}{2} + \epsilon}\paren{\cX}} \leq C \norm{g}^{ {\frac{d+2\epsilon}{{\beta p}}}}_{W_{2}^{\beta p/2}\paren{\cX}} \norm{g}^{1- \frac{d+2\epsilon}{{\beta p}}}_{L_{2}\paren{\cX}} \,,
	\end{equation}
	for any function $g \in W_{2}^{\beta}\paren{\cX}$.
\end{corollary}

\begin{proof}
	First, notice that $\cX = [-1,1]^{d}$ is a Lipschitz bounded domain. The first inequality is obtained by choosing $p_1=p_2=p\geq 1$, $r = d/p+\epsilon$, $s_1 = \beta$, and $s_2 = 0$ in Proposition~\ref{thm:gagl_nierenberg_ineq}; checking that $\paren{s_{1},s_{2},p_1,p_2} \notin A$; and noting that for any $\beta>0$ we have $W_{2}^{\beta}\paren{\cX} \cap L_{2}\paren{\cX} = W_{2}^{\beta}\paren{\cX}$.
	The second inequality stems from the choice $p = p_1 = p_2 = 2$ (note that this is for the $p$ in the Proposition which is different from the $p$ in the inequality), $s_{2}=0$, $s_{1} = \frac{\beta p}{2}$ and noting the inclusion  $\smash{W_{2}^{\beta}\paren{\cX} \subset W_{2}^{\beta p}\paren{\cX} \subseteq W_{2}^{\beta p/2} \paren{\cX} = L_{2}\paren{\cX} \cap W_{2}^{\beta p/2}\paren{\cX}}$ which holds true since $p\geq 2$.
\end{proof}

\subsection{Results from interpolation theory on Sobolev spaces}

To provide a sharp upper bound on the effective dimension (Proposition~\ref{prop:eff_dim_sob}), we also need the following general interpolation result on Sobolev spaces (stated in Theorem 3.8 in \cite{Narcowich:04}). Recall (see \cite{Wendlandt:05}, p.172 ) that the fill distance 
of a set of points $\cZ \subset \cX$ is defined as  $h_{Z,\cX} := \sup_{x \in \cX} \inf_{z \in Z} \norm{x-z}_{2}$.

\begin{proposition}[Theorem 3.8 in \cite{Narcowich:04}]
	\label{prop:interp_est}
	
	Suppose $\Phi:\mbr^d \to \mbr$ to be a positive definite function such that its Fourier transform $\smash{\cF\paren{\Phi}}$ satisfies
	\begin{align}
	\label{eq:dec_cond}
	c_1 \paren[1]{1 + \norm{\omega}_{2}^{2}}^{-q} \leq \smash{\cF\paren{\Phi}}\paren{\omega} \leq c_2 \paren[1]{1 + \norm{\omega}_{2}^{2}}^{-q}
	\end{align}
	where $q\geq s\geq r\geq 0$ and $c_1,c_2$ are some constants. Assume that $\cX \subset \mbr^d$ is bounded domain, has Lipschitz boundary and satisfies the interior cone condition (see Chapter 4 in \cite{Adams:03}) with parameters $\paren{\varphi, R_{0}}$. Let $k=\lfloor q \rfloor$ and $\cZ \subset \cX$ be such that its mesh norm $h\bydef h_{\cZ,\cX}$ satisfies  
	\begin{equation}
	\label{eq:meshnorm_assumption}
	h_{\cZ,\cX} \leq k^{-2}Q\paren{\varphi}R_{0}, \quad \text{where} \quad Q\paren{\varphi} \bydef \frac{\sin\paren{\varphi}\sin\paren{\theta} }{8\paren{1+\sin\paren{\theta}}(1+\sin\paren{\varphi}) }
	\end{equation}
	and $\theta = 2 \arcsin\big({\sin\paren{\varphi}}/({4\paren{1+\sin{\varphi}}}) \big)$.
	If $f \in W_{2}^{s}\paren{\cX}$ then there exists a function $v \in \mathrm{span}\{\Phi\paren{\cdot - x_{j}}, x_{j} \in \cZ \}$ such that for every real $0\leq r \leq s$
	\begin{align}
	\label{eq:interpol_ineq}
	\norm{f - v}_{W_{2}^{r}\paren{\cX}} \leq C h_{\cZ,\cX}^{s-r}\norm{f}_{W_{2}^{s}\paren{\cX}} \,,
	\end{align}
	where $C$ is some constant independent of $h_{Z,\cX}$ and $f$.

\end{proposition}

Let us now instantiate the above Proposition to the specific cases we are interested in by choosing $\cX, \Phi, \cZ$, and $r$. Let $T \in \mbn$ be fixed; set $\cX \bydef [-1,1]^d$, $\Phi$ being the feature map of Sobolev RKHS $W_{}^{s}\paren{\mbr^d}$. In this case  (see 3.1 in \cite{Narcowich:04}) $\Phi$ satisfies decay rate from Equation~\eqref{eq:dec_cond} with $q=s$. Choose $\cZ$ to be the set of points of size $T$ such that $h_{Z,\cX} \lesssim T^{-\frac{1}{d}}$ (the latter means that there exists constant $C>0$ such that $h_{Z,\cX} \leq CT^{-\frac{1}{d}}$). To control when then condition~\ref{eq:meshnorm_assumption} is fulfilled, we first notice that $\cX$ is star-shaped (see Definition~11.25 in \cite{Wendlandt:05}, also Proposition~2.1 of \cite{Narcowich:04} ); it includes $\ell_2$ ball centered at origin with radius $r=1$ and can be included in the $\ell_2$ ball centered at $0$ of radius $2\sqrt{d}$. Thus, by Proposition~2.1 in \cite{Narcowich:04}, we obtain that $\cX$ satisfies interior cone condition with the radius $R_{0}=1$ and angle $\varphi = 2\arcsin{\frac{1}{2\sqrt{d}}}$. A straightforward calculation shows that in this case
$$Q\paren{\varphi} = Q\paren{u\paren{\varphi}}= \frac{u}{8} \paren{1 - \frac{8}{8+u\sqrt{16 - u^2}}} = \paren{\frac{u}{8}}^{2} \frac{\sqrt{16-u^2}}{1 + \frac{u}{8}\sqrt{16-u^2}},$$ where $u \bydef \frac{\sin{\varphi}}{1 + \sin{\varphi}} = \frac{\sqrt{4d-1}}{2d+\sqrt{4d-1}}$. Notice that in this case we have
that $\frac{1}{8\sqrt{d}} \leq u \leq \frac{1}{2\sqrt{d}}$. 
We can easily check this by simple inequalities: 
\begin{align*}
u = \frac{\sqrt{4d-1}}{2d+\sqrt{4d-1}} \geq \frac{4d-1}{4d} \geq \frac{1}{8\sqrt{d}},
\end{align*}
and from the other side
\begin{align*}
u \leq \frac{4d-1}{\sqrt{4d-1}} = \frac{1}{\sqrt{4d-1}} \leq \frac{1}{2\sqrt{d}}.
\end{align*}
From these conditions we deduce $Q\paren{u} \geq \frac{1}{2^{12}d}$. Because $h_{\cZ,\cX}= \sup_{x \in \cX} \inf_{z \in \cZ}\norm{x-z}_{2} \lesssim T^{-\frac{1}{d}}$, to satisfy condition~\eqref{eq:meshnorm_assumption} 
we need to have $T \geq \paren{ \frac{k^2}{Q\paren{u}}}^{d}$ where we take $k = \lfloor s \rfloor $ and $R_{0}=1$. 
Notice that the choice $T \geq  \paren{4096 s^{2}d}^{d}$ ensures the last condition, therefore in order to satisfy condition~\eqref{eq:meshnorm_assumption} the size $T$ of the grid $\cZ$ should be of order $\paren{s^{2}d}^{d}$. 
Recall (see \cite{Wendlandt:05}) that the kernel $k(\cdot)$ of the Sobolev space $W_{}^{s}\paren{\mbr^d}$ can be represented by means of Bessel functions of second kind as: 
\begin{align}
\label{eq:kern_sob}
k\paren{x_1,x_2} = \frac{2^{1-s}}{\Gamma \paren{s}} \norm{x_1-x_2}_{2}^{s-\frac{d}{2}}K_{\frac{d}{2}-s}\paren{\norm{x_{1}-x_{2}}_{2}} 
\end{align}
Notice that by Corollary~10.13 in \cite{Wendlandt:05} 
the norm $\norm{\cdot}_{W^{s}\paren{\mbr^d}}$ is equivalent to $\norm{\cdot}_{W_{2}^{s}\paren{\mbr^d}}$. By Theorem 7.13 in \cite{Schaback:07} (see also Corollary 10.48 on p. 170 in \cite{Wendlandt:05} ) a restriction of RKHS $W^{s}\paren{\mbr^d}$ to the domain $\cX \bydef [-1,1]^{d}$ is itself a RKHS $W^{s}\paren{\cX}$ such that it is continuously embedded into $W_{2}^{s}\paren{\cX}$ and its kernel $k_{1}$ is a restriction of kernel $k$ to the space $\cX$. Thus, we can always consider $W^{s}\paren{\cX}$ as a RKHS with reproducing kernel $k_{1}\paren{\cdot}$ obtained by the restriction of the kernel $k\paren{\cdot}$ given by \eqref{eq:kern_sob} to the domain $\cX$. Notice that it can be written as $k_{1}\paren{x_1,x_2} = \Phi_{1}\paren{x_{1}-x_{2}}$ and since $\Phi\paren{\cdot}$ satisfies Assumption~\ref{eq:dec_cond} so also $\Phi_{1}\paren{\cdot}$.

Then, applying Proposition~\ref{prop:interp_est} twice, with $r=0$ and $r=s$ and the above choices of $\cX$, $\Phi$ and $\cZ$ entails the following corollary. 

\begin{corollary}
	\label{cor:interpolationSobolev1}
	Let $\cX \bydef [-1,1]^{d}$, $s > d/2$ and $\cZ \subset \cX^{T}$ be a set of points such that fill distance $h_{Z,\cX} \lesssim T^{-\frac{1}{d}}$, $T \geq T_{0}$, $T_{0} = \paren{4096s^{2}d}^d$. Then, for any $f \in W^s(\cX)$, there exists
	$
	\hat f \in \mathrm{span}\{k\paren{x,\cdot}, x \in Z \},
	$
	such that 
	\[
	\big\|f - \widehat{f}\big\|_{L_{2}\paren{\cX}} \leq C_1 T^{-\frac{s}{d}} \norm{f}_{W_{2}^{s}\paren{\cX} }, \qquad \big\|f - \widehat{f} \big\|_{W_{2}^{s}\paren{\cX}} \leq C_2 \norm{f}_{W_{2}^{s}\paren{\cX} } \,,
	\]
	and $\hat f(x) = f(x)$  for any $x \in \cZ$, where the constants $C_1$ and $C_2$ depend on $d$ and $s$ but are independent of the set $\cZ$ and function $f$.
\end{corollary}

The latter proposition together with Gagliardo-Nierenberg inequality yield the following approximation result of functions $f \in W^s(\cX)$ by low ranked projections $Pf$.

\begin{lemma}[Projection approximation]
	\label{lem:proj_sob}
	Let $\cX \bydef [-1,1]^{d}$, $s > d/2$, $T > T_{0}$, $T_{0}$ is given as in Lemma~\ref{cor:interpolationSobolev1} and
	$\cZ \subset \cX^{T}$ be a set of points $T$ points $\{x_1,\ldots,x_T\}$ such that the fill distance $h_{\cZ,\cX} \lesssim T^{-\frac{1}{d}}$ and $P_{\cZ}: W^{s}(\cX) \to W^{s}(\cX)$ be the orthogonal projection on $\mathrm{span}\{k_x : x \in \cZ\}$. Then, for any $f \in W^{s}(\cX)$ and for any $\epsilon_{} >0$
	\begin{align}
	\label{eq:main_proj}
	\norm{f- P_{\cZ}f}_{L_{\infty}\paren{\cX}} = \sup_{x \in \cX}\abs{f\paren{x} - (P_{\cZ}f)\paren{x}} \leq C T^{-\frac{s-\epsilon}{d}+\frac{1}{2}}\norm{f}_{W^{s}\paren{\cX}},
	\end{align}
	where $C$ is a constant independent of $f$ and $T$. 
	Furthermore, if $s \in \mbn_{*}$ then Equation~\eqref{eq:main_proj} holds with $\epsilon =0$.
\end{lemma}

\begin{proof}
	Let $f \in W^s(\cX)$ and $\epsilon > 0$. The first inequality follows from inclusion $f-P_{\cZ}f \in W^{s}\paren{\cX} \subset C\paren{\cX}$ when $s>\frac{d}{2}$. Define
	\begin{align}
	\label{eq:minimizer_subset}
	\hat{f}_{Z} \bydef \argmin_{g \in \mathrm{span}\{k_{x},x \in \cZ \} } \norm{f-g}_{W^s(\cX)}^{2}.
	\end{align}
	Because $W^{s}\paren{\cX}$ is a Hilbert space, $\hat{f}_{Z} = P_{\cZ}f \in W^{s}\paren{\cX}$. Furthermore, through reproducing property in RKHS $W^{s}\paren{\cX}$ and from the definition of an orthogonal projector, we have for any $x \in \cZ$ that $P_{\cZ}f\paren{x} = \inner{P_{\cZ}f,k_{x}} = \inner{f,P_{\cZ}k_{x}} = \inner{f,k_{x}}= f\paren{x}$. 
	By using the Sobolev embedding Theorem between the spaces $\smash{W_{2}^{{d}/{2}+\epsilon}\paren{\cX}}$ and $\smash{L_{\infty}\paren{\cX}}$ (Equation~(9) on page 60 in \cite{Edmunds:96}, applied with $s_{1}=d/2+\epsilon$, $s_{2}=0$, $n=d$, $p_1 =2$, and $p_2 =\infty$), and by using Gagliardo-Nierenberg Inequality \eqref{eq:spec_case}, we get 
	\begin{align*}
	\norm{f- P_{\cZ}f}_{L_{\infty}\paren{\cX}} 
	& \leq C_1\norm{f - P_{\cZ}f}_{W_{2}^{\frac{d}{2} +\epsilon}\paren{\cX}} \\
	& \leq C_2\norm{f-P_{\cZ}f}^{\frac{d}{2s} + \frac{\epsilon}{s}}_{W_{2}^{s}\paren{\cX}}\norm{f-P_{\cZ}f}^{1 - \frac{d}{2s} - \frac{\epsilon}{s}}_{L_{2}\paren{\cX}} \\
	& \leq  C_3 \big(\norm{f-P_{\cZ}f}_{W_{2}^{s}\paren{\cX}}\big)^{\frac{d/2 +\epsilon}{s}} T^{-\frac{s}{d}+\frac{1}{2}+\frac{\epsilon}{d}} \norm{f}^{1-\frac{d}{2s}-\frac{\epsilon}{s}}_{W_{2}^{s}\paren{\cX}}\\
	& \leq C_4 T^{-\frac{s}{d}+\frac{1}{2}+\frac{\epsilon}{d}} \norm{f}_{W_{2}^{s}\paren{\cX}} \,
	\end{align*}
	where the constants $C_1, C_2, C_3$, and $C_4$ are independent of $f$ and $T$. Finally in the specific case $s \in \mbn$ we directly apply Corollary 11.33 from \cite{Wendlandt:05} with $m=0$, $\tau=s$, $q=\infty$ to $f-P_{\cZ}f$ and obtain directly bound \eqref{eq:main_proj} with $\epsilon=0$.
\end{proof}

\section{Proof of Theorem.~\ref{prop:eff_dim_sob}. Upper bound on the effective dimension of the Sobolev RKHS}
\label{sec:d_eff_bounds}

Notice that the effective dimension can be rewritten as: 
\begin{align*}
            d_{eff}^{n}\paren{\tau} = \tr{ \paren{T_{n} +\tau \mbi}^{-1}T_{n}},
\end{align*}
where $T_{n}$ - (empirical )covariance operator.
We provide below some auxiliary results that control the tail of the trace of the kernel integral operator.  These results are provided in Lemmata 2,3 by \cite{Pagliana:20} and are just formulated here for completeness of the narrative. 

\begin{lemma}
\label{lm:supx_Linfty} 
Let $\cH_k$ be some RKHS over domain $\cX \subseteq \mbr^d$ with continuous reproducing kernel $k: \cX \times \cX \to \mbr$. Let $A: \cH_k \to \cH_k$ be a bounded linear operator and $A^{\star}$ be its adjoint. Then
$$ 
    \sup_{x \in \cX} \norm{A k_x}_{\cH_k}^{2} \leq \sup_{\norm{f}_{\cH_k}\leq 1} \norm{A^*f}^{2}_{L_{\infty}\paren{\cX}} \,.
$$ 
\end{lemma}

\begin{lemma} \label{lem:eig_tail_bound}
Let $\cH_k$ be some RKHS over domain $\cX \subseteq \mbr^d$ with reproducing kernel $k: \cX \times \cX \to \mbr$ and $\mu$ be any $\sigma-$finite measure on $\cX$. Let $\ell \in \mbn_{+}$ and $P : \cH_{k} \mapsto \cH_{k}$ be a projection operator with rank less than or equal to $\ell \in \mbn_{+}$.   Then
\begin{align*}
    \sum_{t > \ell} \la_t(L) \leq \int_{\cX} \|(I-P)k_x\|^2_{\cH_{k}} d\mu(x)& \leq \sup_{x \in \cX} \|(I-P)k_x\|^2_{\cH_{k}} \,,
\end{align*}
where $L: L_{2}\paren{\cX,\mu}\mapsto L_{2}\paren{\cX,\mu}$ is the kernel integral operator as defined in Equation~\eqref{eq:kern_int} and  $\la_{t}\paren{L}$ are its $t$-th eigenvalues. 
\end{lemma}

%
%

Notice that the effective dimension upper bound for the Sobolev RKHS can be also recovered from a more general result of Lemma~4 in \cite{Pagliana:20} when taking scale $\gamma = n^{\frac{1}{2s-d}}$ therein.  We provide here the proof for completeness.

\begin{proof}{\bf of Theorem~\ref{prop:eff_dim_sob}.}
Let $s > d/2$, $t\geq T_{0}$.  By Lemma~\ref{lem:proj_sob} for the orthogonal projector $P$ on the set of $t$ points $\cZ = \{x_{1},\ldots,x_{t}\} \in \cX^{t}$ such that fill distance $h_{\cZ,\cX} \lesssim t^{-\frac{1}{d}}$ for  any $\epsilon^{'} \in \mbr_{+}$ holds
$$\sup_{\|f\|_{W^{s}\paren{\cX}} \leq 1} \|f - Pf\|_{L_{\infty}\paren{\cX}} \leq C t^{-\frac{s'}{d} + \frac{1}{2}},$$
where $s'=s-\epsilon'$ and  $C$ is a constant that depends on $\cX, d, s,\epsilon$, but not on $t$.
Applying Lemma~\ref{lm:supx_Linfty} with $A = I - P$ we obtain:
$$
\sup_{x \in \cX} \|(I-P)k_x\|_{\cH_k} \leq C t^{-\frac{s'}{d} + \frac{1}{2}}.
$$
Let $ \{x_{i}\}_{i\geq 1}^{n}$ be the sequence of inputs in $\cX$. Then with the choice $\mu \bydef (1/n) \sum_{i=1}^n \delta_{x_i}$, the kernel integral operator $L$ equals $K_n/n$; combining Lemma~\ref{lem:eig_tail_bound} with the last inequality yields
\begin{equation}\label{eq:bound-eig-preliminary}
    \sum_{\ell > t} \la_t\big({K_n}/{n}\big) \leq \frac{1}{n} \sum_{i=1}^n \big\|(I-P)k_{x_i}\big\|^2_{\cH_{k}}  \leq \sup_{x \in \cX} \big\|(I-P)k_x \big\|^2_{\cH_{k}} \leq  C t^{-\frac{2s'}{d} + 1}\,.
\end{equation}
From the definition of the effective dimension (see Def.~\ref{def:eff_dim}), we have
    \begin{equation}
        \label{eq:deff_firstbound}
        d^n_{eff}(\tau) \bydef \sum_{j=1}^{n} \frac{\lambda_{j}\paren{K_{n}}}{\lambda_{j}\paren{K_n}+ \tau} \leq \sum_{j=1}^{t}\frac{\lambda_{j}\paren{K_{n}}}{\lambda_{j}\paren{K_n}+ \tau} + {\tau^{-1}}\sum_{j\geq t} {\lambda_{j}\paren{K_{n}}}\,,
    \end{equation}
where we used that  since since $K_{n}$ is positive semidefinite, $\la_j(K_n) \geq 0$ for all $j \geq 1$. Furthermore,  ${\lambda_{j}\paren{K_{n}}}/{(\lambda_{j}\paren{K_n}+\tau)} \leq 1$ for all $j \geq 1$, which implies 
\[
    \sum_{j=1}^{t}\frac{\lambda_{j}\paren{K_{n}}}{\lambda_{j}\paren{K_n}+ \tau} \leq t \,.
\]
By homogeneity of the eigenvalues we have $\la_j(K_n) = n \la_j(K_n/n)$, and therefore
\[
    {\tau^{-1}}\sum_{j\geq t} {\lambda_{j}\paren{K_{n}}} = n{\tau^{-1}}\sum_{j\geq t} {\lambda_{j}\paren{K_{n}/n}}.
\]
Combining the last two inequalities with Inequalities~\eqref{eq:bound-eig-preliminary} and~\eqref{eq:deff_firstbound}, we upper bound the effective dimension as
\[
    d^n_{eff}(\tau) \leq t + C n \tau^{-1} t^{-2s'/d + 1}.
\]
Choosing $t$ to balance the terms in the above equation, i.e. $t = n^{\frac{d}{2s'}}\tau^{-\frac{d}{2s'}}$, we get 
\[
    d_{eff}^{n}\paren{\tau} \leq C_{1} \paren{\frac{n}{\tau}}^{\frac{d}{2s'}}= C_{1} \paren{\frac{n}{\tau}}^{\frac{d}{2\paren{s-\epsilon'}}} \,.
\] 
Then assuming $\epsilon' < s/2$, and using $1/(1-x) \leq 1+2x$ for $0\leq x \leq 1/2$, we have
\[
    d_{eff}^{n}\paren{\tau} \leq  C_{1} \paren{\frac{n}{\tau}}^{\frac{d}{2s} \frac{1}{1 -\frac{\epsilon'}{s}}} \leq C_{1} \paren{\frac{n}{\tau}}^{\frac{d}{2s} (1 + \frac{2 \epsilon'}{s})} = C_{1} \paren{\frac{n}{\tau}}^{\frac{d}{2s} + \frac{d}{s^2} \epsilon'} \leq C_{1} \paren{\frac{n}{\tau}}^{\frac{d}{2s} + \frac{2 \epsilon'}{s}} \,.
\]
For any $\epsilon \in (0,1)$, the choice $\epsilon' = \epsilon s /2$ concludes the proof in the case $s \in \mbr$. 

Finally, to satisfy condition $t\geq T_{0}$ it is sufficient to have $n,\tau$ such that $\frac{n}{\tau} \geq CT_{0}^{\frac{2s^{'}}{d}}$. The latter can be alleviated by additional additive constant in the final bound. The result for $s \in \mbr_{+}$ follows. Lastly, the result implies also the particular case with $s \in \mbn$ by taking $\epsilon =0$.
\end{proof}

\section{Proof of Theorem~\ref{thm:reg_tar}}

\begin{proof} Recall that \algo, when competing against some function $f$ in an arbitrary RKHS $\cH_{k}$ with a bounded reproducing kernel, attains the general regret upper bound as given in Equation~\eqref{eq:reg_bound}. Plugging in the bound on the effective dimension of Theorem.~\ref{prop:eff_dim_sob} with $\cH_{k}=W^{s}\paren{\cX}$ into the regret upper bound~\eqref{eq:reg_bound} gives
    \begin{align}
        R_{n}\paren{\cF} \leq \tau \norm{f}_{\cH_{k}}^2 + M^2 C_1\paren[1]{ 1 +\log\paren[1]{1 + \frac{n\kappa^2}{\tau}}}\paren[1]{\widetilde{C} \paren{n \tau^{-1}}^{\frac{d}{2s} + \epsilon} +1}\,,
    \end{align}
    for any $\epsilon >0$.
    Balancing the first and second terms in order to minimize the right hand size (by choosing an appropriate value of $\tau$), i.e. by setting $\tau \bydef n^{\frac{d}{2s+d}}$ we have
    \begin{align*}
        R_{n}\paren{\cF} \leq C n^{\frac{d}{2s+d} + \epsilon}\log\paren{n},
    \end{align*}
    where a constant $C$ depends only on $d,s,R,M,\cX$ and does not depend on $n$.
\end{proof}

\section{Proof of Theorem~\ref{thm:higher_beta_p}}
\label{app:proof_higher_beta}
We start by introducing a general lemma for the regret of \algo when competing against continuous function and then proceed with the proof of the main theorem.

\begin{lemma}
\label{lm:regret-proxy-dec}
Let $f \in C(\cX)$ and $g \in \cH_k$. Assume that $\paren{x_i}_{i=1}^{n} \in \cX^{n}$ and $y_i \in [-M, M],$ for some $M > 0$. Then the regret of algorithm \eqref{eq:kern_ridge} when competing against function $f$ is bounded by
\[
    R_n(f) \leq \tau \norm{g}^{2}_{\cH_k}+ M^2\paren{1 + \log\paren{1+\frac{n\norm{k}^{2}_{\infty}}{\tau}}}d_{eff}^{n}\paren{\tau} + 2n\norm{f-g}_{L_{\infty}\paren{\cX}}\paren[1]{M+\norm{g}_{L_{\infty}\paren{\cX}}}
\]
\end{lemma}
\begin{proof}
    Let $\epsilon \in (0,1)$ and let $g \in \cH_k$ be some function which is to be chosen later. Denote by $v$ the vector $v = (f(x_1), \dots, f(x_n)) \in \mbr^n$ and $w = S_ng = (g(x_1),\dots,g(x_n)) \in \mbr^n$.
We can decompose the regret in the following way: 

\begin{align}
\label{eq:start_decom}
\begin{aligned}
    R_{n}\paren{f} & = \norm[1]{Y_n - \hat{Y}_n}^{2} - \norm{Y_{n}-v}^2\\
    & = \norm[1]{Y_n - \hat{Y}_{n}}^2 - \norm{Y_n-w}^2 - \norm{v-w}^2 + 2\inner{Y_n-w,v - w} \\
    & \leq \norm[1]{Y_n - \hat{Y}_{n}}_{2} - \norm{Y_n-w}^2 + 2\inner{Y_n-w,v-w} \\
    & \leq R_{n}\paren{g} + 2\inner{Y_n-w,v-w}. 
\end{aligned}
\end{align}
Applying the regret upper bound \eqref{eq:reg_bound}
 to the element $g$ we get: 
\begin{align*}
    R_{n}\paren{g} \leq \tau \norm{g}^{2}_{\cH_k}+ M^2\paren{1 + \log\paren{1+\frac{n\norm{k}^{2}_{\infty}}{\tau}}}d_{eff}^{n}\paren{\tau},
\end{align*}
where we recall that $d_{eff}^{n}\paren{\tau}$ is the effective dimension of the RKHS $\cH_k$ with respect to the sample $\cD \subset {\cX}^{n}$. 
For the second term on the right hand side in inequality \eqref{eq:start_decom} we have: 
\begin{align}
\label{eq:2nd_term}
    \begin{aligned}
        \inner{Y_n-w,v-w} &\leq \sum_{t=1}^{n}\abs{ \paren{y_t-g(x_t)}\paren{f\paren{x_t}-g\paren{x_t}}}  \\
        &\leq \sum_{t=1}^{n}\paren{\abs{y_t}+\abs{g\paren{x_t}}}\abs{f\paren{x_t}-g\paren{x_t}}\\
        & \leq n\norm{f-g}_{L_{\infty}\paren{\cX}}\paren[1]{M+\norm{g}_{L_{\infty}\paren{\cX}}} 
    \end{aligned}
\end{align}
Putting together the aforementioned bounds we obtain our final result.
\end{proof}

\begin{proof}{\bf of Theorem~\ref{thm:higher_beta_p}.}

Let $\sigma > 0$ be some fixed bandwidth. By Proposition~\ref{prop:sig_approx_sobolev} for any function $f \in W_{p}^{\beta}\paren{\cX} \subset W_{2}^{\beta}\paren{\cX}$, $p \geq 2$ and $\sigma >0$ there exists $f_{\sigma} \in B_{\sigma}$ such that for $0 \leq r \leq \beta$ we have: 

\begin{align}
\begin{aligned}
\label{eq:f_sig}
    \norm{f-f_{\sigma}}_{L_{2}\paren{\cX}} \leq C_1\sigma^{-\beta} \norm{f}_{W_{2}^{\beta}\paren{\cX}}, \qquad \norm{f_{\sigma}}_{W_{2}^{r}\paren{\cX}} \leq C_2 \sigma^{\paren{r-\beta}}\norm{f}_{W_{2}^{\beta}\paren{\cX}}.
\end{aligned}
\end{align}

Since $f \in W_{p}^{\beta}\paren{\cX}$ and $p \geq 2$ so the inclusion implies that we have $\norm{f}_{W_{2}^{\beta}\paren{\cX}} \leq C \norm{f}_{W_{p}^{\beta}\paren{\cX}}$ with some constant $C$. 

Let $\epsilon_1 >0$ be any positive number. Applying Sobolev embedding Theorem (see  Equation~(9) on page 60 in \cite{Edmunds:96} with $s_{1}=d/2+\epsilon_1$, $s_{2}=0$, $n=d$, $p_1 =2$, and $p_2 =\infty$), Proposition~\ref{thm:gagl_nierenberg_ineq} for a function $f-f_{\sigma} \in W_{p}^{\beta}\paren{\cX}$ and the fact that for $p\geq 2$ $W_{p}^{\beta }\paren{\cX} \subset W_{2}^{\beta p/2 }\paren{\cX}, W_{p}^{\beta}\paren{\cX} \subset W_{2}^{\beta}\paren{\cX}$ we get

\begin{align}
\begin{aligned}
\label{eq:approx_error}
    \norm{f-f_{\sigma}}_{L_{\infty}\paren{\cX}} & \leq C_{1}\norm{f-f_{\sigma}}_{W_{2}^{\frac{d}{2}+\epsilon_1}\paren{\cX}}  \hspace*{4.5cm} \leftarrow \text{by Sobolev embedding Theorem} \\
    & \leq C_2 \norm{f-f_{\sigma}}_{W_{2}^{\beta p /2}\paren{\cX}}^{\frac{d+2\epsilon_1}{\beta p}}\norm{f-f_{\sigma}}_{L_{2}\paren{\cX}}^{1 -\frac{d+2\epsilon_1}{\beta p}} \hspace*{2.1cm} \leftarrow \text{by Inequality~\eqref{eq:spec_case_p}}  \\ 
    & \leq C_{4}\norm{f-f_{\sigma}}_{W_{2}^{\beta p /2}\paren{\cX}}^{\frac{d+2\epsilon_1}{\beta p}}\paren{\sigma^{-\beta} \norm{f}_{W_{2}^{\beta}\paren{\cX}}}^{1 -\frac{d+2\epsilon_1}{\beta p}}  \hspace*{1.0cm} \leftarrow \text{ by Proposition~\ref{prop:sig_approx_sobolev}} \\
    & \leq C_{5}\norm{f-f_{\sigma}}^{\frac{d+2\epsilon_1}{\beta p}}_{W_{p}^{\beta}\paren{\cX}} \sigma^{-\beta + \frac{d}{p} +\frac{2\epsilon_1}{p}}\norm{f}_{W_{p}^{\beta}\paren{\cX}}^{1-\frac{d +2\epsilon_1}{\beta p}},\hspace*{1.8cm} \leftarrow \text{ by inclusion }
\end{aligned}
\end{align}
with a constant $C_{5}$ which does not depend on $f,f_{\sigma}$ or $\sigma$. 
Since $f_{\sigma}$ satisfies \eqref{eq:f_sig} we obtain for any $r \in \mbr_{+}$, $r \geq \beta $:
\begin{align}
\label{eq:band_norm}
    \norm{f_{\sigma}}_{W_{2}^{r}\paren{\cX}} \leq \tilde{C_1}\sigma^{\paren{r-\beta}}\norm{f}_{W_{2}^{\beta}\paren{\cX}} \leq \tilde{C_{2}}\sigma^{\paren{r-\beta}}\norm{f}_{W_{p}^{\beta}\paren{\cX}},
\end{align}
where we obtain the second inequality by inclusion of the Sobolev spaces ($W_{p}^{\beta}\paren{\cX} \subset W_{2}^{\beta}\paren{\cX}$) and the constant $\tilde{C_2}$ depends only on $\cX,d,\beta$ but not $\sigma$.
Notice  that by the triangle inequality and \eqref{eq:band_norm} with $r=\beta$ we have: 
\begin{align}
\label{eq:bound_Lp_space}
    \norm{f-f_{\sigma}}_{W_{p}^{\beta}\paren{\cX}} \leq \norm{f}_{W_{p}^{\beta}\paren{\cX}} + \norm{f_{\sigma}}_{W_{p}^{\beta}\paren{\cX} } \leq \paren{1+\tilde{C}^{\frac{1}{p}}}\norm{f}_{W_{p}^{\beta}\paren{\cX}}.
\end{align}

Thus, plugging \eqref{eq:bound_Lp_space} in the Equation~\eqref{eq:approx_error} we deduce: 
\begin{align}
\label{eq:approx_error_2}
    \norm{f - f_{\sigma}}_{L_{\infty}\paren{\cX}} \leq C_{5}\sigma^{-\beta + \frac{d}{p} +\epsilon_1} \norm{f}_{W_{p}^{\beta}\paren{\cX}}.
\end{align}

Note also that by using \eqref{eq:band_norm} with $r=s\geq \beta $ we have: 
\begin{align}
\label{eq:norm_bound}
    \norm{f_{\sigma}}_{W_{2}^{s}\paren{\cX}} \leq C_2 \sigma^{\paren{s-\beta}}\norm{f}_{W_{2}^{\beta}\paren{\cX}} \leq C_3 \sigma^{\paren{s-\beta}}\norm{f}_{W_{p}^{\beta}\paren{\cX}},
\end{align}
where the last inequality  holds since $W_{p}^{\beta}\paren{\cX} \subset W_{2}^{\beta}\paren{\cX}$. Notice that $f_{\sigma}$ as in Proposition~\eqref{prop:sig_approx_sobolev} is of limited bandwidth
and is continuous on $\cX$, therefore
$\norm{f_{\sigma}}_{L_{\infty}\paren{\cX}} = \norm{f_{\sigma}}_{C\paren{\cX}}$. Now since $f \in W_{p}^{\beta}\paren{\cX}$ and $\beta > \frac{d}{p}$, so by Sobolev Embedding Theorem $f \in C\paren{\cX}$; for the $f_{\sigma}$ chosen as in Proposition~\eqref{prop:sig_approx_sobolev} we have $$\|f_\sigma\|_{L_{\infty}\paren{\cX}} \leq \|f_\sigma\|_{W^{\beta}_p(\cX)} \leq \tilde{C}^{1/p} \|f\|_{W^{\beta}_p(\cX)},$$
where the last step is true due to \eqref{eq:band_norm}.  

Lemma~\ref{lm:regret-proxy-dec} with $g = f_\sigma \in C\paren{\cX}$, $\cH_k = W^s(\cX)$implies that for any $f \in W_{p}^{\beta}\paren{\cX}$ it holds:
\begin{align}
\label{eq:regr_bound}
\begin{split}
R_n(f) &~~\leq~~ \tau \norm{f_\sigma}^{2}_{W^s_2(\cX)}+ M^2 \log\paren{e+\frac{e n\norm{k}^{2}_{\infty}}{\tau}}d_{eff}^{n}\paren{\tau} \\
& \qquad\qquad + 2n\norm{f-f_\sigma}_{L^{\infty}\paren{\cX}}\paren[1]{M+\norm{f_\sigma}_{L_{\infty}\paren{\cX}}}.
\end{split}
\end{align} 
Denote $\epsilon = \sigma^{-1}$, $s^{'} = s - \epsilon_1$,$\beta^{'} = \beta - \epsilon_1$, By plugging \eqref{eq:approx_error_2}, \eqref{eq:norm_bound}, 
and the bound for $d^n_{eff}(\tau)$ from Theorem~\ref{prop:eff_dim_sob} in \eqref{eq:regr_bound} while noticing that $s^{'}-\beta^{'} = s-\beta$ 
we obtain for any $f$: 
\begin{align*}
    R_{n}\paren{f} & ~~\leq~~  \tilde{C}_1~\tau\epsilon^{-2\paren{s^{'}-\beta^{'}}}\norm{f}^2_{W_{p}^{\beta}\paren{\cX}} + \tilde{C}_2 M^2\paren[1]{1 + \log\paren[1]{1+\frac{n\norm{k}^2_{\infty}}{\tau}}} n^{\frac{d}{2s^{'}}}\tau^{-\frac{d}{2s^{'}}} \\
    & \qquad\qquad + \tilde{C_3} n\epsilon^{\beta^{'} - d/p}\norm{f}_{W_{p}^{\beta}\paren[1]{\cX}}\paren[1]{M + \norm{f}_{W^\beta_p\paren{\cX}}}
\end{align*}
where $\tilde{C}_1, \tilde{C}_2, \tilde{C}_3$ are constants depend on $d, \beta, s, d$, but not $n, M, \tau, \epsilon, f$.
By setting
$$
\epsilon = n^{-\frac{2s^{'}}{2s^{'}(\beta^{'} + d - d/p) - d(\beta^{'} + d/p)}}, \quad  \tau = n \epsilon^{2s^{'} - \beta^{'} - d/p} = n^{1 - \frac{2s^{'}(2s^{'} - \beta^{'} - d/p)}{2s^{'}(\beta^{'} + d - d/p) - d(\beta^{'} + d/p)} }
$$
and noticing that with such choice of $\tau,\epsilon$ for any $f \in \cF$ we have $R_{n}\paren{f} \leq C \tau \epsilon^{-2\paren[1]{ s^{'}-\beta^{'} }} =n\epsilon^{\beta^{'} - \frac{d}{p}}$
we obtain for all $f \in \cF \bydef \{ f \in W_{p}^{\beta}\paren{\cX}: \norm{f}_{W_{p}^{\beta}\paren{\cX}} \leq R \}$
$$ R_{n}\paren{\cF} = \sup_{f \in \cF} R_{n}\paren{f} \leq C n^{1 - \frac{2s^{'} (\beta^{'} - d/p)}{2 s^{'} (\beta^{'} + d - d/p) - d (\beta^{'} + d/p)}} = C n^{1 - \frac{ \beta^{'}p - d}{ \paren{\beta^{'}p +d }\paren{1-\frac{d}{2s^{'}}} + d(p-2) }},$$
where $C$ depends on $d, \beta, s, d, R, M,\cX$, but not $n$. 
Now to obtain the final claim we choose $s = \frac{d}{2} + \epsilon_1$ thus $s^{'} = \frac{d}{2}$ and we have:
$1 - \frac{ \beta^{'}p - d}{ \paren{\beta^{'}p +d }\paren{1-\frac{d}{2s^{'}}} + d(p-2) }  = 1 - \frac{\beta}{d} \frac{p - \frac{d}{\beta}}{p-2} + \frac{\epsilon_1 p}{d\paren{p-2}}$, from which the final claim follows. 
%
%
\end{proof}

\section{Proof of the lower bounds (Theorem~\ref{thm:lower_bounds_2})}
\label{app:proof_lower_bound}

To prove the lower bounds, we use the notion of the \textit{sequential} fat-shattering dimension (see Definition~12 in \cite{Rakhlin:14}). Recall (see \cite{Rakhlin:14a}) that a \textit{$\cZ$-valued tree $\bz$ of depth $n$} is a complete rooted binary tree with nodes labeled by the elements of the set $\cZ$. More rigorously, $\bz$ is a set of labeling functions $\paren{\bz_1,\ldots,\bz_n}$ such that $\bz_{t} : \{-1,1\}^{t-1} \mapsto \cZ$ for every $t\leq n$ . For any $\smash{\epsilon \in \{-1,1\}^{n}}$, we denote $\smash\{\bz_{t}\paren{\epsilon} \bydef \bz_{t}\paren{\epsilon_1,\ldots,\epsilon_{t-1}} \}$ to be the label of the node at the level $t$, which is obtained by following the path $\epsilon$. 

\begin{definition}[Fat-shattering dimension, see Definition~7 in \cite{Rakhlin:14a} ]
	\label{def:fat_dim}
	Let $\gamma>0$. An $\cX$-valued tree $\bx$ of depth $d$ is said to be $\gamma$-shattered  by $\cF=\{f: \cX \mapsto \mbr\}$ if there exists an $\mbr-$valued tree $\bs$ of depth $d$ such that 
	\[
	\forall \epsilon \in \{-1,1\}^{d}, \quad \exists f^{\epsilon} \in \cF, \quad \text{s.t.} \quad \epsilon_t \paren{ f^{\epsilon}\paren{\bx_{t}\paren{\epsilon}} -\bs_{t}\paren{\epsilon} } \geq \frac{\gamma}{2} \,,
	\]
	for all $t \in \{1,\dots,d\}$. The tree $\bf{s}$ is called a \textit{witness}. The largest $d$ such that there exists a $\gamma$-shattered tree $\bx$ is called the ({sequential}) fat-shattering dimension of $\cF$ and is denoted by $\fatt_{\gamma}\paren{\cF}$. 
\end{definition}

If the last inequality becomes equality, we say that the tree $\bx$ is \textit{exactly} shattered by the elements of $\cF$ or (alternatively) that class $\cF$ exactly shatters the tree $\bx$. We recall also the notion of sequential covering numbers and the sequential entropy of class $\cF$. 
\begin{definition}
	A set $V$ of $\mbr-$valued trees of depth $n$ forms a $\gamma-$ cover (with respect to the $\ell_q$ norm, $1\leq q < \infty$) of a function class $\cF \subset \mbr^{\cX}$ on a given $\cX-$valued tree $\bx$ of depth $d$ if
	\begin{align*}
	\forall f \in \cF, \forall \epsilon \in \{\pm 1\}^{d}, \exists \bv \in V, \text{ s.t. } \paren{\frac{1}{n}\sum_{t=1}^{d} \abs{f\paren{\bx_t \paren{\epsilon}} - \bv_{t} \paren{\epsilon}}^{q} }^{1/q} \leq \gamma.
	\end{align*}
	In the case $q=\infty$, we have that $\abs{f\paren{\bx_{t} \paren{\epsilon}} - \bv_{t} \paren{\epsilon}} \leq \gamma$ for all $t \in \{1,\ldots,d\}$. The size of the smallest $\gamma$-cover of a tree $\bx$ is denoted by $\cN_{q}\paren{\gamma, \cF,\bx}$; and $\cN_{q}\paren{\gamma,\cF,d} = \sup_{\bx}\cN\paren{\gamma, \cF,\bx}$ where the last supremum is taken over all trees of depth $d$. Finally, the sequential entropy of class $\cF$ is $\sup_{\bx}\log \cN_{q}\paren{\gamma, \cF,\bx}$. 
\end{definition}

To derive the main results of Theorem~\ref{thm:lower_bounds_2}, we use the following consequences of Lemmata 14,15 in Section~5, \cite{Rakhlin:14}.

\begin{lemma}[Variant of \bf{Lemma~14} in \cite{Rakhlin:14}]
	\label{lem:RS_14}
	Let $n \in \mbn_{*}$, $\cY = [-M,M]$ and 
	$\cF \subseteq \big\{ f : \cX \to [-M/4,M/4]\big\}$ for some $M>0$. 
	If 
	$\gamma >0$ such that $n \leq \fatt_{\gamma}\paren{\cF}$ then
	\begin{align*}
	\tilde{R}_{n}\paren{\cF} \geq \frac{M}{4} n \gamma.
	\end{align*}
\end{lemma}

\begin{proof}
	Since,
	$\gamma>0$ such that $n \leq \fatt_{\gamma}\paren{\cF}$, by definition of the fat-shattering dimension there exists an $\cX-$valued tree $\bx$ of depth $n$ (and a witness of shattering $\boldsymbol{\mu}$), which is shattered by the elements of $\cF$. Further proof follows the same lines as in the original argument of Lemma~14 of \cite{Rakhlin:14} with the tree $\bx$, witness of shattering $\boldsymbol{\mu}$, $\beta := \gamma$ and functions (as well as witness of shattering bounded in $[-\frac{M}{4},\frac{M}{4}]$ instead of $[-1,1]$) therein.
\end{proof}

\begin{lemma}[Variant of \textbf{Lemma~15 in \cite{Rakhlin:14}}]
	\label{lem:RS_15} 
	Let $n \in \mbn_*$, $\gamma>0$, and $\cF'$ be a class of functions from $\cX$ to $[-M/4,M/4]$
	which exactly $\gamma$-shatters some tree $\bx$ of depth $\fatt_{\gamma}\paren[0]{\cF'} < n$. Then the minimax regret with respect to $\cF'$ is lower-bounded as
	\begin{align}
	{\tilde{R}_{n}\paren[1]{\cF'}} \geq \frac{M}{4}C\paren[2]{2\sqrt{2}\gamma \sqrt{n {\fatt_{\gamma}\paren{\cF'}}} - n\gamma^2}.
	\end{align}
	
\end{lemma}
\begin{proof}
	The lemma is proved in the same way as Lemma~15 in \cite{Rakhlin:14}, by noting that since $\cF'$ exactly shatters $\bx$, we can consider $\cF=\cF'$ in the original proof. The argument follows then the same lines by noticing that the target functional class is a subset of $\{ f: \cX \mapsto [-\frac{M}{4},\frac{M}{4}]\}$ (instead of $\{ f: \cX \mapsto [-1,1]\} $ as in the original argument). 
\end{proof}

To prove the lower bounds, we provide a tight control of $\fatt_{\gamma}\paren{\cF}$ (in terms of the scale $\gamma$, while constants may depend on the range $\cY$,$d$,$\beta$) 
for $\cF$ being the bounded ball in Sobolev space $W_{p}^{\beta}\paren{\cX}$.

We recall the notion of sequential Rademacher complexity (see \cite{Rakhlin:14}): 
\begin{align*}
\cR_{n}\paren{\cF} = \sup_{\bx} \ee{\epsilon}{n^{-1}\sup_{f \in \cF} \sum_{t=1}^{n}\epsilon_{t}f\paren{\bx_{t}\paren{\epsilon}}},
\end{align*}
where $\ee{\epsilon}{\cdot}$ denotes the expectation under the product measure $\mbp =( \frac{1}{2}\delta_{-1} + \frac{1}{2}\delta_{1})^{\otimes n}$, the supremum is over all $\cX-$ valued trees of depth $n$. 
Firstly we provide an auxiliary Lemma which provides an upper bound of the fat-shattering dimension of the Sobolev ball $B_{W_{p}^{\beta}\paren{\cX}}(0,1)$.

\begin{lemma}
	\label{lem:rad_up_bound}
	Let $n \in \mbn, n \geq 1$, $M>0$ an let $\cF \bydef B_{W_{p}^{\beta}\paren{\cX}}\paren{0,M/4}$, $p \geq 2$.
	For the fat-shattering dimension $\fatt_{\gamma}\paren{\cF}$ on the scale $\gamma >0$ when $\beta \neq \frac{d}{2}$ it holds
	\begin{align*}
	\fatt_{\gamma}\paren{\cF} 
	\leq \max\{ \tilde{C}_1\gamma^{- \paren{\frac{d}{\beta} \lor 2 }},1\},
	\end{align*}
	where ${C_1}$ is some constant which depends on $\beta,d,M$
	but not on $\gamma$. In the case $\frac{\beta}{d} = 1/2$ we have
	\begin{align*}
	\fatt_{\gamma}\paren{\cF} 
	\leq \max\{ \tilde{C}_2\paren{\frac{\gamma}{\log(\gamma)}}^{-2},1\},
	\end{align*}
	where ${C}_{2}$ is some constant which depends on $\beta,d,M$.
	but not on $\gamma$. 
\end{lemma}

\begin{proof}{\textbf{of Lemma~\ref{lem:rad_up_bound}} }
	Following from Definition~\eqref{def:fat_dim}, if $\bx$ of depth $n$ is $\gamma-$shattered by the elements of $\cF$, then $n\leq \fatt_{\gamma}\paren{\cF}$. For an arbitrary functional class $\cF$  from the definition of the fat-shattering dimension for any $\gamma >0$ such that $\fatt_{\gamma}\paren{\cF} > n$  we have that $\cR_{n}\paren{\cF} \geq \frac{\gamma}{2}$ (one readily checks this by considering Rademacher complexity over the set of $n$ shattered points). Therefore, $\cR_{n}\paren{\cF} \geq \sup\{ \frac{\gamma}{2}: \fatt_{\gamma}\paren{\cF} > n\}$, which is equivalent to $\fatt_{\gamma}\paren{\cF} \leq \min\{n : \cR_{n}\paren{\cF} \leq \frac{\gamma}{2}\}$.
	By Proposition~1 and Definition~3 in \cite{Rakhlin:14a} for all $c \in \mbr$, we have $\cR_{n}\paren{c\cF} = \abs{c} \cR_{n}\paren{\cF}$, where $c\cF = \{cf: f \in \cF\}$. Taking $c=\frac{4}{M}$, we have for $\cF^{'} = B_{W_{\infty}^{\beta}\paren[1]{\cX}}(0,\frac{M}{4})$ that $\cR_{n}\paren[1]{\cF^{'}} = \frac{M}{4}\cR_{n}\paren[1]{\cF}$, where $\cF=B_{W_{\infty}^{\beta}\paren[1]{\cX}}(0,1)$. 
	From the definition of $\norm{\cdot}_{W_{\infty}^{\beta}\paren{\cX}}$, it follows that if $f \in B_{W_{\infty}^{\beta}\paren{\cX}}\paren{0,1}$, then $\max_{x \in \cX}{\abs{f(x)}} \leq 1$.
	By Theorem~3 in \cite{Rakhlin:14a} we have for any functional class
	$\cF \subset [-1,1]^{\cX}$
	\begin{align}
	\label{eq:rad_comp1}
	\cR_{n}\paren[1]{\cF} \leq \sup_{\bx} \inf_{\rho \in (0,1]} \paren{4 \rho + \frac{12}{\sqrt{n}}\int_{\rho}^{1}\sqrt{\log_2 \cN_{2}\paren{\delta,\cF,\bx}}d\delta }.
	\end{align}
	It is straightforward to check that for any tree $\bz$ it holds that 
	\begin{align}
    	\label{eq:entr_ineq}
	    \cN_{2}\paren{\gamma,\cF,\bz} \leq \cN_{\infty}\paren{\gamma,\cF,\bz}.
	\end{align}
	
	Furthermore, if $\cN_{\infty}\paren{\cF,\gamma}$ is a metric entropy of class $\cF$ on scale $\gamma>0$, then it is easy to check that for any tree $\bz$ of depth $d\geq 1$ and any scale $\gamma>0$, $\cN_{\infty}\paren{\gamma,\cF,\bz} \leq \cN_{\infty}\paren{\gamma,\cF}$. Indeed, this follows trivially by taking for any tree $\bz$ witness $v(\cdot)= g\paren{\bz\paren{\cdot}}$, where $g\paren{\cdot}$ is the element of $\gamma-$net such that $\norm{f-g}_{\infty} \leq \gamma$.
	Furthermore, for $\cF = B_{W_{p}^{\beta}\paren{\cX}}(0,1), \beta > d/p$ the metric entropy of $\cF$ on the scale $\delta$ is (up to some constant $C$ which does not depend on $\delta$) upper bounded by $\delta^{-\frac{d}{\beta}}$. The latter bound is a well-known result and it can be deduced from the general result for Besov spaces stated in Theorem 3.5 in \cite{Edmunds:96} (see also Equation~(38) on page 19 in \cite{Vovk:06}). 
	Thus, using Equations~\eqref{eq:rad_comp1} and \eqref{eq:entr_ineq}, the fact that metric entropy uniformly bounds sequential entropy, properties of Rademacher complexity (see Lemma 3 in \cite{Rakhlin:15}) and the upper bound on the metric entropy of the Sobolev ball 
	$\cF = B_{W_{\infty}^{\beta}\paren{\cX}}(0,M/4)$,
	we get
	\begin{align}
	    \begin{aligned}
	    \label{eq:rad_int_bound}
	    \cR_{n}\paren{\cF} &= \frac{M}{4}\cR_{n}\paren[1]{ \frac{4}{M}B_{W_{\infty}^{\beta}}\paren{0,\frac{M}{4}}} \leq \frac{M}{4} \cR_{n}\paren[1]{B_{W_{\infty}^{\beta}\paren{\cX}}(0,1) } \\ 
	    &\leq \frac{M}{4}\inf_{\rho \in (0,1]}\paren[1]{4\rho + \frac{12 C_{}}{\sqrt{n}}\int_{\rho}^{1}\delta^{-\frac{d}{2\beta}}d\delta  }
	    \leq C_{1}\inf_{\rho \in (0,1]}\paren{4\rho + \frac{12}{\sqrt{n}}\int_{\rho}^{1}\delta^{-\frac{d}{2\beta}}d\delta  },
	    \end{aligned}
	\end{align}
	where we use $C_{1}=\frac{M}{4}\max\{1,C\}$ for completeness. 
	Notice that if $\beta > \frac{d}{2}$, then integral $\int_{0}^{1}t^{-\frac{d}{2\beta}}dt$ is finite, thus in this case in \eqref{eq:rad_int_bound} we can take $\rho = 0$, which implies $\cR_{n}\paren{\cF} \leq \frac{12C_1}{\sqrt{n}} \frac{1}{1-\frac{d}{2\beta}}$. When $\beta < \frac{d}{2}$, then the choice $\rho = \rho_{min} = (9n^{-1})^{\frac{\beta}{d}}$ leads to the bound $\cR_{n}\paren{\cF} \leq 12 C_1 n^{-\frac{\beta}{d}}\frac{1}{1-\frac{2\beta}{d}}$. Finally, in the case when $\beta=\frac{d}{2}$ with the choice $\rho =\frac{3}{\sqrt{n}}$, one gets $\cR_{n}\paren{\cF} \leq 6 \frac{C_{1}\ln (n)}{\sqrt{n}}$. 
	
	Thus we obtain 
	\begin{align}
	\label{eq:rad_bound_ex}
	\cR_{n}\paren[1]{\cF} \leq 12 C_{1}K_{}n^{-\paren{\frac{\beta}{d} \wedge \frac{1}{2} }},
	\end{align}
	where in Equation~\eqref{eq:rad_bound_ex} $K = \frac{1}{1-(\frac{2\beta}{d} \wedge \frac{d}{2\beta})}$ if $\beta \neq \frac{d}{2}$ otherwise $K = \frac{\ln(n)}{2}$.
	If $\frac{\beta}{d} \neq \frac{1}{2}$ then we have
	\begin{align*}
	\fatt_{\gamma}\paren{\cG} \leq \fatt_{\gamma}\paren{\cF}
	& \leq \min\{n : \cR_{n}\paren{\cF} \leq \frac{\gamma}{2}\} \\
	& \leq \min \bigg \{ n : 12C_{1}Kn^{-\paren{\frac{\beta}{d}\land \frac{1}{2}}} \leq \frac{\gamma}{2} \bigg \} \\
	& \leq \lceil  \paren{\frac{\gamma }{ 24 C_{1} K }}^{-\paren{\frac{d}{\beta}\lor 2}} \rceil \\
	& \leq \max\{ C_{2} \gamma^{- \paren{\frac{d}{\beta}\lor 2}} ,1\}.
	\end{align*}
	with $C_{2} = 2\cdot(24C_{1}K)^{\frac{d}{\beta} \lor 2}$. In the case, when $\frac{\beta}{d} = \frac{1}{2}$ we have that by any $n \geq  \lceil \paren{\frac{\gamma/24C_1}{\log(\gamma/24C_1)}}^{-2}\rceil $ ensures that $\frac{6C_{1}\ln(n)}{\sqrt{n}} \leq \frac{\gamma}{2}$, from which we deduce $\fatt_{\gamma}\paren{\cG} \leq \max\{ C_{2} \paren{\frac{\gamma}{\log(\gamma)}}^{-2},1\}$.
\end{proof}

To derive the first statement of Theorem~\ref{thm:lower_bounds_2}  we construct a class $\cG \subset B_{W_{p}^{\beta}\paren{\cX}}\paren{0,M}$ which satisfies Lemmata~\ref{lem:RS_14} and \ref{lem:RS_15} and deduce the final bound for the minimax regret $\tilde{R}_{n}\paren{B_{W_{p}^{\beta}\paren{\cX}}\paren{0,M}}$ by inclusion argument.

\paragraph{Class construction.} We provide a class construction, taking inspiration from the nonparametric regression in the statistical learning scenario (see, for example, Theorem 3.2 in \cite{Gyorfi:02}). 
Recall that $\cX = [-1,1]^d$; for a given $n \in \mbn$ denote $b \bydef  n^{-\frac{1}{d}}$. 
Consider the following set of half--open intervals
$$A = \{ A_{\ell}=[-1+\ell b,-1+(\ell+1)b), 0 \leq \ell \leq \lfloor 2n^{1/d} \rfloor -1\},$$ and let $\cP = A^{d}$ be its $d-$th power. Let $I \bydef \{0,\ldots,\lfloor 2n^{\frac{1}{d}}\rfloor -1\}^{d}$, $N = \abs{I} = \lfloor 2n^{\frac{1}{d}}\rfloor^{d}$ and $\pi : I \mapsto \{1,\ldots,N\}$ be a function which maps an element $k \in I$ to its index in the lexicographic order among the elements in $I$. Because lexicographic order is a total order, we have that $\pi(\cdot)$ is a bijection.  
For each $k \in I \bydef \{0,\ldots,\lfloor 2n^{\frac{1}{d}} \rfloor -1\}^{d}$ such that $\pi(k) = j$, we denote $B_{j} = \prod_{i=1}^{d}[-1+k_{i}b,-1+(k_{i}+1)b)$. Notice that $\cup_{j=1}^{N}B_{j} \subset \cX$ and for $i\neq j$ obviously $B_{i} \cap B_{j} = \emptyset$.
For a cube $B_{t}$, $t \in \{1,\ldots,N\}$ we denote  $a_t \in \mbr^d$ to be its center. One can show explicitly that $a_{t}=\paren{b\paren{\frac{1}{2}+\paren{\pi^{-1}\paren{t}}_{1}}-1,\ldots, b\paren{\frac{1}{2}+\paren{\pi^{-1}\paren{t}}_{d}}-1}$. Consider the following set of functions:
\begin{align}
\label{eq:func_set}
\cF_{\beta,d,n} = \bigg \{f: f\paren{x} = \frac{M n^{-\frac{\beta}{d}}}{4 {\norm{g}_{W_{\infty}^{\beta}\paren{\cX}}}}\sum_{t=1}^{N} c_{t} {g_{n,t}\paren{x}}, c_{j} \in \{-1,1\} \bigg \},
\end{align}
where $g_{n,t}\paren{x} = g\paren[1]{n^{\frac{1}{d}}\paren{x-a_{t}} }$, and $g$ such that
$g\paren{x} = \frac{1}{2}\paren{1 - \sigma(\frac{\norm{x}_{2}^2 - a^2}{c^2-a^2} )}$, $c =\frac{1}{2}$, $a = \frac{1}{4}$ and  $\sigma\paren{t}= \frac{h\paren{t}}{h\paren{t} + h\paren{1-t}}$, $h\paren{t} = e^{-1/t^{2} }\mbi_{t > 0}$ for $t \in \mbr$, $x \in \mbr^d$. 
We need the following Lemma, which shows that the functional class $\cF_{\beta,d,n}$ defined by Equation~\eqref{eq:func_set} is included in the ball of the space 
$W_{\infty}^{\beta}\paren{\cX}$.
\begin{lemma}
	\label{lem:f_incl}
	Let $\beta>0$, $d\geq 1$, $n \in \mbn$; consider $N \bydef \lfloor 2n^{\frac{1}{d}}\rfloor^{d}$ and the class $\cF_{\beta,d,n}$, as defined in \eqref{eq:func_set}. It holds that 
\[
	\cF_{\beta,d,n} \subset B_{L_{\infty}\paren{\cX}}\paren[2]{0,\frac{M}{4}}.
\]
	Moreover, a stronger inclusion holds, namely, that
\[
    \cF_{\beta,d,n} \subset B_{W_{\infty}^{\beta}\paren{\cX}}\paren[2]{0,\frac{M}{4}}.
\]
\end{lemma}
\begin{proof}
	First, notice that $g\paren{\textbf{0}} = \frac{1}{2}\paren[1]{1-\sigma(-\frac{a^2}{c^2-a^2})}$.  Because $t_{0}\bydef -\frac{a^2}{c^2-a^2} < 0$,  $h\paren{t_0} = 0$ and consequently $\sigma\paren{t_0} =0$ from which we have $g\paren{0}= \frac{1}{2}(1-\sigma(t_{0})) = \frac{1}{2}$. 
	For a cube $B_{j}$, if $x \notin B_{j}$, then we have $g_{n,j}\paren{x} =0$. Indeed, as $x \notin B_{j}$, for $a_{j}$ center of $B_{j}$ holds $\max_{i\leq d}\abs{x^{(i)}-a^{\paren{i}}_{j}} \geq \frac{n^{-\frac{1}{d}}}{2}$. Therefore, because $\norm[1]{n^{\frac{1}{d}}\paren{x^{(i)}-a^{(i)}_{j}}}_{2} \geq n^{\frac{1}{d}}\max_{i\leq d}\abs{x-a_{j}} \geq \frac{1}{2}$ and because $g\paren{\cdot}$, as constructed above is a mollifier from $\mbr^d$ to $\mbr$ with non-zero support on $B_{\mbr^d}\paren{0,1/2}$ (see paragraph 13 in \cite{WTu:11}), we have $g_{n,j}\paren{x} = g\paren{n^{\frac{1}{d}}\paren{x-a_{j}}}=0$.
	From the definition of the norm in the functional class $W_{\infty}^{\beta}\paren{\cX}$, it follows that
	for any $x \in \cX$ we have
	$ \abs{g\paren{x}} \leq \norm{g}_{L_{\infty}\paren{\cX}} \leq  \norm{g}_{W_{\infty}^{\beta}\paren{\cX}}$.
	Furthermore, for any element $f \in \cF_{\beta,d,n}$, for any $x \in \cX \setminus \cup_{k =1 }^{N} B_{k}$ we have $f\paren{x}=0$. If $x \in \cup_{k=1}^{N}B_{k}$,
	then there exists some cube $B_{j}$ with $x \in B_j$. Thus we get 
	\begin{align*}
	    \abs{f\paren{x}} &= \left| \frac{M n^{-\frac{\beta}{d}}}{4 {\norm{g}_{W_{\infty}^{\beta}\paren{\cX}}}}\sum_{t=1}^{N} c_{j} {g_{n,t}\paren{x}}\right|  \leq \frac{M}{4}{n^{-\frac{\beta}{d}} \frac{ \abs{g_{n,j}\paren{x}}}{\norm{g}_{W_{\infty}^{\beta}\paren{\cX}}}} \\
	    & \leq \frac{M}{4}n^{-\frac{\beta}{d}} \frac{\norm{g}_{L_{\infty}\paren{\cX} }}{\norm{g}_{W_{\infty}^{\beta}\paren{\cX}}} \leq \frac{M}{4},
	\end{align*}
	so that $\cF_{\beta,d,n} \subset B_{L_{\infty}\paren{\cX}}\paren{0,\frac{M}{4}}$ and the first part of the claim is proved.
	Let $\beta = m +\sigma$. For every $r \leq m$, $r \in \mbn$ and $x \in \cX$, we notice that if $x\ \in \cX \setminus \cup_{k=1}^{N} B_{k}$, then because it is a finite linear combination of mollifiers we have $D^{r}f\paren{x} = 0$. By a chain rule for every $f \in \cF_{\beta,d,n}$, $x \in \cX$, $k \leq N$ such that $x \in B_{j}$:
	\begin{align*}
	\sup_{x \in \cX} \abs{D^{r}f\paren{x}} & = \sup_{B_{j} \in \cP} \sup_{x \in B_{j}} \abs{D^{r}f\paren{x}} \\
	&=  \sup_{B_{j} \in \cP} \sup_{x \in B_{j}} \frac{M}{4\norm{g}_{W_{\infty}^{\beta}\paren{\cX}}} \abs{D^{r}n^{-\frac{\beta}{d}}g_{n,j}\paren{x}}\\
	& = \sup_{B_{j} \in \cP} \sup_{x \in B_{j}} \frac{M n^{-\frac{\beta}{d}}}{4 \norm{g}_{W_{\infty}^{\beta}\paren{\cX}} }\abs{D^{r}g\paren{n^{\frac{1}{d}}\paren{x-a_j} }}\\
	&=  \frac{M}{4\norm{g}_{W_{\infty}^{\beta}\paren{\cX}}}n^{\frac{r-\beta}{d}} \sup_{B_{j} \in \cP} \sup_{x \in B_{j}} \abs{D^{r}g\paren{x}}\\
	& \leq \frac{M}{4} \frac{\sup_{x \in \cX } \abs{D^{r}g\paren{x}}}{\norm{g}_{W_{\infty}^{\beta}\paren{\cX}}} = \frac{M}{4}\frac{\norm{D^{r}g}_{L_{\infty}\paren{\cX}}}{\norm{g}_{W_{\infty}^{\beta}\paren{\cX}}} \leq \frac{M}{4}
	\end{align*}
	Consider $D^{\gamma}f$ of a function $f \in \cF_{\beta,d,n}$. For some $1 \leq j \leq N$ we have for any
	$x,z, \in \overline{B_{j}}$ ( here  $\overline{B_{j}} = B_{j} \cup \partial B_{j}$)  it holds
	\begin{align*}
	\frac{\abs{D^{\gamma}f\paren{x} - D^{\gamma}f\paren{z}}}{\norm{x-z}^{\sigma}} &=\frac{Mn^{-\frac{\beta}{d}} }{4 \norm{g}_{W_{\infty}^{\beta}\paren{\cX}}}\frac{\abs{D^{\gamma}g_{n,j}\paren{x} - D^{\gamma}g_{n,j}\paren{z}}}{\norm{x-z}^{\sigma}} \\
	& = \frac{M n^{-\frac{\beta}{d}}}{4 \norm{g}_{W_{\infty}^{\beta}\paren{\cX}} } \frac{\abs{D^{\gamma}g\paren[1]{n^{\frac{1}{d}}(x-a_{j})} - D^{\gamma}g\paren[1]{n^{\frac{1}{d}}(z-a_{j})} } }{ \norm{x-z}^{\sigma} } \\
	& = \frac{M n^{-\frac{\beta}{d} }}{4 \norm{g}_{W_{\infty}^{\beta}\paren{\cX}} } \frac{\abs{D^{\gamma}g\paren{\overline{x}}\frac{\partial^{\gamma}}{\partial x_{1}\ldots \partial x_{d}}n^{\frac{1}{d}}\paren{x-a_{j}} - D^{\gamma}g\paren{\overline{z}} \frac{\partial^{\gamma}}{\partial z_{1}\ldots \partial z_{d}}n^{\frac{1}{d}}\paren{z-a_{j}} } }{ n^{-\frac{\sigma}{d}}\norm{\overline{x}-\overline{z}}^{\sigma} } \\
	& \leq \frac{M n^{-\frac{\beta}{d} + \frac{m}{d}+\frac{\sigma}{d} }}{4 \norm{g}_{W_{\infty}^{\beta}\paren{\cX}} } \sup_{x,z \in \cX, x\neq z } \frac{\abs{D^{\gamma}g\paren{x} - D^{\gamma}g\paren{z} } } {\norm{x-z}^{\sigma} } \\
	& = \frac{M}{4}\frac{1}{\norm{g}_{W_{\infty}^{\beta}\paren{\cX} }}\sup_{x,z \in \cX,x\neq z} \frac{ \abs{ D^{\gamma} g\paren{x} -D^{\gamma} g\paren{z}}}{\norm{x-z}^{\sigma}} \leq \frac{M}{4}
	\end{align*}
	Furthermore, if $B_{j},B_{k} \in \cP$ are two different cubes then for $x \in \overline{B_{j}}$ and $z \in \overline{B_{k}}$ consider elements $\overline{x} \in \partial B_{j}$ and $\overline{z} \in \partial B_{k}$, which lie on the line between $x$ and $z$. Notice that if $\overline{B_{j}}$ and $\overline{B_{k}}$ have common $d-1$ hyperplane (i.e., they are the neighbour cells) then $\overline{x} = \overline{z}$. In all cases, it follows from the construction of $f \in \cF_{\beta,d,n}$ that $D^{\gamma}f\paren{\overline{x}}=D^{\gamma}f\paren{\overline{z}}=0$. Therefore, we have
	\begin{align*}
	\frac{\abs{D^{\gamma}f\paren{x} - D^{\gamma}f\paren{z}}}{\norm{x-z}^{\sigma}} &=      \frac{Mn^{-\frac{\beta}{d}} }{4 \norm{g}_{W_{\infty}^{\beta}\paren{\cX}}}\frac{ \abs{D^{\gamma}g_{n,j}\paren{x} - D^{\gamma}g_{n,j}\paren{\overline{x}}-D^{\gamma}g_{n,k}\paren{\overline{z}}+ D^{\gamma}g_{n,k}\paren{z}} }{\norm{x-z}^{\sigma}} \\
	& \leq \frac{Mn^{-\frac{\beta}{d}} }{4 \norm{g}_{W_{\infty}^{\beta}\paren{\cX}}}\frac{ \abs{D^{\gamma}g_{n,j}\paren{x} - D^{\gamma}g_{n,j}\paren{\overline{x}}}+\abs{D^{\gamma}g_{n,k}\paren{\overline{z}}- D^{\gamma}g_{n,k}\paren{z}} }{\norm{x-z}^{\sigma}}\\
	& \leq \frac{Mn^{-\frac{\beta}{d}}}{4\norm{g}_{W_{\infty}^{\beta}\paren{\cX} } } \frac{ \norm{g}_{W_{\infty}^{\beta}\paren{\cX} } n^{\frac{\beta}{d}} \paren{\norm{x-\overline{x}}^{\sigma} + \norm{z-\overline{z}}^{\sigma} } }{\norm{x-z}^{\sigma} } \\
	& \leq \frac{M}{4} 2^{\sigma}\frac{\norm{\overline{x}-x}^{\sigma}+\norm{\overline{z}-z}^{\sigma}}{ 2\norm{x-z}^{\sigma} }\\
	& \leq \frac{M}{4} 2^{\sigma} \paren{ \frac{\norm{x - \overline{x}} + \norm{z - \overline{z}}}{2}}^{\sigma} \frac{1}{\norm{z-x}^{\sigma}} \\
	& \leq \frac{M}{4} \frac{\norm{x-z}^{\sigma}}{\norm{x-z}^{\sigma}} = \frac{M}{4} 
	\end{align*}

	If for any pair $(x,z) \in \cX^2$, $x \neq z$ one element (without losing of generality let it be $z$) does not belong to the union of the cubes $\cup_{B \in \cP} B$, then we can substitute this point by the  point $\overline{z}$, which is the intersection of the segment $[x,z]$ and the boundary of the closest cube to the point $z$. Notice that in this case $D^{\gamma}f\paren{z} = D^{\gamma}f\paren{\overline{z}} =0$ by construction of $f$ and $\norm{x-z}_{2}^{\sigma} \geq \norm{x-\overline{z}}_{2}^{\sigma}$. Applying aforementioned analysis to a pair $(x,\overline{z})$ which lies in some (different) cubes $B_j, B_k$, we get 
	$$\frac{\abs{D^{\gamma}f(x) - D^{\gamma}f(z)}}{\norm{x-z}^{\sigma}} \leq \frac{\abs{D^{\gamma}f(x) - D^{\gamma}f(\overline{z})}}{\norm{x-\overline{z}}^{\sigma}} \leq \frac{M}{4}.$$
	Finally, case $(x,z) \in \cX^2$, where none of the points belong to the union of the cubes, is trivial. 
	
	Considering these cases together we have $\sup_{x,y \in \cX, x\neq y}\frac{\abs{D^{\gamma}f(x)-D^{\gamma}f(y)}}{\norm{x-y}^{\sigma}} \leq \frac{M}{4}$ for any $f \in \cF_{\beta,d,n}$. Therefore, $\cF_{\beta,d,n} \subset  B_{W_{\infty}^{\beta}\paren{\cX}}\paren{0,\frac{M}{4}}$.   
	
\end{proof}

\begin{proof}{\textbf{{ of Theorem~\ref{thm:lower_bounds_2}}}.}
	For $n \geq 1$ consider the functional class $\cF_{\beta,d,n}$ as given by Equation~\ref{eq:func_set}. Consider a $\cX-$valued tree $\bx$ of depth $N := \lfloor 2n^{\frac{1}{d}} \rfloor ^{d}$ constructed as follows: for any $\epsilon \in \{-1,1\}^{N}$, any $t \leq N$ we set $x_t \paren{\epsilon} = a_{t}$, where $a_{t}$ is the center of the correspondent cube. Now, for any $\epsilon \in \{-1,1\}^{n}$ consider $f^{\epsilon}\paren{\cdot} \in \cF_{\beta,d,n}$
	where $\cF_{\beta,d,n}$ as in \eqref{eq:func_set}  and
	$f^{\epsilon}\paren{x} = \frac{M}{4 \norm{g}_{W_{\infty}^{\beta}\paren{\cX}}}\sum_{j=1}^{N}\epsilon_{j}n^{-\frac{\beta}{d}}g_{n,j}\paren{x}$. Then for the tree $\bx$, for every $\epsilon \in \{-1,1\}^{N}$, $1 \leq t \leq N$ and a real-valued (witness of shattering) $s_{t}\paren{\cdot}\bydef 0$, we have
	
	\begin{align}
	\begin{aligned}
	\label{eq:expres_shat}
	\epsilon_{t} \paren{ f^{\epsilon}\paren{x_{t}\paren{\epsilon}} -s_{t}\paren{\epsilon} } = \epsilon_t f^{\epsilon}\paren{x_t\paren{\epsilon}} = \epsilon_t f^{\epsilon}\paren{a_t}&= \frac{M}{4 \norm{g}_{W_{\infty}^{\beta}\paren{\cX}}}\epsilon_t \sum_{j=1}^{N}\epsilon_j n^{-\frac{\beta}{d}}g_{n,j}\paren{a_t} \\
	&=  \frac{M}{4 \norm{g}_{W_{\infty}^{\beta}\paren{\cX}}}n^{-\frac{\beta}{d}}g_{n,t}\paren{a_{t}} \\
	& = \frac{M}{4 \norm{g}_{W_{\infty}^{\beta}\paren{\cX}}}n^{-\frac{\beta}{d}}g_{}\paren{0}
	= C_{M,g}\frac{n^{-\frac{\beta}{d}}}{2},
	\end{aligned}
	\end{align}
	where $C_{M,g}:=\frac{M}{4 \norm{g}_{W_{\infty}^{\beta}\paren{\cX}} }$.
	{Thus, class $\cF_{\beta,d,n}$ with $\tilde{\gamma} = \tilde{\gamma}\paren{n} \bydef C_{M,g}{n^{-\frac{\beta}{d}}}$ (exactly) shatters the tree $\bx$}. Notice that $N = \lfloor 2n^{\frac{1}{d}}\rfloor^{d} \leq 2^{d}n$; from the other side we have $N \geq \paren[1]{ n^{\frac{1}{d}}}^{d} \geq n $. Thus, from the definition of fat-shattering dimension, it follows,
	\begin{align}
	\label{eq:fatt_lwr_bound}
	\fatt_{\tilde{\gamma}}\paren[2]{\cF_{\beta,d,n}}\geq  N \geq n.
	\end{align}
	
	All conditions of Lemma~\ref{lem:RS_14} are fulfilled for the class $\cF_{\beta,d,n}$; by Lemma~\ref{lem:f_incl} $\cF_{\beta,d,n} \subset B_{W_{\infty}^{\beta}\paren{\cX}}\paren{0,\frac{M}{4}}$. Applying Lemma~\ref{lem:RS_14} to the class $\cF_{\beta,d,n}$, using Lemma~\ref{lem:f_incl} and simple inclusion $B_{W_{\infty}^{\beta}\paren{\cX}}\paren{0,\frac{M}{4}} \subset B_{W_{p}^{\beta}\paren{\cX}}\paren{0,\frac{M}{4}} \subset  B_{W_{p}^{\beta}\paren{\cX}}\paren{0,M}$, we obtain for the Sobolev ball $\cF \bydef B_{W_{p}^{\beta}\paren{\cX}}\paren{0,M}$ 
	\begin{align*}
	\tilde{R}_{n}\paren{ \cF } \geq                                     \tilde{R}_{n}\paren{\cF_{\beta,d,n}} \geq 
	\frac{M}{4}n\tilde{\gamma} \geq 
	\frac{
		\cdot M^2}{16\norm{g}_{W_{\infty}^{\beta}\paren{\cX}}}n^{1-\frac{\beta}{d}},
	\end{align*}
	so that the case $ \frac{d}{p} < \beta\leq \frac{d}{2}$  is proved.
	
	To prove the second bound, notice that by Lemma~\eqref{lem:f_incl} for any $n \in \mbn_{*}$, $\cF_{\beta,d,n} \subset B_{W_{\infty}^{\beta}\paren{\cX}}\paren{0,\frac{M}{4}}$, which implies that $\fatt_{\gamma}\paren{B_{W_{\infty}^{\beta}}\paren{0,\frac{M}{4}}} \geq \fatt_{\gamma}\paren{\cF_{\beta,d,n}}$. In particular, this holds  if we choose $n_{0} \bydef \lfloor \paren{\frac{\gamma}{C_{M,g}}}^{-\frac{d}{\beta}  } \vee 1\rfloor $, then $n_{0} < \paren{\frac{\gamma}{C_{M,g}}}^{-\frac{d}{\beta}} \vee 1$, which is equivalent to $C_{M,g}n^{-\frac{\beta}{d}}_{0} \leq \gamma$. Notice that if $\gamma_1 < \gamma_2$ then $\fatt_{\gamma_1}\paren{\cF} \geq  \fatt_{\gamma_2}\paren{\cF}$. Applying the first property to the classes $\paren{B_{W_{\infty}^{\beta}}\paren{0,\frac{M}{4}}}$ and $\cF_{\beta,d,n_{0}}$ on the scale $\gamma$ and the second property for the class $\cF_{\beta,d,n_o}$ on the scales $\gamma$ and $C_{M,g}n_{0}^{-\frac{\beta}{d}}$ we consequently get
	\begin{align}
	\fatt_{\gamma}\paren{B_{W_{\infty}^{\beta}\paren{\cX}}\paren{0,\frac{M}{4}} } \geq  \fatt_{\gamma}\paren{\cF_{\beta,d,n_{0}}} \geq \fatt_{C_{M,g}n_{0}^{-\frac{\beta}{d}}}\paren{\cF_{\beta,d,n_{0}}} \geq n_{0}. 
	\end{align}
	
	Finally, because $n_{0} \geq 1$, so by using elementary $\lfloor a \rfloor \geq \frac{a}{2}$ we have $n_{0} \geq \frac{1}{2}\paren{\frac{\gamma}{C_{M,g}}^{-\frac{d}{\beta}} \vee 1 }$; therefore,
	$$\fatt_{\gamma}\paren[1]{B_{W_{\infty}^{\beta}}\paren{0,M/4} } \geq \fatt_{\gamma}\paren{ \cF_{\beta,d,n_{0}}} \geq \frac{1}{2} \paren[1]{\paren{ \frac{\gamma}{C_{M,g}}}^{-\frac{d}{\beta}} \vee 1}$$. 
	
	Choose $\gamma \bydef C_{M,g}^{\frac{d}{\beta +d}}n^{-\frac{\beta}{2\beta+d}}$, $n_{0} \bydef \lfloor \paren{\frac{\gamma}{C_{M,g}}}^{-\frac{d}{\beta}} \rfloor $, where $C_1$ is a constant as in Lemma~\ref{lem:rad_up_bound} and $C_{M,g}$ is a constant as in Equation~\eqref{eq:expres_shat}. For $\beta > \frac{d}{2}$, we have by inclusion and by Lemma~\ref{lem:rad_up_bound} that for any $n$ with the choice of $\gamma$ as before it holds: $ \fatt_{\gamma}\paren{\cF_{\beta,d,n}} \leq  \fatt_{\gamma}\paren[1]{B_{W_{\infty}^{\beta}\paren{\cX}}\paren{0,\frac{M}{4}}} \leq \tilde{C}n^{\frac{2\beta}{2\beta +d}} < \tilde{C}n$. Furthermore, as for any $n \in \mbn$, $B_{W^{\beta}_{\infty}\paren{\cX}  }\paren{0,\frac{M}{4}} \supset \cF_{\beta,d,n}$, so, in particular, $B_{W^{\beta}_{\infty}\paren{\cX}}\paren{0,\frac{M}{4}} \supset F_{\beta,d,n_0}$, which implies $R_{n}\paren{B_{W^{\beta}_{\infty}\paren{\cX}}\paren{0,\frac{M}{4}}} \geq R_{n} \paren[1]{\cF_{\beta,d,n_0}} $. Thus, applying Lemma \ref{lem:RS_15} to the class $\cF_{\beta,d,n_{0}}$ with any $n$ and $\gamma,n_{0}$ as above, we obtain:
	
	\begin{align*}
    	\tilde{R}_{n}\paren{\cF_{\beta,d,n_{0}}} &\geq C \gamma\sqrt{n}\paren[1]{                2\sqrt{2}\sqrt{\fatt_{\gamma}\paren{\cF_{\beta,d,n_{0}}}} - \sqrt{n}\gamma} \\
	    & \geq C n^{-\frac{\beta}{2\beta+d}} \sqrt{n} \paren[1]{2 C_{M,g}^{\frac{d}{d+\beta}} n^{\frac{d}{2\paren{2\beta+d}} } - C_{M,g}^{\frac{d}{\beta +d}} n^{\frac{1}{2} - \frac{\beta}{2\beta+d}}} \\
    	& \geq \tilde{C_{1}} n^{-\frac{\beta}{2\beta+d} + \frac{1}{2} + \frac{d}{2\paren{2\beta +d}}}
    	= \tilde{C_1}{n^{\frac{d}{2\beta +d}}},
	\end{align*}
	where $C_{1}$ is some constant independent of $n$.
	Now, the final bound for $\beta<\frac{d}{2}$ follows from the inclusion $\cF_{\beta,d,n_{0}} \subset B_{W_{\infty}^{\beta} }\paren{0,\frac{M}{4}} \subset B_{W_{p}^{\beta}\paren{\cX}}\paren{0,\frac{M}{4}}$ which implies
	$\tilde{R}_{n}\paren[1]{B_{W_{\infty}^{\beta}\paren{\cX}}\paren{0,M} } \geq \tilde{R}_{n}\paren[1]{B_{W_{p}^{\beta}\paren{\cX}}\paren{0,M} } \geq \tilde{R}_{n}\paren[1]{\cF_{\beta,d,n_{0}}}$.
\end{proof}

\section{Regret rates comparison}
\label{app:regret_rates_comp}
\begin{figure}[htp]
\centering
\includegraphics[width=.3\textwidth]{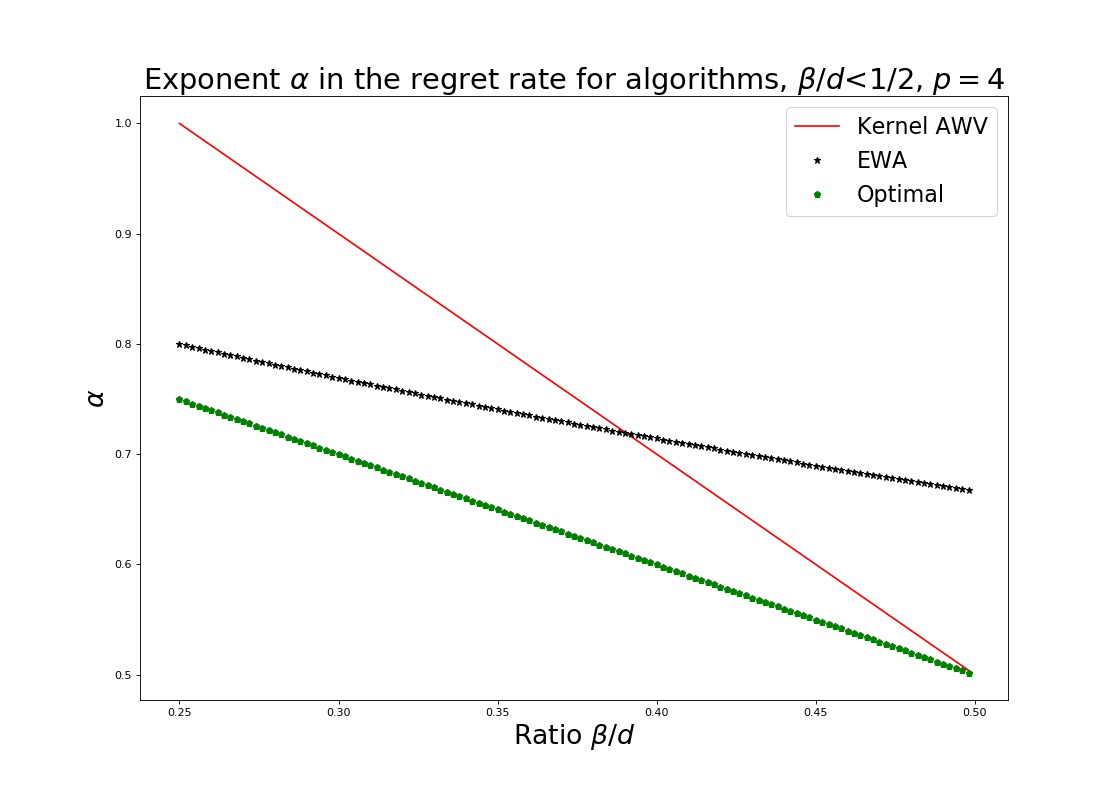}\hfill
\includegraphics[width=.3\textwidth]{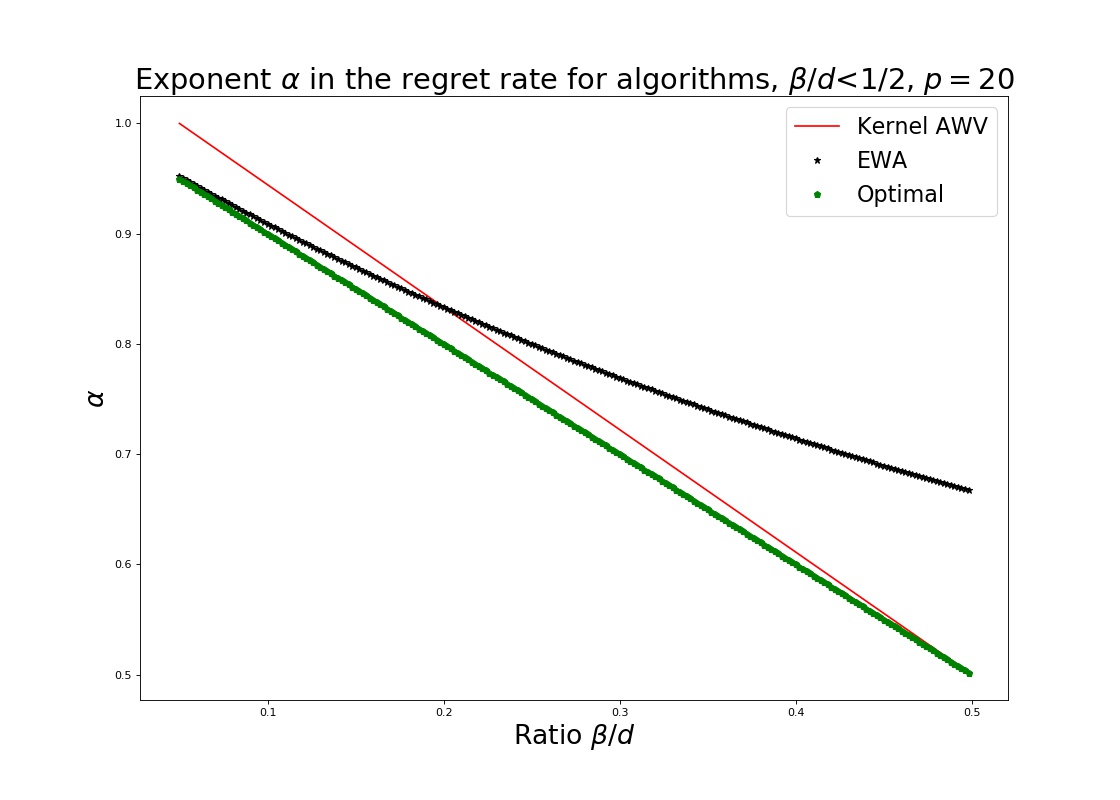}\hfill
\includegraphics[width=.3\textwidth]{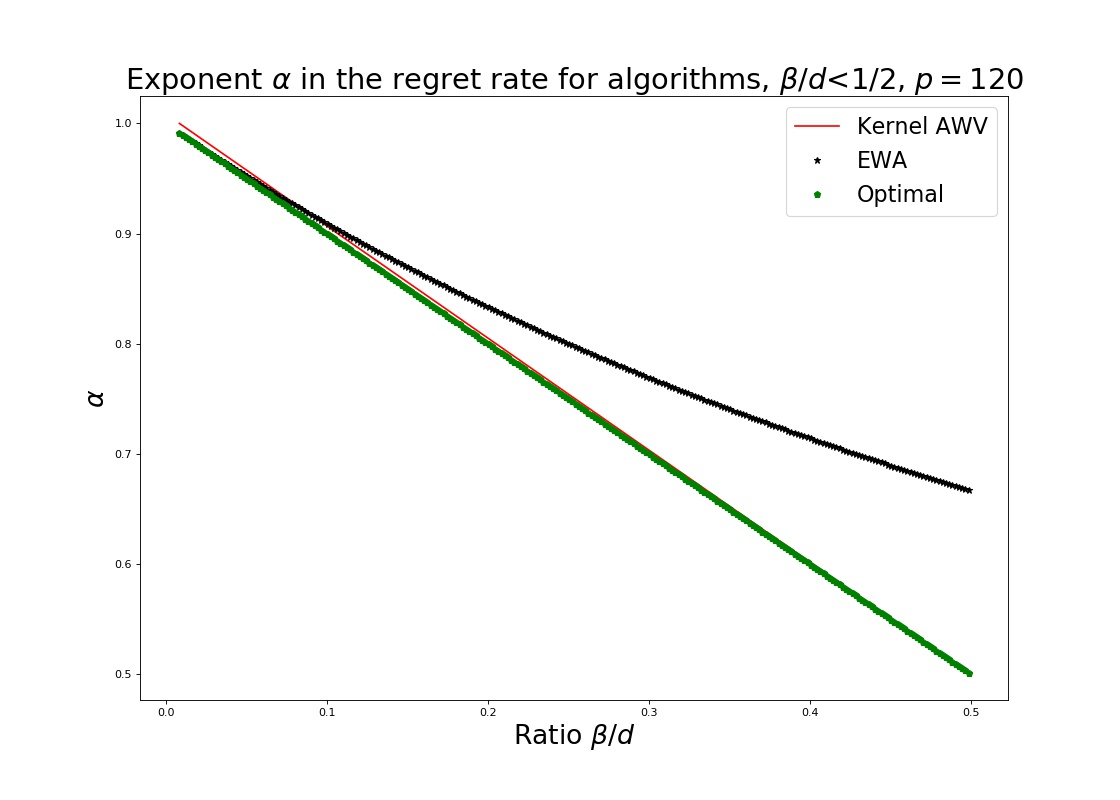}

\caption{Exponent of the regret in the case $W_{p}^{\beta}\paren{\cX}$, $\frac{1}{p}< \frac{\beta}{d}\leq \frac{1}{2}$, $p=4,20,120$.}
\label{fig:figure_pneq2i}

\end{figure}

Here we provide a short comparison of the exponents of theoretical regret rates between \algo~\eqref{eq:kern_ridge} and EWA  (\cite{Vovk06Metric}). 
One can check that when $\frac{\beta}{d}<\frac{\sqrt{1+4p}-1}{2p}$, EWA provides better rate then \algo given by~\eqref{eq:kern_ridge} with $s=\frac{d}{2}+\epsilon$, $\epsilon >0$ and $\tau_{n}$ chosen as in the Theorem~\ref{thm:higher_beta_p}. For a fixed pair $(\beta,d)$ this means that with increasing regularity of the function $f$ in terms of its integral $p-$norm, \algo estimates its behaviour better then EWA for larger range of possible values $(\beta,d)$.  This effect is illustrated in Figure~\ref{fig:figure_pneq2i}. 

\end{document}